\documentclass[12pt]{article}

\usepackage[scale=0.75]{geometry}
\usepackage{times}
\usepackage[utf8x]{inputenc}
\usepackage[T1]{fontenc} 
\usepackage{color}
\usepackage{soul}

\usepackage{enumitem}
\usepackage{interval}
\usepackage{graphicx}
\usepackage{amsmath,amsfonts,amssymb}
\usepackage{mathrsfs} 

\definecolor{vertfonce}{rgb}{0., 0.66, 0.25}
\definecolor{rougefonce}{rgb}{0.75, 0.09, 0.}
\definecolor{chocolat}{rgb}{0.823,0.412,0.118}
\definecolor{purp}{rgb}{0.66,0.,0.99}
\definecolor{cyanf}{rgb}{0.,0.5,0.99}

\usepackage[breaklinks=true,
colorlinks=true,
linkcolor=rougefonce,
citecolor=vertfonce]%
{hyperref}

\newcommand{\plus}[1]{\mathop{\,+\,}\limits_{#1}}
\newcommand{\Lh}{\mathcal{L}_\h}

\newcommand{\cI}{\mathcal{I}}

\newcommand{\restr}{\raisebox{-0.5ex}{$\upharpoonright$}}
\newcommand{\trsp}{\raisebox{.6ex}{${\scriptstyle t}$}}
\newcommand{\phy}{\varphi}

\def\hequa{\renewcommand{\theequation}{$\h$-\arabic{equation}}}
\def\nequa{\renewcommand{\theequation}{\arabic{equation}}}

\newcommand{\TM}{\mathbb{T}}
\newcommand{\RM}{\mathbb{R}}
\newcommand{\ZM}{\mathbb{Z}}
\newcommand{\QM}{\mathbb{Q}}
\newcommand{\NM}{\mathbb{N}}

\newcommand{\SM}{\mathbb{S}}
\newcommand{\PM}{\mathbb{P}}

\newcommand{\RP}{\mathbb{RP}}

\newcommand{\Cinf}{C^\infty}
\newcommand{\ham}[1]{\mathcal{X}_{#1}}
\newcommand{\deriv}[2]{\frac{\partial #1}{\partial #2}}
\newcommand{\h}{\hbar}
\newcommand{\ku}{{k}}
\newcommand{\eu}{{e}}

\newcommand{\abs}[1]{\left|#1\right|}
\newcommand{\norm}[1]{\left\|#1\right\|}

\newcommand{\theor}{Theorem}
\newcommand{\defin}{Definition}
\newcommand{\lemma}{Lemma}
\newcommand{\remar}{Remark}
\newcommand{\corol}{Corollary}
\newcommand{\propo}{Proposition}
\newcommand{\notat}{Notation}
\newcommand{\demon}{Proof}

\newtheorem{theo}{\theor}[section]
\newtheorem{defi}[theo]{\defin}
\newtheorem{prop}[theo]{\propo}
\newtheorem{lemm}[theo]{\lemma}
\newtheorem{coro}[theo]{\corol}
\newtheorem{nota}[theo]{\notat}

\newcommand{\demons}[1][$\!\!$]{\noindent\textbf{\demon\ }\textsl{#1}\textbf{.}~}
\newcommand{\cqfd}{\hfill $\square$\par\vspace{1ex}}
\newcommand{\dd}[1]{\ensuremath{\operatorname{d}\!{#1}}}
\renewcommand{\O}{\mathscr{O}}

\newenvironment{rema}
{\par\vspace{1ex}\refstepcounter{theo}%
\noindent\textbf{\remar~\thetheo} }
{~\hfill\mbox{$\triangle$}\par\vspace{1ex}}

\newenvironment{demo}[1][$\!\!$]
{\demons[#1]\ }
{\cqfd}

\title{Asymptotic Lattices, Good Labellings,\\ and the Rotation Number
  for Quantum Integrable Systems}

\author{Monique \textsc{Dauge}\footnote{Univ Rennes, CNRS, IRMAR - UMR
    6625, F-35000 Rennes, France} \and Michael
  \textsc{Hall}\footnote{University of Southern California, Los
    Angeles (USA)} \and \textsc{Vũ Ngọc} 
  San$\,\!^*$}

\date{}

\begin{document}
\maketitle

\begin{abstract}
  This article introduces the notion of good labellings for asymptotic
  lattices in order to study joint spectra of quantum integrable
  systems from the point of view of inverse spectral theory. As an
  application, we consider a new spectral quantity for a quantum
  integrable system, the quantum rotation number. In the case of two
  degrees of freedom, we obtain a constructive algorithm for the
  detection of appropriate labellings for joint eigenvalues, which we
  use to prove that, in the semiclassical limit, the quantum rotation
  number can be calculated on a joint spectrum in a robust way, and
  converges to the well-known classical rotation number.  The general
  results are applied to the semitoric case where formulas become
  particularly natural.
  \end{abstract}

\begin{footnotesize}
  \noindent \textbf{MS Classification :}
  81S10, 
  81Q20, 
  58J40, 
  58J50, 
  65L09. 
\end{footnotesize}

\tableofcontents

\section{Introduction}
\label{s:intro}
\subsection{Motivations and aims}
Our motivations mainly come from the connection between integrable
two-degree of freedom Hamiltonian systems and the spectrum of related
quantum systems.  Let $M$ be $\RM^4$ or, more generally, a
4-dimensional symplectic manifold.  If a two-degree of freedom
Hamiltonian system, given by a Hamiltonian $H\in\Cinf(M)$, is
integrable in the classical Liouville sense, then it is well known
that most of the dynamics takes place on invariant tori of dimension
2. On each such torus $\Lambda$, the motion is particularly simple: in
suitable angle coordinates $(\theta_1,\theta_2)$ on $\Lambda$, the
Hamiltonian vector field $\ham{H}$ is constant:
\[
  (\ham{H}){\restr_\Lambda} = \alpha_1 \deriv{}{\theta_1} + \alpha_2
  \deriv{}{\theta_2}.
\]
This is the content of the classical action-angle theorem (see for
instance~\cite{duistermaat} and the references therein). The direction
in $\RP^1$ given by the frequency vector $w:=(\alpha_1,\alpha_2)$ is
called the \emph{rotation number}. If $w$ is rational (by this we mean
that $\alpha_1$ and $\alpha_2$ are linearly dependent over $\QM$),
then the trajectory of the Hamiltonian system on $\Lambda$ is
periodic. On the contrary, when $w$ is irrational, this trajectory is
dense on $\Lambda$. Thus, the knowledge of $w$ gives important
information on the nature of the dynamics. Understanding the variation
of $w$ is also crucial for the study of perturbations of $H$, via the
various ``KAM'' theorems~(see for instance the review
article~\cite{poeschel-lecture-kam}).

The first goal of this article is to investigate the effect of this
dynamical quantity on a \emph{quantum} system. Assume for instance
that we consider a Schrödinger operator $P=-\h^2\Delta + V$, with a
smooth potential $V$ on a 2-dimensional Riemannian manifold $X$;
assume moreover that this operator is \emph{quantum
  integrable} 
in the sense that there exists another selfadjoint operator $\hat J$,
commuting with $P$, whose principal symbol $J$ is an independent
integral of motion for the Hamiltonian $H(x,\xi)= \norm{\xi}^2+V(x)$.
Therefore $H$ is completely integrable, and we can ask the question:
what is the manifestation of the underlying classical rotation number
$w$ on the spectrum of $P$? It has been acknowledged for a long time
that so-called non-degeneracy hypothesis, or diophantine hypothesis on
the rotation number, à la KAM, can be crucial to obtain a good
description of the spectrum of quantum integrable systems, see for
instance~\cite{bleher-revolution, bleher}. But, is there a direct, or
recognizable, signature of the classical rotation number on the
spectrum?  A first answer was given by Hitrik and Sjöstrand in a
series of papers~\cite{hitrik-sjostrand-I,san-hitrik-sjostrand} where
they study the case of weakly non-selfadjoint operators (\emph{i.e.},
$V$ is perturbed by a small imaginary term), and they proved that the
asymptotics of the spectrum, in the semiclassical limit $\h\to 0$,
exhibit very different behaviors depending on the rationality of the
rotation number $w$.

In the purely selfadjoint case, the construction of quasi-modes in the
various situations where $w$ is strongly irrational or not was known
for a long time, see for instance~\cite{lazutkin-book,
  colin-quasi-modes}; however, perhaps surprisingly, only the
non-selfadjoint case, where the spectrum, instead of being
one-dimensional, is deployed in the complex plane, gives some hope to
recognize useful geometric structures from the eigenvalues
themselves. For quasi-periodic dynamics, this idea was exploited to
recover the Birkhoff normal form from the complex spectrum
in~\cite{hall-13}, and to define quantum monodromy in the
non-selfadjoint case~\cite{phan-14}.  A strong motivation for our work
is the recent paper~\cite{hitrik-sjostrand-band}, which naturally
leads to this intriguing question: can you detect the rationality of
the rotation number from the spectrum?

It is precisely this type of inverse spectral problem that we study in
this paper. Here we stick to the simpler ``normal'' case, which means
that instead of considering the spectrum of a truly non-selfadjoint
operator, we consider the joint spectrum of a pair of commuting
selfadjoint operators. The fully non-selfadjoint case is still largely
open. We prove that one can define a \emph{quantum rotation number},
in a very natural and concrete way, from the joint spectrum of two
commuting operators as soon as we know appropriate labelling of the
joint spectrum by pairs of integers, given for instance by the
Bohr-Sommerfeld quantization rule that provides us ``quantum
numbers''; moreover we show that, in a suitable sense, this quantum
rotation number converges to the classical rotation number in the
semiclassical limit $\h\to 0$. This result, which we actually prove in
any dimension, is however more delicate to use than one could think at
first sight if one does not have an a priori knowledge of appropriate
labellings of the joint spectra as $\h\to0$.  A fully relevant
question is the possibility to construct such labellings with suitable
regularity properties from the bare data of joint spectra.


Precisely because of this issue, the second goal of this article is to
set a rigorous ground for the theory of labellings of joint
spectra. Given a discrete set of points in $\RM^2$ (obtained as the
joint spectrum of some unknown system), how can you label these points
is a sensible manner, in a similar way as the usual Bohr-Sommerfeld
rules would tell you to do so? This turns out to be a delicate
question, even for systems with global action-angle variables. Such
discrete sets are called asymptotics lattices; in this article, the
theory is build up from scratch up to the point where it is finally
able to provide the proof of the inverse spectral problem for the
quantum rotation number. However, we expect further developments
and generalizations in the future, for instance in relation to trace
formulas.

After going through this constructive and combinatorial issue, we
finally show that one can detect the classical rotation number from
the joint spectrum of an associated quantum integrable system. We also
believe that the quantum rotation number will prove to be a useful
object in the study of quantum integrable (or near-integrable)
systems, and we hope to apply this idea on concrete systems in a near
future. The recent article~\cite{hamraoui2018} on asymmetric-top
molecules nicely supports this idea.

The second motivation of this paper comes from the theory of semitoric
systems, see~\cite{san-polytope,san-alvaro-I, san-alvaro-II,
  san-alvaro-first-steps}. The general conjecture for quantum
semitoric systems is that one can always recover the underlying
classical system from the spectrum of a quantum semitoric
system~\cite[Conjecture
9.1]{san-alvaro-survey}. In~\cite{san-alvaro-first-steps}, a sketch of
proof was provided, which builds on a number of geometric and spectral
invariants, like quantum monodromy~\cite{cushman-duist,
  san-mono}. When we wrote this paper, the inverse result was known in
the toric case~\cite{san-charles-pelayo}, and recent advances showed
that a general conjecture like this was not out of
reach~\cite{san-lefloch-pelayo:jc,san-alvaro-inverse-focus,lefloch-pelayo19},
although probably under suitable genericity assumptions. Thus, we were
motivated to exhibit a new invariant that can be recovered from the
spectrum. It now turns out that our work, and in particular the
precise study of asymptotics lattice and their labellings
(Sections~\ref{sec:asympt-latt} and \ref{sec:labelling-algorithm}),
has been used crucially to finally completely solve the semitoric
conjecture, in a constructive way~\cite{san-yohann21}.

While the rotation number can be defined for any completely integrable
system, it bears a particular nice form in the case of semitoric
systems, see Section~\ref{sec:semitoric-case}.  We believe that our
semitoric theory and algorithms should find natural applications in
the study of axisymmetric Schrödinger operators, a question that we
hope to investigate in a future work.

Our work contributes to the general inverse spectral theory in the
semiclassical limit. Of course, in such generality, this question has
a long history, especially when restricted to the particular case of
Laplace-Beltrami or Schrödinger operators; see for instance the
survey~\cite{datchev-hezari-survey13} and the references therein. The
semiclassical inverse problem for more general Hamiltonians, which is
our concern here, was also considered, albeit by less numerous
studies, see in particular~\cite{iantchenko-sjostrand-zworski} for a
very general treatment. In the case of Liouville integrable systems, a
line of program, closely following recent development in the
classification of Lagrangian fibrations, was proposed
in~\cite{san-inverse, san-alvaro-first-steps, san-sepe-poisson16}, but
prior works have been produced in the Riemannian case, see for
instance~\cite{zelditch-revolution}.

\subsection{Outline and sketch of results}
Our paper consists of two main sections, \ref{sec:classical} and
\ref{sec:quantum}, devoted to classical and quantum rotation numbers,
respectively.  Correspondences between the two are proved in
Section~\ref{sec:quantum}.

In Section~\ref{sec:classical} we recall and generalize, in the
setting of completely integrable Hamiltonians $H$ on $2n$-dimensional
symplectic manifolds, the classical definition of the rotation number:
for such $H$, associated with its constants of motions
$(f_1,\ldots,f_{n-1})$, there exist action-angle coordinates
$(I,\theta)\in\ T^*\TM^n$, in which $H=g(I_1,\ldots,I_n)$ and a local
definition of the rotation number is, for each chosen Liouville
torus $\Lambda$
\[
  [w_I](\Lambda) = [\partial_1 g(I(\Lambda)) : \dots : \partial_n g(I(\Lambda))] \in \RP^{n-1}.
\] 
We explain under which conditions it is well defined (which seems to
be an information that is difficult to locate in the literature),
proving in particular that the rotation number $[w_{I'}]$ associated
with another set $I'$ of action variables is deduced from $[w_I]$ by a
$\textup{SL}(n,\ZM)$ transformation (Lemma \ref{lemm:ItoI'}).  We also
exhibit its relationship with Hamiltonian monodromy
(Proposition~\ref{prop:monodromy}).  Then we restrict to the case
$n=2$ to address the semitoric case in which the rotation number is
well-defined modulo integers (Proposition~\ref{prop:semitoric}).

In Section~\ref{sec:quantum}, we first recall known facts about the
joint spectrum of semi-classical commuting $\h$-pseudodifferential
operators $(P_1,\ldots,P_n)$, arriving to the fundamental result that
the joint spectrum $\Sigma_\h$ is (locally) a deformed lattice modulo
$\O(\h^\infty)$:
\[
  G_\h(\h k_1,\ldots, \h k_n) + \O(\h^\infty),\quad (k_1,\ldots, k_n)\in\ZM^n
\]
in which $G_\h$ is a smooth map admitting an asymptotic expansion in
integer powers of $h$ (Theorem \ref{theo:bs}). These Bohr-Sommerfeld
quantization rules tell two things: first, the spectrum as a set
resembles a lattice, and second, its elements can be numbered (or
``labelled") by multi-integers $(k_1,\ldots, k_n)$ in a coherent way
as $\h\to0$. This is what we call a ``good labelling''.

In the special case when $P_n=\hat H$ has the Hamiltonian $H$ as
principal symbol, and $J_j$ is the principal symbol of $P_j$, so that
the map $(J_1,\ldots,J_{n-1},H)$ form a $n$-degree of freedom
completely integrable system $F$, having at hands the joint spectra
$\Sigma_\h$ and their good labellings
\[
  \ZM^n \ni (k_1,\ldots, k_n)=:\ku\mapsto\lambda_{\ku}(\h) = 
  (\lambda^{(1)}_{\ku}(\h),\ldots,\lambda^{(n)}_{\ku}(\h))\in\Sigma_\h\,,
\]
we define quantum rotation numbers as projectivized finite differences
of the eigenvalues $E_{\ku}(\h)=\lambda^{(n)}_{\ku}(\h)$ of $\hat H$
\[
  [\hat w_\h](\ku) = [E_{\ku+\eu_1}(\h)-E_{\ku}(\h): \ldots : 
  E_{\ku+\eu_n}(\h)-E_{\ku}(\h)] \in \RP^{n-1},
\]
with $\eu_1,\ldots,\eu_n$ unit coordinate vectors in $\ZM^n$. Then we
prove (Theorem~\ref{theo:q-rotation-number}) that the quantum rotation
number is a semiclassical deformation of the classical rotation
number, namely: if $c$ is a regular value of $F$, there exist action
variables $I$ so that 
\[
  [\hat w_\h](\ku_\h) = [w_I](\Lambda) + \O(\h),\quad \h\to0,
\]
where the labels $\ku_\h$ are such that
$\lambda_{\ku_\h}(\h)=c+\O(\h)$ and $\Lambda$ is the Liouville torus
$F^{-1}(c)$.

Though satisfying, this result is not the end of the story, but rather
its starting point.  The question that we address in the rest of the
paper is, from the sole knowledge of the joint spectra $\Sigma_\h$ in
some window $\Omega$ for some family of values of $\h$ tending to $0$,
can we construct labellings of these sets so that a quantum rotation
number can be computed?  Before introducing an algorithm, we analyze
in more detail the properties of different sorts of labellings of an
``asymptotic lattice" (a family of sets $\Lh$ displaying the same
structure as the joint spectra of commuting $\h$-pseudodifferential
operators). We prove that good labellings are not unique, but can be
deduced from each other by transformations in the affine group
$\textup{SL}(n,\ZM)\ltimes\ZM^n$ (Proposition \ref{prop:A}). We then
face the issue that, in general, good labellings are not
algorithmically accessible. Indeed, the asymptotic lattice does not
necessarily concentrate to a point in the observation window, but
rather goes through. That is why we introduce a relaxed version of
labelling, the \emph{linear labellings}, which coincide with a good
labelling up to a translation $\varkappa_\h\in\ZM^n$ (depending on $\h$
in general). Fortunately, linear labellings retain enough information
to allow for the determination of a rotation number.

The core of the paper is to obtain a robust and constructive
determination of linear labellings for arbitrary 2-dimensional
asymptotic lattices.  To achieve this goal, we propose a two part
algorithm. In the first part of the algorithm, for a chosen $\h$,
small enough, we find a labelling of $\Lh$ by a set of bi-integers
$(m,n)\in\ZM^2$ containing $\ZM^2\cap B(0,\rho\h^{-1})$ for a
$\rho>0$. The label $(0,0)$ is associated with a closest element
$\lambda\in\Lh$ to a distinguished point $c\in\Omega$. The conclusive
statement of this part is Theorem \ref{theo:basis} which, despite its
apparent simplicity, necessitates an intricate proof. The most
delicate part is to prove that all points of $\Lh$ (in the window
$\Omega$) are effectively labelled, and for this we use the dimension
$n=2$ in a crucial way. In the second part of the algorithm, we
combine this construction for any $\h$ belonging to a sequence tending
to $0$, with appropriate $\textup{SL}(n,\ZM)$ transformations, in
order to finally obtain a linear labelling (Theorem
\ref{theo:algo}). The determination of a quantum rotation number from
the sole knowledge of a joint spectrum, and its convergence to a
classical rotation number come as natural consequences of this
algorithm (Theorem \ref{theo:A}).

In the semitoric case, the Hamiltonian $J$ induces a periodic
flow. This has remarkable consequences both on (quantum) rotation
numbers and joint spectra. On the one hand, the classical rotation
number can be identified with an angle with a natural dynamical
interpretation (Proposition~\ref{prop:semitoric}). On the other hand,
the joint spectrum locally displays a nice structure of vertical bands
(Proposition~\ref{prop:J-spectrum}); this makes the process of
detecting labellings easier. We show in
Section~\ref{sec:semitoric-algorithm} how to correct the first part of
the general algorithm in order to obtain a labelling that respects the
semitoric structure. As a consequence, we obtain that the semitoric
version of the rotation number can be computed from the joint spectrum
(Corollaries~\ref{coro:semitoric-elliptic} and~\ref{coro:semitoric}).

\paragraph{Acknowledgements}
We would like to thank the anonymous referees for their insights,
which have substantially improved the presentation and clarified the
geometric background.  The authors are happy to acknowledge the
support of the Centre Henri Lebesgue, program ANR-11-LABX-0020-0.
M.H. benefited from ANR NOSEVOL (ANR 2011 BS01019 01) during the start
of this project.

\section{The Classical Rotation Number}
\label{sec:classical}

In this section, we first recall the notion of rotation number for
integrable Hamiltonian systems of any dimension. A general discussion
for $n=2$ can also be found in~\cite{bolsinov-fomenko-book}. Then, we
show the relationship between the globalized version of rotation
numbers with Hamiltonian monodromy
(Section~\ref{sec:monodromy}). Finally, we explain how the usual
rotation number (which is an angle) fits within the framework of
semitoric systems when $n=2$ (Section~\ref{sec:semitoric-case}).

Let $H$ be a completely integrable Hamiltonian on a symplectic
manifold $M$ of dimension $2n$: there exists smooth functions
$f_1,\dots, f_{n-1}$ on $M$ such that, letting $f_n:=H$, then for all
$j$ and $k$, $\{f_j, f_k\}=0$ and the differentials $df_1,\dots, df_n$
are almost everywhere independent. From the dynamical viewpoint, the
Hamiltonian $H$ defines a dynamical system through its Hamiltonian
vector field $\ham{H}$, and the functions $f_j$ are constants of
motion. On the geometric side, we have a foliation of the phase space
$M$ by (possibly singular) Lagrangian leaves given by the common level
sets of the map $F=(f_1,\dots, f_n)$.

\subsection{Local rotation number}

Let us assume that the map $F:=(f_1,\dots, f_n):M\to \RM^n$ is proper,
and let $c$ be a regular value of $F$. By the action-angle
theorem~\cite{mineur-action-angle} (see~\cite{duistermaat} for a more
modern proof), the level set $F^{-1}(c)$ is a finite union of
$n$-dimensional tori, called \emph{Liouville tori}; and near each
Liouville torus, there exist action-angle coordinates
$(I,\theta) := (I_1,\dots,I_n,\theta_1,\dots,\theta_n)\in
\textup{neigh}(\{0\}\times \TM^n, \RM^n \times \TM^n)$
such that
\begin{equation}
  H=g(I_1,\dots,I_n),
  \label{equ:H}
\end{equation}
for some smooth map $g:\RM^n\to\RM$, and the symplectic form of $M$ is
$\dd I_1\wedge d\theta_1 + \cdots + \dd I_n \wedge d\theta_n$. More
precisely, the action-angle theorem states that there is a local
diffeomorphism:
\[
  G: (\RM^n, 0) \to (\RM^n, c)
\]
such that, in the new coordinates, $F=G(I)$. We shall always assume
(as we may) that the actions $I$ are \emph{oriented} with respect to
$F$, which means that $\det d G(0)>0$. In particular, $d g$ does not
vanish near the origin, which enables the following definition:
\begin{defi}
  \label{defi:rotation}
  The \emph{rotation number} of $H$ relative to the oriented action
  variables $I:=(I_1,\dots,I_n)$ on the Liouville torus
  $\Lambda\subset M$ is the projective vector
  \[ [w_{I}](\Lambda) := [\partial_1 g(I(\Lambda)) : \dots : \partial_n
    g(I(\Lambda))] \in \RP^{n-1}.
  \]
\end{defi}
Probably the easiest dynamical interpretation of
$[w_{I}](\Lambda)$ is the following. One deduces
from~\eqref{equ:H} that the Hamiltonian vector field of $H$, which is
tangent to $\Lambda$, has the form
\[
  \ham{H} = \partial_1 g\deriv{}{\theta_1} + \cdots + \partial_n g
  \deriv{}{\theta_n}.
\]
Therefore, the flow of $\ham{H}$ is a ``straight line'' winding
quasi-periodically on the affine torus $\Lambda$ in the coordinates
$(\theta_1,\dots,\theta_n)$, and the rotation number $[\nabla g]$ is
simply the direction of the trajectory. The map
$\Lambda\mapsto \nabla g(\Lambda) := (\partial_1 g(I(\Lambda)),\dots,
\partial_n g(I(\Lambda)))$ is generally referred to as the
\emph{frequency map} of the system. The rotation number is nothing but
the projectivized version of the frequency map. 

When $n=2$, it is also
convenient (and usual) to define the rotation number as the ratio:
\begin{equation}
  w_I(\Lambda):=\frac{\partial_1 g(I(\Lambda))}{\partial_2
    g(I(\Lambda))} \in \overline{\RM}
  \label{equ:ratio}
\end{equation}
with values in the 1-point compactification of the real line
$\overline{\RM}=\RM \cup \{\infty\}$.  The diffeomorphism
$\overline{\RM} \to \RP^1, z \mapsto [z]$ is defined by
$[z] := [z : 1 ]$ if $z\in\RM$ and $[z] := [1 : 0]$ if $z=\infty$, and
$[w_I](\Lambda) = [w_I(\Lambda)]$.

Note that the rotation number is not well defined by $H$ only: it
depends on the choice of actions, which in turn might depend on the
choice of constants of motion $f_j$. In this paper, we shall always
assume that $f_1,\dots,f_{n-1}$ are given, as part of the data of the
integrable Hamiltonian $H$; the question of the consequences of the
choice of $f_j$ is interesting and apparently not widely spread; we
hope to return to this problem in a future paper.

Thus, we assume that $F=(f_1,\dots,f_{n-1},, f_n=H)$ is fixed.  We
formalize the result on the change of action variables $I\to I'$ in
the next lemma for further reference.
\begin{lemm}
  \label{lemm:ItoI'}
  If $I'=(I'_1,\dots,I'_n)$ is another set of action variables, then
  there is a matrix $A \in \textup{SL}(n,\ZM)$ such that
  \begin{equation}
    \dd I:=A\, \dd I',
    \label{equ:changement-actiond}
  \end{equation}
  and the new rotation number related to $I'$ is
  \begin{equation}
    [w_{I'}]=\trsp A \circ  [w _I] \;,
    \label{equ:changement-action}
  \end{equation}
  where $\circ$ denotes the natural action of $\textup{SL}(n,\ZM)$ on
  $\RP^{n-1}$ (elements of $\RP^{n-1}$ are viewed here as equivalence
  classes of vectors in $\RM^n$). When $n=2$, if instead we write
  $w_I$ as the quotient~\eqref{equ:ratio}, then $w_{I'}$ is obtained
  by the Möbius transformation:
  \begin{equation}
    w_{I'} = \frac{a w_I + c}{b w_I + d}\quad\mbox{in which}\quad
    A = \begin{pmatrix}
      a & b \\ c & d
    \end{pmatrix}.
    \label{equ:moebius}
  \end{equation}
\end{lemm}

An important fact is that the rank of the components of $w_I$ over the
\emph{rationals} (when $n=2$, this is simply the rationality of $w_I$)
is well defined (i.e. preserved by $\textup{SL}(n,\ZM)$
transformations). This is easy to understand from the dynamical
viewpoint; when $n=2$, a rational rotation number is equivalent to a
periodic Hamiltonian flow for the vector field $\ham{H}$ on the
Liouville torus, which is of course independent of the choice of
action-angle coordinates. In higher dimensions, the rank over $\QM$
gives the dimension of the sub-torus where the dynamics takes place.



\begin{nota}
\label{notaBr}
  Let $\textup{B}_r$ be the set of regular Liouville tori. 
\end{nota}
If $F$ has
connected level sets, $\textup{B}_r$ can be identified with the open
subset in $\RM^n$ of regular values of $F$. In general, it follows
from the action-angle theorem that $\textup{B}_r$ is a smooth covering
above the open set of regular values of $F$ in $\RM^2$. In order to
make our inverse statement precise, we introduce the following
definition.

\begin{defi}
  \label{defi:arotation}
  A function $[w]:W\to \RP^{n-1}$ is called a \emph{rotation number}
  for the system $F=(f_1,\dots,f_{n-1},, f_n=H)$ on the open set
  $W\subset \textup{B}_r$ if for every $\Lambda\in W$, there exist a
  neighborhood $U$ of $\Lambda$ in $W$ and a set of action variables
  $I:=(I_1,\dots,I_n)$ on $U$ such that $[w] = [w_I]$ on $U$.
\end{defi}
(When $n=2$, we will use the same terminology for a function
$w:W\to\overline{\RM}$).  As a consequence, a rotation number is
always a smooth function.

\begin{rema}
  Let $E\in\RM$ be a regular value
  of the Hamiltonian $H$. The submanifold $\Sigma_E:=H^{-1}(E)$ is
  foliated by the level sets of $f:=(f_1,\dots,f_{n-1})$. Regular
  level sets of $f{\restr_{\Sigma_E}}$ correspond to Liouville tori;
  when $n=2$, they form smooth one-dimensional families inside
  $\Sigma_E$. The restriction of the rotation number to this family is
  usually called the \emph{rotation function}. It is important in many
  situations to know whether this function is a local
  diffeomorphism. This property (the so-called ``isoenergetic KAM
  condition'', see for instance~\cite{arnold-KAM}) is invariant under
  Möbius transformations.
\end{rema}

\subsection{Global rotation number and monodromy}
\label{sec:monodromy}

We have seen that rotation numbers exist in a neighborhood of any
regular Liouville torus. A natural question arises, whether it is
possible to define a rotation number on the whole set of regular tori,
$\textup{B}_r$. In fact, since the rotation number depends on the
choice of local action coordinates, it is naturally related to the
so-called Hamiltonian monodromy of the system. We show in this section
the relation between these two objects, which explains how to define a
global rotation number.  We were not able to locate this statement in
the literature, however it is implicitly underlying Cushman's argument
for the non-triviality of the monodromy of the Spherical Pendulum
reported in~\cite{duistermaat}; see also~\cite{cushman-book}. In the
presence of global $\mathbb{S}^1$ actions, this argument can be
generalized to the co-called ``fractional
monodromy''~\cite{broer-efstathiou-lukina10}.

In the realm of classical integrable systems (of any number of degrees
of freedom $n$), the Hamiltonian monodromy was discovered by
Duistermaat in~\cite{duistermaat}. It is the obstruction to the
existence of global action coordinates (modulo constants). In
topological terms, it is merely the $\textup{SL}(n,\ZM)$-holonomy of
the flat bundle over the set of regular tori of the system whose fiber
is the homology of the corresponding Liouville torus. Since then, a
number of references have explained this monodromy and its
relationships with the many interesting geometric data
involved. Nowadays, the prominent object which contains the
Duistermaat monodromy is the integral affine structure on the set of
regular tori. It shows up for instance in Mirror Symmetry questions,
see~\cite{kontsevich-soibelman}. The relationships with quantum
integrable systems was pointed out in~\cite{cushman-duist} and proved
in~\cite{san-mono}. See also~\cite{san-dijon} for an overview of
various aspects of monodromy and semiclassical analysis.

Recall that, by the action-angle theorem, given a Liouville torus
$\Lambda$, the oriented action variables $I=(I_1,\dots, I_n)$ give
local coordinates for $\textup{B}_r$, and hence the differentials
$\dd I_1(\Lambda),\dots,\dd I_n(\Lambda)$ span the cotangent space of
$\textup{B}_r$ at $\Lambda$, and Lemma~\ref{lemm:ItoI'} implies that
$T^*\textup{B}_r$ can be seen as a flat bundle over $\textup{B}_r$
with structure group $\textup{SL}(n,\ZM)$. This defines the integral
affine structure of $\textup{B}_r$. The Hamiltonian monodromy is the
holonomy of this flat bundle. 

For our main algorithm in Section~\ref{sec:labelling-algorithm}, it
will be useful to follow action variables when we move along
continuous paths in $\textup{B}_r$. The Hamiltonian monodromy gives a
homomorphism:
\[
  \mu : \pi_1(\textup{B}_r) \to \textup{SL}(n,\ZM),
\]
defined up to a global conjugacy.  Here $\pi_1(\textup{B}_r)$ is the
fundamental group of the set of regular tori $\textup{B}_r$. Since
$\textup{B}_r$ may not be simply connected, we introduce a simply
connected covering
\[
  \pi:\tilde{\textup{B}}_r \to \textup{B}_r.
\]
An element $\tilde\Lambda\in \tilde{\textup{B}}_r$ can be seen as a
homotopy class of a path from a fixed torus $\Lambda_0$ to
$\pi(\Lambda)$. Recall that a loop $\gamma$ in $\textup{B}_r$ acts
fiberwise on the covering $\tilde{\textup{B}}_r$ by concatenation of
paths: given $\tilde\Lambda\in \pi^{-1}(\gamma(0))$, we define
$\gamma\cdot\tilde\Lambda$ to be (the homotopy class of)
$\tilde\Lambda\circ \gamma$, \emph{i.e.} the path corresponding to
$\tilde \Lambda$ followed by $\gamma$.

Let us now turn to the global definition of rotation numbers.  Let $g$
be as in~\eqref{equ:H}, which we now view as a map on $\textup{B}_r$;
the rotation number $[w_I(\Lambda)]$ can now be identified with the
projectivized covector
\begin{equation}
  \label{equ:global_rotation_number_I}
  [w_I(\Lambda)] := [\dd{}_{\Lambda} g] \in \PM(T^*_\Lambda \textup{B}_r)\,.
\end{equation}
In view
of~\eqref{equ:changement-action}, the collection of rotation numbers
$[w_I]$, for all local action variables $I$, forms a (local) section
of the projectivized bundle $\mathbb{P}(T^*\textup{B}_r)$.  Because of
monodromy, a global section may not exist. Rather, a global rotation
number should be seen as a function on $\tilde{\textup{B}}_r$,
\emph{i.e.} a \emph{multivalued} function on $\textup{B}_r$.
\begin{prop}
  \label{prop:monodromy}
  Given a regular torus $\Lambda\in \textup{B}_r$ for the integrable
  system $F$, a choice of action variables $I$ near $\Lambda$, and
  $\tilde\Lambda\in\pi^{-1}(\Lambda)$, there exists a unique smooth function
  $[\tilde w]:\tilde{\textup{B}}_r\to \RP^{n-1}$ such that:
  \begin{enumerate}
  \item \label{prop:item:1} for any simply connected open set
    $U\subset \textup{B}_r$, the function $[w]:U\to \RP^{n-1}$ defined
    by
    \[
      [\tilde w] = [w]\circ \pi
    \]
    is a rotation number for the system $F$
    (Definition~\ref{defi:arotation}).
  \item \label{prop:item:2} near $\tilde\Lambda$,
    $[\tilde w]=[w_{I}]\circ \pi$.
  \end{enumerate}
  Moreover, for any loop $\gamma\in \pi_1(W)$,
  \begin{equation}
    \gamma \cdot [\tilde w] = \trsp \mu(\gamma) \circ [\tilde w],
    \label{equ:rotation_number_holonomy}
\end{equation}
  where
  $\gamma \cdot [\tilde w] := \left(\tilde\Lambda' \mapsto [\tilde w]
    (\gamma\cdot\tilde\Lambda')\right)$.
\end{prop}
\begin{demo}
  The section of $\mathbb{P}(T^*\textup{B}_r)$ given by the collection
  of all $[w_I]$'s can be lifted to $\tilde{\textup{B}}_r$ by $\pi^*$.
  Since $\tilde{\textup{B}}_r$ is simply connected, the lifted section
  is actually global. Given the local trivialization determined by the
  choice of $I$ near $\Lambda$, this global section gives rises to a
  unique map $[\tilde w] : \tilde{\textup{B}}_r \to \RM P^{n-1}$
  satisfying \emph{\ref{prop:item:1}.}  and \emph{\ref{prop:item:2}},
  obtained by parallel transport;
  equation~\eqref{equ:rotation_number_holonomy} follows
  from~\eqref{equ:changement-action}.
\end{demo}

Here is how this proposition spells out more concretely: we can find
``global action variables'' on $\tilde{\textup{B}}_r$, \emph{i.e.}\ a
map $\tilde I:\tilde{\textup{B}}_r\to\RM^n$ such that for each small
open ball $\tilde U$ on $\tilde{\textup{B}}_r$, the restriction of
$\tilde I$ to $\tilde U$ descends to action variables $I_{\tilde U}$
for the initial integrable system $F$ on $U=\pi(\tilde U)$.  Moreover,
we may assume that $I_{\tilde U_0}=I$ if $\tilde U_0$ is a
neighborhood of the initial torus $\tilde \Lambda$. For each
$\tilde U$, $[\tilde w]$ is defined on $\tilde U$ by
  \begin{equation}
    [\tilde w] = [w_{I_{\tilde U}}] \circ \pi \;,
    \label{equ:wU}
  \end{equation}
  where $[w_{I_{\tilde U}}]$ given by
  Definition~\ref{defi:rotation}. The fact that $[\tilde w]$ is a
  smooth, single-valued function on $\tilde{\textup{B}}_r$ follows from the
  transition formula~\eqref{equ:changement-action}, where by
  definition here all transition maps $A$ are the identity matrix,
  because $I_{\tilde U} = I_{\tilde U'}$ on
  $\pi(\tilde U) \cap \pi(\tilde U')$. 

\subsection{The semitoric case}
\label{sec:semitoric-case}

A semitoric system, in a broad sense, is a two degree of freedom
integrable Hamiltonian system with a global $\SM^1$ symmetry; more
precisely, the Hamiltonian $H$ commutes with a smooth function $J$
whose Hamiltonian flow $t\mapsto \phy_J^t$ is $2\pi$-periodic, except
at fixed points. Under some additional hypothesis (which ensures, for
instance, that the level sets of the energy-momentum map $(J,H)$ are
connected), such systems have been introduced in~\cite{san-polytope}
and classified in~\cite{san-alvaro-I,san-alvaro-II}. Semitoric
systems, in particular under their quantum version (see
Section~\ref{sec:q-semitoric}) are ubiquitous in physics, when
investigating the coupling between spin, atoms, q-bits, small
molecules, etc. For a discussion and more references on this topic,
see~\cite{san-yohann21}. The symplectic theory of semitoric systems
borrows both from almost-toric
systems~\cite{symington-four,leung-symington} and Hamiltonian $\SM^1$
spaces~\cite{karshon-S1}. Recently, the initial definition has been
generalized in several ways, allowing the inclusion of famous examples
like the Spherical Pendulum,
see~\cite{san-pelayo-ratiu-affine,hohloch-sabatini-sepe-symington-18,
  ratiu-wacheux-zung-17}. 

\begin{figure}[ht]
  \centering \includegraphics[width=0.45\linewidth]{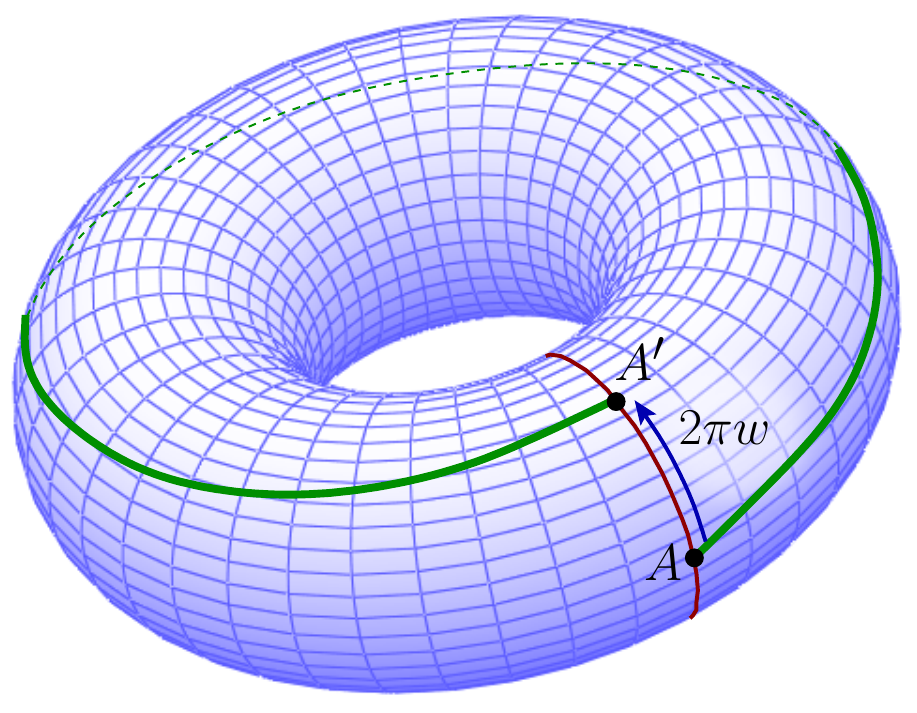}
  \caption{The semitoric rotation number.}
  \label{fig:semtoric}
\end{figure}

In the present work, we don't need the classification
of~\cite{san-alvaro-I,san-alvaro-II}, and hence we will use the
terminology ``semitoric'' is a very general acceptance (which
includes, for instance, lagrangian leaves with several focus-focus
singularities~\cite{palmer-pelayo-tang2019semitoric}):
\begin{defi}
  An integrable system $F=(J,H)$ on a 4-dimensional symplectic
  manifold $M$ will be called \emph{semitoric} if the map
  $F:M\to\RM^2$ is proper, and the function $J$ is a momentum map for
  an effective Hamiltonian $\SM^1$-action on $M$.
\end{defi}
The existence of such a global symmetry has the interesting
consequence that the rotation number for the system $F=(J,H)$ is
defined in a more natural way, and can be interpreted as an angle in
$\RM/\ZM$:

\begin{defi}\label{defi:semitoric-rotation}
  If $\Lambda$ is a Liouville torus for the semitoric system $(J,H)$,
  the (semitoric) \emph{rotation number} $w(\Lambda)\in\RM/\ZM$ is
  defined as follows. Take a point $A\in\Lambda$ and let $\SM^1_A$ be
  the orbit of this point under the flow of $J$. Consider the $H$-flow
  of $A$ and denote by $A'$ the first return point of the trajectory
  starting from $A$ on the orbit $\SM^1_A$. (See
  Figure~\ref{fig:semtoric}.) Then $2\pi w(\Lambda)$ is the time
  necessary to flow, under the action of $J$, from $A$ to $A'$:
  \[
    \varphi_J^{2\pi w(\Lambda)}(A) = A'.
  \]
\end{defi}
This definition is of course not new. It was popularized in particular
by Arnold in the treatment of perturbations of integrable systems, see
for instance~\cite[Section 3 G]{arnold-ordinary}.  In order to state
in which sense this new rotation number coincides with the general
rotation number from Definition~\ref{defi:rotation} (more specifically
Equation~\eqref{equ:ratio}), we recall the following basic fact from
the theory of semitoric systems (see for instance~{\cite[Lemma
  5.1]{san-pelayo-ratiu-affine}}):
\begin{lemm}
  \label{lemm:semitoric-action}
  Let $(J,H)$ be a semitoric system, and let $\Lambda$ be a Liouville
  torus. Then there exists a set of oriented action integrals near
  $\Lambda$ of the form $(J,I_2)$.
\end{lemm}
It will be convenient to call such sets of action integrals
``semitoric action variables''.
\begin{prop}
  \label{prop:semitoric}
  Let $(J,H)$ be a semitoric system, and let $\Lambda$ be a Liouville
  torus.  Let $I$ be semitoric action variables near $\Lambda$.  Let
  $w(\Lambda)$ be the semitoric rotation number of
  Definition~\ref{defi:semitoric-rotation}, and let $w_{I}(\Lambda)$
  be the rotation number from Equation~\eqref{equ:ratio}. Then
  \[
    w(\Lambda) = (w_{I}(\Lambda) \mod \ZM).
  \]
\end{prop}
\begin{demo}
  Since $I=(J,I_2)$, we have $\ham{H} = a_1 \ham{J} + a_2 \ham{I_2}$,
  for some smooth functions $a_j=a_j(\Lambda)$, such that
  $w_{I}(\Lambda)=a_1(\Lambda)/a_2(\Lambda)$. Notice that, since
  $\ham{H}$ and $\ham{J}$ must be linearly independent on $\Lambda$,
  we must have $a_2(\Lambda)\neq 0$. Hence $w_{I}$ is a standard real
  number.

  From Definition~\ref{defi:semitoric-rotation}, there exists a
  positive time $\tau=\tau(\Lambda)$ such that the flow of $\ham{H}$
  at time $\tau$ followed by the flow of $\ham{J}$ at time
  $-2\pi w(\Lambda)$ sends $A$ to itself. In other words, the time-1
  flow of
  \begin{equation}
    \ham{\textup{per}}:= \tau(\Lambda)\ham{H} - 2\pi w(\Lambda)
    \ham{J}
    \label{equ:Hper}
  \end{equation}
  is the identity on $\Lambda$. By the theory of action-angle
  variables, $\ham{\textup{per}}$ is the Hamiltonian vector field of
  the Hamiltonian
  \[
    H_{\textup{per}} := \int_\gamma \alpha,
  \]
  where $\alpha$ is a primitive of the symplectic form in a
  neighborhood of $\Lambda$, and the cycle $\gamma$ is the homology
  class of a periodic trajectory of $H_{\textup{per}}$ on
  $\Lambda$. Moreover, the action of $\ham{\textup{per}}$ must be
  effective, otherwise there would be a periodic trajectory of
  $\ham{\textup{per}}$ of period $1/2$, which would imply
  by~\eqref{equ:Hper} that the image of $A$ by the flow of $\ham{H}$
  at time $\tau(\Lambda)/2$ would return to the $\SM^1$-orbit of $A$,
  contradicting the definition of $\tau(\Lambda)$. Hence the pair
  $(J,H_{\textup{per}})$ is a set of action variables near $\Lambda$,
  which implies $H_{\textup{per}} = \pm 2\pi I_2 + 2\pi k J$ for some
  $k\in\ZM$.  From~\eqref{equ:Hper} we get
  \[
    \tau(\Lambda)\ham{H} = \pm 2\pi \ham{I_2} + 2\pi (w(\Lambda) -
    k)\ham{J}.
  \]
  The assumption that $I=(J,I_2)$ is oriented means
  $\deriv{I_2}{H}>0$, which implies that only the sign $+ 2\pi I_2$
  can occur, and hence we get $w_{I}(\Lambda) = w(\Lambda) - k$.
\end{demo}

  In concrete examples, computing the rotation number can be a hard
  task. But, even in the purely classical theory, it is an important
  piece of information about the system, related to the
  ``tennis-racket effect''~\cite{cushman-book}, and is intimately
  related to the symplectic invariants classifying semitoric systems,
  as demonstrated in the case of the spherical pendulum
  by~\cite{dullin-pendulum}.

\section{The Quantum Rotation Number}
\label{sec:quantum}

The quantum version of the image of the map $F=(f_1,\dots, f_n)$,
where $f_j$ are commuting Hamiltonians, is the joint
spectrum~\cite{colinI} of commuting operators $P_1,\dots, P_n$. With
this in mind, it is very natural to try to define a ``quantum rotation
number'' by replacing energies by ``quantized energies''. That this is
indeed possible will follow from the mathematical version of the
semiclassical Bohr-Sommerfeld rules, which requires a rigorous
framework establishing the relationship between quantum theory
(functional analysis of operators on Hilbert spaces) and classical
mechanics (symplectic geometry of Hamiltonian systems). In this paper,
for this purpose, we use the well-established theory of semiclassical
pseudodifferential quantization.

In Section~\ref{sec:bs}, we recall the definition of the joint
spectrum, and state the Bohr-Sommerfeld rules (Theorem~\ref{theo:bs}).
From this, in Section~\ref{sec:asympt-latt}, we extract the crucial
notion of asymptotic lattices and their good (and linear) labellings,
which we believe is one of the main contributions of our paper. This
leads, in Section~\ref{ss:quantrot}, to a very natural definition of
quantum rotation number, which converges to the classical rotation
number as $h\to 0$. From this point, we specialize to $n=2$. The
important case of semitoric systems is discussed in
Section~\ref{sec:q-semitoric}.  Section~\ref{sec:labelling-algorithm}
can be thought of as the technical core of the paper; we construct
robust algorithms for detecting linear labellings of general
2-dimensional asymptotic lattices. Theses algorithms are applied in
Section~\ref{sec:inverse} to finally solve the inverse problem for the
rotation number.

\subsection{Bohr-Sommerfeld quantization rules for quantum integrable
  systems}\label{sec:bs}

Bohr-Sommerfeld rules bear several meanings in mathematics. We use
here the original formulation coming from quantum theory, which
describes the joint spectrum of commuting operators in terms of the
underlying symplectic geometry of their symbols.

Let $P$ be a (possibly unbounded) selfadjoint operator acting on a
Hilbert space $\mathcal{H}$. The spectral theorem constructs from $P$
its so-called spectral measure, which is a projector-valued measure on
$\RM$. The support of the spectral measure is the spectrum of $P$,
denoted by $\sigma(P)$. The spectrum of $P$ in an interval
$\interval{E_1}{E_2}$ is called discrete if any $\lambda\in \sigma(P)$
is isolated in $\sigma(P)\cap \interval{E_1}{E_2}$, and for any
$\epsilon>0$ small enough, the spectral projector associated with the
interval $\interval{\lambda-\epsilon}{\lambda+\epsilon}$ has finite
rank.

The selfadjoint operators $P_1,\dots,P_n$ are said to pairwise commute
if their spectral measures commute, which we denote by
$[P_j,P_k]=0$. In this case, one can define the joint spectral
measure, and by definition the \emph{joint spectrum} of
$(P_1,\dots,P_n)$ is the support (in $\RM^n$) of this measure. We
shall denote it by $\Sigma(P_1,\dots,P_n)$. The joint spectrum is
called discrete if all $\lambda\in \Sigma(P_1,\dots,P_n)$ are isolated
and, for any small enough compact neighborhood $K$ of $\lambda$, the
joint spectral projector on $K$ has finite rank.

Let $P$ be a semiclassical $\h$-pseudodifferential operator on
$X=\RM^n$, or on a compact manifold $X$ of dimension $n$, where the
semiclassical parameter $\h$ varies in an interval
$\interval[open left]{0}{\h_0}$ for some $\h_0>0$. More precisely, in
the $\RM^n$ case, we will assume that $P$ is the semiclassical Weyl
quantization of a symbol in $S(m)$ where $m$ is an order function in
the sense of~\cite[Definition 7.4]{dimassi-sjostrand}: there exist
$C>0,N>0$ such that
\[
  1 \leq m(X) \leq C \langle X - Y \rangle ^N m(Y) \qquad \forall X,Y
  \in \RM^{2n},
\]
where $\langle X \rangle := (1+\norm{X}^2)^{1/2}$. For instance, one
can take $m(X):=\langle X \rangle^d$, for some $d\in\RM$.  A function
$a=a(x,\xi;\h)\in\Cinf(\RM^{2n})$ belongs to the class $S(m)$ if, by
definition,
\[
  \forall \alpha\in\NM^{2n}, \forall \h\in(0,\h_0],\quad
  \abs{\partial^\alpha a(x,\xi;\h)}\leq C_\alpha m(x,\xi),
\]
and, in this paper, we shall always assume that such a symbol $a$
admits an asymptotic expansion $a\sim \sum_{j\geq 0} \h^j a_j$ in
non-negative integral powers of $\h$ in this topology: for any $N>0$,
\begin{equation}
  a(x,\xi;\h) - \sum_{j=0}^N \h^j a_j(x,\xi) \in \h^{N+1}S(m).
  \label{equ:expansion}
\end{equation}

We then say that a pseudodifferential operator $P$ belongs to $S(m)$
(using the same notation, since there will be in general no ambiguity)
if its Weyl symbol $a$ belongs to $S(m)$. The first term $a_0$
in~\eqref{equ:expansion} is called the principal symbol of $P$.  Other
classes of pseudodifferential operators can be used,
see~\cite{charbonnel}, but we stick to $S(m)$ for its ease of use and
ample documentation (see for
instance~\cite{dimassi-sjostrand,martinez-book,zworski-book-12}).

In the manifold case, we will assume that $P$ belongs to the
Kohn-Nirenberg class $S^d(X)$, ($d\in\RM$), which means that in local
coordinates, after cut-off in the position variable $x\in X$, $P$ can
be written as a pseudodifferential operator with a symbol
$a(x,\xi;\h)$ such that (see for instance~\cite[Chapter
14]{zworski-book-12})
\[
  \forall (\alpha,\beta)\in\NM^{n}\times\NM^n, \forall
  \h\in(0,\h_0],\quad \abs{\partial^\alpha_x \partial^\beta_\xi
    a(x,\xi;\h)} \leq C_{\alpha,\beta} \langle \xi
  \rangle^{d-\abs{\beta}}.
\]
In this case, to unify notations, we call
$m(x,\xi):=\langle \xi\rangle^d$ the corresponding order function, and
denote $S(m):=S^d(X)$.

The operator $P$ is said to be elliptic at infinity in $S(m)$ if
$P\in S(m)$ and $\abs{a_0(x,\xi)} \geq c m(x,\xi)$ for some $c>0$ and
for large $(x,\xi)$ in the $\RM^n$ case; and similarly
$\abs{a_0(x,\xi)} \geq c \langle \xi \rangle^d$ in the manifold case.

A \emph{quantum integrable system} is the data of $n$ commuting
selfadjoint pseudodifferential operators $P_1,\dots,P_n$ in $S(m)$,
such that the map of the principal symbols $F:=(p_1,\dots,p_n)$
defines a classical integrable system in $T^*X$. For simplicity, we
shall always assume that $F$ is globally proper and that the operator
$Q:=P_1^2+\cdots + P_n^2$ is elliptic at infinity in $S(m^2)$. This
implies that $m$ is unbounded. (We could also include the case of
bounded $m$ by imposing properness of $F$ only onto some compact set,
and considering the intersection of the joint spectrum with a strict
compact subset of this set.) Note that, for pseudodifferential
operators, the commutation property $[P_j,P_k]=0$ is equivalent to the
weak commutation (see~\cite{charbonnel})
\[
  \forall u\in \Cinf_0(\RM^n), \quad P_j P_k u = P_k P_j u.
\]
In this case, for any $f\in\Cinf_0(\RM^n)$, the operator
$f(P_1,\dots,P_n)$ (obtained via the joint spectral measure) is a
pseudodifferential operator in $S(m^{-N})$ for any $N\in\NM$
(\cite[Chapter 8]{dimassi-sjostrand}).

Typical examples of quantum integrable systems (with unbounded
symbols) are given by two degree of freedom Schrödinger operators with
an axi-symmetric potential; these can occur either on
$\RM^2$~\cite{charbonnel}, or on a Riemannian surface of
revolution~\cite{colinI,colinII}, and can be used to obtain efficient
numerical methods, see~\cite{BDM-99}. Other classical examples of
quantum integrable systems, which do not feature a global $\SM^1$
symmetry, include the C-Neuman problem and the Laplace-Beltrami
operator on an ellipsoid, see~\cite{toth-quadratic}.

\begin{prop}[\cite{charbonnel}]
  If $P_1,\dots,P_n$ is a quantum integrable system, satisfying the
  above hypothesis, then its joint spectrum $\Sigma(P_1,\dots,P_n)$ is
  discrete when $\h$ is small enough.
\end{prop}
\begin{demo}
  Let $\lambda:=(\lambda_1,\dots,\lambda_n)\in \Sigma(P_1,\dots,P_n)$,
  and let $K\subset\RM$ be a compact neighborhood of $\lambda$. Since
  $F$ is proper, the map $f\circ F$ has compact support for any
  $f\in\Cinf_0(\RM^n)$. Hence the pseudodifferential operator
  $f(P_1,\dots,P_n)$ is compact for small $\h$. Choosing $f\equiv 1$
  on a neighborhood of $K$, we see that the spectral projector
  $1_K(P_1,\dots,P_n) = f(P_1,\dots,P_n)1_K(P_1,\dots,P_n)$ is
  compact; hence its rank is finite.
\end{demo}

The following result, which is a mathematical justification of the old
Bohr-Sommerfeld rule from quantum mechanics, states that the joint
spectrum of a quantum integrable system takes the form of an
approximate lattice in any neighborhood of a regular value of $F$.  It
was first obtained by Colin de Verdière in the homogeneous
setting~\cite{colinII}, and then by Charbonnel in the semiclassical
setting~\cite{charbonnel}. A purely microlocal approach was proposed
in~\cite{san-focus} (see also~\cite{san-panoramas}). When $n=1$ and
$P$ is a Schrödinger operator, explicit methods for second order ODEs
can be employed, see~\cite{yafaev-multiple}.

\begin{theo}[\cite{colinII,charbonnel}]\label{theo:bs}
  Let $P_1,\dots,P_n$ is a quantum integrable system satisfying the
  above hypothesis. Let $c\in\RM^n$ be a regular value of the
  principal symbol map $F=(p_1,\dots,p_n)$. Assume that the fiber
  $F^{-1}(c)$ is connected. Then we can describe the joint spectrum
  around $c$ as follows.
  \begin{enumerate}

  \item\textbf{ (joint eigenvalues have multiplicity one) }There
    exists a non empty open ball $B\subset\RM^n$ (independent of $\h$)
    around $c$ and $\h_0>0$ such that for any
    $\lambda\in\Sigma(P_1,\dots,P_n) \cap B$ and for all $\h<\h_0$,
    the joint spectral projector of $(P_1,\dots,P_n)$ onto the ball
    $B(\lambda,\h^2)$ has rank 1.

  \item\textbf{ (the joint spectrum is a deformed lattice) } There
    exists a bounded open set $U\subset\RM^n$ and a smooth map
    $G_\h:U\to\RM^n$ admitting an asymptotic expansion in the $\Cinf$
    topology:
    \[
      G_\h = G_0 + \h G_1 + \h^2 G_2 + \cdots
    \]
    such that the joint eigenvalues in $\Sigma(P_1,\dots,P_n) \cap B$
    are given by the quantities
    \begin{equation}
      \lambda = G_\h(\h k_1, \dots, \h k_n) + \O(\h^\infty),
      \label{equ:bs}  
    \end{equation}
    where $(k_1,\dots,k_n) \in \ZM^n$ is such that
    $G_\h(\h k_1, \dots, \h k_n)\in B$, and the $\O(\h^\infty)$
    remainder is uniform as $\h\to 0$. Moreover, $G_0$ is a
    diffeomorphism from $U$ onto a neighborhood of $\overline B$ given
    by the action-angle theorem:
    \begin{equation}
      F = G_0(I_1,\dots, I_n),
      \label{equ:G0}    
    \end{equation}
    where $(I_1,\dots, I_n)$ is a set of action variables for the
    classical system, defined in a neighborhood of $F^{-1}(c)$, and
    $U$ is a neighborhood of $I_c = I(F^{-1}(c))$.
  \end{enumerate}
\end{theo}
It is important to notice that, contrary to the classical case, one
can not always assume that the actions $(I_1,\dots, I_n)$ take values
in a neighborhood of the origin (in other words, $I_c$ is not
necessarily zero). In fact $I_c$ is given by the integrals of the
Liouville 1-form of the cotangent bundle $T^*X$ over a basis of cycles
of the torus $F^{-1}(c)$.
\begin{rema}
  That this statement can be generalized to multiple connected
  components of $F^{-1}(c)$ is ``well known to specialists''; however
  to our knowledge that generalization cannot be found in the
  literature. We haven't included it here because it will be important
  for us to make the connectedness hypothesis throughout the
  paper. However, we expect it to become necessary for the future
  study of integrable systems with hyperbolic singularities, which are
  very common in the spectroscopy of small molecules (eg. LiNC/NCLi,
  see~\cite{joyeux-sado-tenny}), and related to the so-called
  fractional monodromy, see for instance~\cite{giacobbe08}. On the
  mathematical side, several authors have already attacked the
  semiclassical analysis near the hyperbolic separatrix, for instance,
  for 1D Schrödinger operators~\cite{marz,bleher94}, Liouville
  surfaces~\cite{bleher}, for general 1D pseudodifferential
  operators~\cite{colin-p,colin-p2}, and transversally hyperbolic 2D
  systems~\cite{san-colin}.
\end{rema}
\begin{rema}
  The term $G_1$ can be expressed in terms of Maslov indices and the
  subprincipal symbols of the operators $P_1$, see~\cite{san-focus};
  we won't need this particular formula here.
\end{rema}
\begin{rema}\label{rema:toeplitz}
  We have not attempted here to extend our analysis to commuting
  Berezin-Toeplitz operators on prequantizable symplectic manifolds
  (see~\cite{lefloch-book} for a nice introduction to Berezin-Toeplitz
  quantization). However, the notion of asymptotic lattices that we
  develop in this paper is quantization-agnostic, and hence, thanks to
  the work of Charles~\cite{charles-bs}, we believe that most of our
  results should be adaptable to that situation.
\end{rema}

\subsection{Asymptotic lattices and good labellings}
\label{sec:asympt-latt}

Sets of points in $\RM^n$ that are described by an equation of the
type of~\eqref{equ:bs} can be called 'asymptotic lattices',
see~\cite{san-mono}. Indeed, when $\h$ is small enough, the map
$G_\h: U \to \RM^n$ is a diffeomorphism onto its image, and hence it
sends the local lattice $\h\ZM^n\cap U $ one-to-one into the joint
spectrum $\Sigma(P_1,\dots,P_n) \cap B \mod \O(\h^\infty)$. This also
implies that for any strictly smaller ball $\tilde B\Subset B$, the
map obtained by the restriction to lattice points
\[
  G_\h: \h\ZM^2 \cap G_\h^{-1}(\tilde B) \to \Sigma(P_1,\dots,P_n)
  \cap \tilde B \mod \O(\h^\infty)
\]
is onto.  The goal of this section is to make these observations
precise by introducing the formal definition of asymptotic lattices
and deriving some of their properties; this analysis partly builds on
and extends~\cite{san-mono}.

\subsubsection{Asymptotic charts and asymptotic lattices}

\begin{defi}
  \label{defi:lattice}
  Let $B\subset\RM^n$ be a simply connected bounded open set. Let
  $\Lh\subset B$ be a discrete subset of $B$ depending on the small
  parameter $\h\in\cI$, where $\cI\subset\RM^*_+$ is a set of positive
  real numbers admitting $0$ as an accumulation point. We say that
  $(\Lh, \cI, B)$ is an \emph{asymptotic lattice} if the following two
  statements hold:
  \begin{enumerate}
  \item
    \label{item:lattice-distance}
    there exist $\h_0>0$ , $\epsilon_0>0$, and $N_0\geq 1$ such that
    for all $\h\in\cI\cap\interval[open left]{0}{\h_0}$,
    \[
      \h^{-N_0}\min_{\scriptstyle (\lambda_1,\lambda_2)\in \Lh^2 \atop
        \scriptstyle \lambda_1\neq\lambda_2} \norm{\lambda_1-\lambda_2} \geq
      \epsilon_0;
    \]
  \item
    there exist a bounded open set $U\subset\RM^n$ and a family of
    smooth maps $G_\h:U\to\RM^n$, for $\h\in\cI$, admitting an
    asymptotic expansion in the $\Cinf(U)$ topology:
    \begin{equation}
      G_\h = G_0 + \h G_1 + \h^2 G_2 + \cdots
      \label{equ:G}
    \end{equation}
    with $G_j\in\Cinf(U;\RM^n)$, for all $j\geq 0$. Moreover, $G_0$ is
    an orientation preserving diffeomorphism from $U$ onto a
    neighborhood of $\overline B$, and
    \[
      G_\h(\h\ZM^n \cap U) = \Lh + \O(\h^\infty)\quad \text{ inside }
      B,
    \]
    by which we mean: there exists a sequence of positive numbers
    $(C_N)_{N\in\NM}$, such that for all $\h\in\cI$, the following
    holds.
    \begin{enumerate}
    \item For all $\lambda\in \Lh$ there exists $k\in\ZM^n$ such that
      $\h k \in U$ and
      \begin{equation}
        \forall N\in\NM, \qquad \norm{\lambda - G_\h(\h k)} \leq C_N \h^N.
        \label{equ:chart}  
      \end{equation}
    \item \label{item:def-asym-latt2b} For every open set
      $\tilde U_0\Subset G_0^{-1}(B)$, there exists $\tilde \h_0>0$
      such that for all
      $\h\in\cI\cap\interval[open left]{0}{\tilde\h_0}$, for all
      $k\in\ZM^n$ with $\h k\in \tilde U_0$, there exists
      $\lambda\in\Lh$ such that~\eqref{equ:chart} holds.
    \end{enumerate}
  \end{enumerate}
  Such a map $G_\h$ will be called an \emph{asymptotic chart} for
  $\Lh$. We will also denote it by $(G_\h,U)$.
\end{defi}
Notice that the orientation preserving requirement is not a true
restriction; indeed, one may pick an element $A\in\textup{GL}(n,\ZM)$
with $\det A=-1$, and then $(G_\h\circ A, A^{-1}U)$ is a good
asymptotic chart if and only if $(G_\h, U)$ is an ``orientation
reversing asymptotic chart''.

For shortening notation, we shall often simply call $\Lh$ an
asymptotic lattice, instead of the triple $(\Lh, \cI, B)$.  Of course
the simplest asymptotic lattice is the square lattice $\h\ZM^n\cap B$,
with chart $G_\h = \textup{Id}$.  In general, if $G_\h$ is any map
satisfying~\eqref{equ:G} with $G_0$ a diffeomorphism onto a
neighborhood of $\overline{B}$, then $G_\h(\h\ZM^n\cap U)\cap B$ is an
asymptotic lattice.

With this definition, we may now rephrase Theorem~\ref{theo:bs} as
follows.
\begin{theo}
  \label{theo:bs-bis}
  With the hypothesis of Theorem~\ref{theo:bs}, near any regular value
  $c\in\RM^n$ with connected fiber, the joint spectrum of the quantum
  integral system $(P_1,\dots,P_n)$, with
  $\h\in\cI=\interval[open left]0{\h_0}$ is an asymptotic lattice.
\end{theo}

Let us state some elementary properties of asymptotic charts, for
further reference. For subsets $A,B$ of $\RM^n$, we use the notation
$A\Subset B$ to mean that $\overline{A}$ is compact and contained in
$B$. For a map $G:U\to\RM^n$, where $U\subset\RM^d$ is an open set, we
denote by $G':U\to\mathcal{L}(\RM^d; \RM^n)$ the full derivative (or
linear tangent map) of $G$.
\begin{lemm}
  \label{lemm:chart-pot-pourri}
  Let $(\Lh, \cI, B)$ be an asymptotic lattice, and let $(G_\h, U)$ be
  an asymptotic chart for it. Let $\h_0>0$. The following properties
  of $G_\h$ hold.
  \begin{enumerate}
  \item \label{item:1} All derivatives of $G_\h$ are bounded on any
    compact subset of $U$, uniformly for $\h\leq \h_0$:
    \begin{equation*}
      \forall \tilde U\Subset U, \forall\alpha\in\NM^n, 
      \exists C_{\tilde U, \alpha}>0,
      \forall \h\in\cI\cap\interval[open left]{0}{\h_0}, 
      \quad \sup_{\tilde U}\norm{\partial^\alpha G_\h} \leq C_{\tilde U, \alpha}.
    \end{equation*}
  \item \label{item:2} For any $\tilde U\Subset U$, there exist
    non-negative constants $M, M'$ such that
    \begin{equation}
      \label{equ:G0-bounded}
      \forall \h\in\cI\cap\interval[open left]{0}{\h_0}, \qquad 
      \sup_{\tilde U}\norm{G_\h - G_0} \leq \h M, \quad 
      \sup_{\tilde U}\norm{G'_\h - G_0'} \leq \h M'.
    \end{equation}
  \item \label{item:3} For any $\tilde U\Subset U$, for any $N\geq 0$,
    there exists a non-negative constant $R_N$ such that, for all
    $\h_1, \h_2\in\cI\cap\interval[open left]{0}{\h_0}$,
    \begin{equation*}
      \sup_{\tilde U}\norm{G_{\h_2} - G_{\h_1}} \leq 
      R_N(\abs{\h_2 - \h_1} + \h_1^N + \h_2^N).
    \end{equation*}
  \item \label{item:4} For any $\tilde U\Subset U$, there exist
    non-negative constants $L_0, L_1$ such that, for all
    $\h\in\cI\cap\interval[open left]0{\h_0}$,
    $\forall \xi_1,\xi_2 \in \tilde U$, if the segment $[\xi_1,\xi_2]$
    is contained in $\tilde U$ then
    \begin{align}
      \label{equ:G-taylor1}
      \norm{G_\h(\xi_2) - G_\h(\xi_1)}&  \leq L_0\norm{\xi_2 - \xi_1}\\
      \label{equ:G-taylor2}
      \norm{G_\h(\xi_2) - G_\h(\xi_1) - G_\h'(\xi_1)\cdot(\xi_2 -
      \xi_1)}& \leq L_1\norm{\xi_2 - \xi_1}^2.
    \end{align}
  \item \label{item:5} The map $G_\h$ is formally invertible: there
    exists a smooth map $F_\h$, defined on a neighborhood of
    $\overline{B}$, that admits an asymptotic expansion
    \begin{equation}
      F_\h = F_0 + \h F_1 + \h^2 F_2 + \cdots\label{equ:F-expan}
    \end{equation}
    such that
    \begin{equation}
      \label{equ:Fh}
      F_\h \circ G_\h = \textup{Id} + \O(\h^\infty),
      \quad \text{ and } \quad
      G_\h \circ F_\h = \textup{Id} + \O(\h^\infty).
    \end{equation}
    In particular, $F_0 = G_0^{-1}$.
  \item \label{item:6} For any $\tilde U\Subset U$, there exists
    $\h_0>0$ such that for all
    $\h\in\cI\cap\interval[open left]0{\h_0}$, the restriction of
    $G_\h$ to $\tilde U$ is invertible onto its image. Moreover, its
    inverse $G_\h^{-1}$ is smooth and admits an asymptotic expansion
    in powers of $\h$ equal to that of $F_\h$ in~\eqref{equ:F-expan}.
  \item \label{item:7}Given $\xi\in U$, the map $\h\mapsto G_\h(\xi)$
    is not necessarily continuous. However, there exists a smooth map
    $\tilde G \in \Cinf(\interval{-\h_0}{\h_0})\times U)$ such that,
    for all $\xi\in U$,
    $\forall \h\in \cI\cap\interval[open left]0{\h_0}$,
    $G_\h(\xi) = \tilde G(\h, \xi) + R_\h(\xi)$, and
    $R_\h = \O(\h^\infty)$ in the $\Cinf(U)$ topology. In particular,
    $\tilde G(0,\xi) = G_0(\xi)$.
  \item \label{item:8} If $\tilde G_\h = G_\h + \O(\h^\infty)$ in
    $\Cinf(U)$, then $\tilde G_\h$ is also an asymptotic chart for
    $(\Lh, \cI, B)$.
  \end{enumerate}
\end{lemm}
\begin{demo}
  Items~\ref{item:1}, \ref{item:2} and~\ref{item:3} directly follow
  from the asymptotic expansion~\eqref{equ:G}. Taylor's formula and
  Item~\ref{item:1} gives Item~\ref{item:4}. Item~\ref{item:5} is a
  standard consequence of~\eqref{equ:G}: the functions $F_j$ can be
  defined by induction, using Taylor expansions, and $F_\h$ is
  obtained by the Borel lemma. Item~\ref{item:7} also follows from a
  Borel summation of the formal series $\sum \h^k G_k$. To prove
  Item~\ref{item:6}, observe that near any $\xi\in U$, if $\h$ is
  small enough, $G_\h'(\xi)$ is invertible due to the invertibility of
  $G'_0$ and~\eqref{equ:G0-bounded}. Hence on a compact subset of $U$,
  we may find $\h_0$ such the local inversion theorem applies if
  $\h < \h_0$, showing that $G_\h$ is a local diffeomorphism for all
  $\h\in\cI\cap\interval[open left]{0}{\h_0}$. Thus, there exists
  $\delta>0$ such that for any $\xi_0\in \tilde U$, the restriction of
  $G_\h$ to the ball $B(\xi_0,\delta)$ is a diffeomorphism onto its
  image. If follows from~\eqref{equ:Fh} that there exists $\h_0>0$
  depending only on $G_\h$ and $\tilde U$ such that
  \[
    \forall \h \in\cI\cap\interval[open left]0{\h_0}, \quad
    \forall\xi\in\tilde U, \quad \norm{F_\h\circ G_\h(\xi) - \xi} \leq
    \delta/3.
  \]
  Hence $G_\h(\xi_1) = G_\h(\xi_2)$ implies
  $\norm{\xi_1 - \xi_2}\leq 2\delta/3$, which in turn implies
  $\xi_1 = \xi_2$, proving the injectivity of $G_\h$ on $\tilde U$,
  which gives the result. Notice that, by choosing $\tilde U$ large
  enough, we can always ensure that $G_\h(\tilde U)$ contains a
  neighborhood of $\overline{B}$.
\end{demo}

\begin{rema}
  As usual in semiclassical analysis, equality of sets ``modulo
  $\O(\h^\infty)$'' has to be taken with care.  Given an asymptotic
  lattice $(\Lh, \cI, B)$, because of boundary effects, it is not true
  that the Hausdorff distance between $\Lh$ and
  $B\cap G_\h(\h \ZM^n \cap U)$ is $\O(\h^\infty)$. Consider the
  following example: $n=1$, $B=\interval[open]{0}{1}$,
  $\cI=\interval[open left]{0}{1}$, $\Lh = \h \ZM \cap B$. One can
  check that $(\Lh, \cI, B)$ is an asymptotic lattice with chart
  $G_\h=\textup{Id} : U=\interval[open]{-1}{2} \to U$. By virtue of
  Item~\ref{item:8} of Lemma~\ref{lemm:chart-pot-pourri}, the map
  $\tilde G_\h(\xi) = \xi + e^{-1/\h}$ is again an asymptotic chart
  for $\Lh$. However, if $\h$ is small enough,
  $\min G_\h(\h \ZM\cap U)\cap B = e^{-1/\h}$ while $\min \Lh=
  \h$. So, in this case, the Hausdorff distance is $\h-e^{-1/\h}$,
  which is not $\O(\h^\infty)$.
\end{rema}
\begin{rema}\label{rema:h-smooth}
  It follows from Item~\ref{item:7} of
  Lemma~\ref{lemm:chart-pot-pourri} that, up to changing $U$ for a
  smaller $U'\Subset U$, one can always choose an asymptotic chart
  $G_\h$ that is smooth with respect to $\h$. This, in turns, improves
  Item~\ref{item:3} by suppressing the term $\h_1^N + \h_2^N$.
\end{rema}

\subsubsection{Good labellings of an asymptotic lattice}

We now turn to the description of the asymptotic lattices themselves.
Using an asymptotic chart, we can put an integer label on each point
of the lattice:
\begin{lemm}\label{lemm:labelling}
  If $G_\h$ is an asymptotic chart for $\Lh$, then there exists a
  family of maps $k_\h: \Lh \mapsto \ZM^n$, $\h\in\cI$, which is
  unique for $\h$ small enough, and such that~\eqref{equ:chart} holds
  with $k=k_\h(\lambda)$, \emph{i.e.} $\h k_{\h}(\lambda) \in U$ and
  \begin{equation}
    \label{equ:chart-labelling}
    \forall\lambda\in\Lh, \quad \forall N\geq 0, \quad
    \norm{G_\h(\h k_\h(\lambda))-\lambda} \leq C_N\h^N.
  \end{equation}
  Moreover this map is injective.
\end{lemm}
\begin{demo}
  From Item~\ref{item:6} of Lemma~\ref{lemm:chart-pot-pourri}, let
  $\h_0>0$ and $\tilde U\Subset U$ be such that, for all
  $\h\leq \h_0$, $G_\h:\tilde U \to B_1$ is invertible onto $B_1$, a
  fixed bounded open neighborhood of $\overline{B}$.  Let
  $\lambda\in\Lh$. If $k$ and $\tilde k$ in $\ZM^n$
  satisfy~\eqref{equ:chart} (for some arbitrary family
  $(C_N)_{N\geq 0}$), then $G_\h(\h k)$ and $G_\h(\h \tilde k)$ belong
  to $B+B(0,C_N\h^N)$. We can assume that $\h_0$ is small enough so
  that the latter is contained in $B_1$. Thus, $\h k$ and
  $\h \tilde k$ belong to $\tilde U$.  Moreover,
  $\|G_\h(\h k) - G_\h(\h\tilde k)\| \leq 2C_N \h^N$, which entails
  $\|\h k - \h \tilde k\|\leq 2L_FC_N \h^N$, where $L_F$ is a uniform
  upper bound for the Lipschitz norm of $G_\h^{-1}$ on $\tilde
  U$. Choose any $N_1>1$ and take $\h_0$ small enough so that
  \hequa
  \begin{equation}
    2L_FC_{N_1}\h_0^{N_1-1}<1\,;
    \label{equ:h0:k-is-well-defined}
  \end{equation}
  \nequa
  we obtain, for any $\h\leq \h_0$, $\|k - \tilde k\| < 1$ and hence
  $k=\tilde k$. This shows that, for any $\h\leq \h_0$, any
  $\lambda\in\Lh$ is associated with a unique $k\in\ZM^n$ such
  that~\eqref{equ:chart} holds. We call this map $k_\h$. Notice that
  $k_\h = \frac{1}{\h}G_\h^{-1}\restr_{\Lh} + \O(\h^\infty)$. In order
  to prove injectivity, assume $k_\h(\lambda_1) = k_h(\lambda_2)$.
  From~\eqref{equ:chart} again we now get
  $\norm{\lambda_1-\lambda_2} \leq 2 C_N \h^N$ for all $N$. Let
  $(\epsilon_0, N_0)$ be the constants defined by
  item~\ref{item:lattice-distance} of Definition~\ref{defi:lattice};
  if
  \hequa
  \begin{equation}
    2C_N \h^N < \epsilon_0 \h^{N_0}\,,
  \end{equation}
  \nequa
  which will happen if one chooses $N>N_0$ and $\h_0$ small enough, we
  conclude that $\lambda_1=\lambda_2$.
\end{demo}

\begin{rema}\label{rema:N=1}
  Note that a consequence of the above lemma (and its proof) is that,
  once Condition~\ref{item:lattice-distance} in the definition of
  asymptotic lattices (Definition \ref{defi:lattice}) is fullfiled, it
  eventually holds with $N_0=1$, provided $\h_0$ is small
  enough. Indeed, let $\lambda_1,\lambda_2\in \Lh$ and let
  $k_i = k_\h(\lambda_i)$, $i=1,2$: for all $N\geq 0$,
  $\norm{\lambda_i - G_\h(\h k_i)} \leq C_N \h^N$. Then if
  $\norm{\lambda_1-\lambda_2} < \epsilon \h$ for some $\epsilon>0$, we
  obtain
  \[
    \norm{G_\h(\h k_1) - G_\h(\h k_2)} < 2 C_N \h^N + \epsilon\h,
  \]
  and hence
  $\norm{\h k_1 - \h k_2} < L_F(2C_N \h^{N-1} +
  \epsilon)\h$. Therefore, we may choose any $N>1$, (and, as above,
  any $N_1>1$, $N_2>N_0$), and conclude by the injectivity of $k_\h$
  that $\lambda_1 = \lambda_2$ as soon as $\epsilon$ and $\h_0$ verify
  \hequa
  \begin{equation}
    \label{equ:h0-lambda-unique}
    \h_0^{N_1-1}<\frac{1}{2L_FC_{N_1}} \,, \quad \h_0^{N_2-N_0} <
    \frac{\epsilon_0}{2C_{N_2}} \,, \quad \text{ and } \quad \epsilon
    + 2C_N \h_0^{N-1} \leq \frac{1}{L_F}\,.
  \end{equation}
  \nequa In other words, if $\epsilon$ and $\h_0$
  satisfy~\eqref{equ:h0-lambda-unique}, then
  \begin{equation}
    \label{equ:unique-ball}
    \forall \h\in\cI\cap\interval[open left]0{\h_0} \quad
    \forall\lambda\in\Lh,\quad B(\lambda,\epsilon\h) \cap \Lh =
    \{\lambda\}\,.
  \end{equation}

  In order to simplify the statement of our results, we could fix from
  now on $N_1=N=2$ and $N_2=N_0+1$, and~\eqref{equ:h0-lambda-unique}
  would be satisfied if we took $\epsilon=\frac{1}{3L_F}$ and
  \[
    \h_0 = \min\left(\frac{\epsilon_0}{4 C_{N_0+1}}\;;\; \frac{1}{6L_F
        C_2}\right).
  \]
  However, we believe that in some situations, having the possibility
  to optimize such estimates by taking large $N$'s can be useful.
\end{rema}

\begin{rema}
  A completely similar argument shows that for all $\lambda\in\Lh$,
  the integer vector $k_\h(\lambda)$ defined in
  Lemma~\ref{lemm:labelling} is the unique element $k\in\ZM^n$ such
  that $\h k \in U$ and
  \[
    \norm{\lambda - G_\h(k \h)} \leq \epsilon \h
  \]
  as soon as $L_F(C_N \h^{N-1} + \epsilon) < 1$ (which is implied
  by~\eqref{equ:h0-lambda-unique}). Of course, this, in turn, implies
  the much better estimate~\eqref{equ:chart-labelling}.
\end{rema}

\begin{defi}\label{defi:labelling}
  The map $k_\h$ from Lemma~\ref{lemm:labelling} will be called a
  \emph{good labelling} of $\Lh$.
\end{defi}
Let $\tilde B\Subset B$.  Then a good labelling is surjective onto
$(\frac{1}{\h}G_0^{-1}(\tilde B))\cap\ZM^n$, in a uniform way: by
Item~\ref{item:def-asym-latt2b} of Definition~\ref{defi:lattice} with
$\tilde U_0 = G_0^{-1}(\tilde B)$, there exists $\tilde\h_0>0$ such
that, for all $\h\leq \tilde\h_0$ there exists $\lambda\in\Lh$ with
$\lambda=G_\h(\h k) + \O(\h^\infty)$, and hence by
Lemma~\ref{lemm:labelling}, $k_\h(\lambda) = k$, as soon as $\h_0$
satisfies~\eqref{equ:h0:k-is-well-defined}.

It follows that the set $\Lh$ is always ``dense in $B$ as $\h\to 0$'',
by which we mean the following:
\begin{lemm}
  \label{lemm:lattice-dense}
  For any $c\in B$, there exists a family $(\lambda_\h)_{\h\in\cI}$
  with $\lambda_\h\in\Lh$ such that
  \[
    \lambda_\h = c + \O(\h).
  \]
\end{lemm}
\begin{demo}
  If $k(\h):=\lfloor \frac{\xi}{\h} \rfloor$ (the vector obtained by
  taking the integer part of all the components of $\frac{\xi}{\h}$),
  where $\xi:=G_0^{-1}(c)$, then $\norm{\h k(\h) - \xi} \leq
  \h$. Thus, if $\tilde B\Subset B$ is a neighborhood of $c$, then
  $\h k(\h)\in G_0^{-1}(\tilde B)$ for $\h$ small enough, and we may
  define $\lambda_\h$ to be the point in $\Lh$ associated with the
  label $k(\h)$.  Since
  $\norm{G_\h(\h k(\h)) - c} \leq L_0\norm{\h k(\h) - \xi} +
  \norm{G_\h(\xi) - c} \leq L_0 \h + M\h$, where $M$ is defined
  in~\eqref{equ:G0-bounded}, we get, for any $N\geq 0$,
  $\norm{\lambda_\h - c} \leq (L_0+M)\h + C_N\h^N$.
\end{demo}

\subsubsection{Linear labellings of an asymptotic lattice}

Equipped with a good labelling (or, equivalently, with an asymptotic
chart), an asymptotic lattice possesses the interesting feature that
its individual points inherit a well-defined ``smooth'' evolution as
$\h$ varies, in the following sense.  Let $\lambda_\h\in\Lh$. If
$\h_0$ is small enough, then by the injectivity of
Lemma~\ref{lemm:labelling} and~\eqref{equ:unique-ball}, for any
$\h\leq \h_0$, $\lambda_\h$ is the unique closest point to
$G_\h(\h k_\h(\lambda_\h))$ in $\Lh$. Thus, we may now fix the integer
$k_0 = k_{\h_0}(\lambda_{\h_0})$ and consider the evolution of the
corresponding point $k_{\h}^{-1}(k_0)\in \Lh$ as $\h$ varies close to
$\h_0$. Although this evolution may not be continuous, it is
$\O(\h^\infty)$-close to the smooth map
$\h\mapsto \tilde G(\h,\h k_0)$, where $\tilde G$ is the smooth
representative introduced in Lemma~\ref{lemm:chart-pot-pourri},
Item~\ref{item:7}.  Notice that the typical behaviour of the point
$G_\h(\h k_0)$, as $\h\to 0$, is to \emph{leave} the neighborhood
where the chart $G_\h$ is meaningful. Indeed, $\h k_0\to 0\in\RM^n$,
and there is no reason why $0$ should belong to $U$. This is analogous
to the well-known phenomenon that occurs for eigenvalues of operators
depending on a parameter, when we restrict our attention to a fixed
spectral window: trying to follow a particular ``mode'' (eigenvector)
as the parameter varies, the corresponding eigenvalue is likely to
escape the window. Good labellings capture this phenomemon. For our
purposes in this paper, it will be useful to disregard this ``drift'',
by defining the weaker notion of a good labelling ``modulo a
constant'', as follows.
\begin{defi}
  \label{defi:linear-labelling}
  A map $\bar k_\h: \Lh\to \ZM^n$ is called a \emph{linear labelling}
  for $\Lh$ if there exists a family
  $(\varkappa_\h)_{\h\in\interval[open left]{0}{\h_0}}$ of vectors in
  $\ZM^n$ such that $\bar k_\h + \varkappa_\h$ is a good labelling for
  $\Lh$.
\end{defi}
Of course any good labelling is a linear labelling; we shall see in
Section~\ref{sec:labelling-algorithm} below that linear labellings may
be easier to construct.

Another important property of asymptotic lattices is that they can be
equipped with an ``asymptotic $\ZM^n$-action'' in the following sense.
Fix a good labelling $k_\h$ for $\Lh$. Let $\tilde B \Subset B$ be
open, and let $\varkappa\in\ZM^n$ be a fixed integral vector. Let
$\hat B$ be an open set such that $\tilde B\Subset\hat B\Subset B$. If
$\h$ is small enough (depending on $\tilde B$, $\hat B$, and
$\varkappa$), for any $\lambda\in \Lh\cap \tilde B$, we have
$\h(k_\h(\lambda)+\varkappa)\in G_0^{-1}(\hat B)$; therefore, there
exists a unique point in $\Lh$, denoted by
$ \lambda \plus{k_\h} \varkappa$, such that
\begin{equation}
  \label{equ:translation}
  k_\h(\lambda \plus{k_\h} \varkappa) = k_\h(\lambda)+\varkappa.
\end{equation}
This property characterizes linear labellings, as follows.
\begin{prop}
  \label{prop:linear-labelling}
  Let $\Lh$ be an asymptotic lattice in the open set $B$ with good
  labelling $k_\h$.  Let $\tilde B\Subset B$ be open and
  connected. Let a map $\tilde k_\h: \Lh \to \ZM^n$, defined for
  $\h\in\cI$, commute with translations in the following sense: for
  any $\varkappa\in\ZM^n$, there exists $\h_\varkappa>0$ such that
  \begin{equation}
    \label{equ:ktilde}
    \forall \h \leq \h_\varkappa, \quad \forall \lambda\in\Lh\cap\tilde B,
    \qquad  
    \tilde k_\h(\lambda \plus{k_\h} \varkappa) = \tilde
    k_\h(\lambda)+\varkappa.
  \end{equation}
  Then $\tilde k_\h$ is a linear labelling for
  $(\Lh\cap\tilde B, \cI, \tilde B)$ (associated with the good
  labelling $k_\h$ restricted to $\tilde B$).
\end{prop}
In order to prove this proposition, we first show that
$\Lh\cap \tilde B$ is `connected by lattice paths'.
\begin{lemm}
  \label{lemm:lattice-connected}
  Let $\Lh$, $\h\in \cI$, be an asymptotic lattice in an open set
  $B$. Let $\tilde B\Subset \hat B \Subset B$, where $\tilde B$ and
  $\hat{B}$ are open and $\tilde B$ is connected.  Given a good
  labelling $k_\h$ for $\Lh$, there exists $\h_0 >0$ such that, given
  any $\h\in\interval[open left]{0}{\h_0}\cap \cI$, and given any pair
  of points $z_1,z_2\in \Lh\cap\tilde B$, there exists a finite
  sequence $(\epsilon^{(j)})_{j\in\{1,\dots,N\}}$ in $\{-1,0,1\}^n$,
  and a finite sequence of points $\lambda_j\in \Lh\cap\hat B$,
  $j=0,\dots, N$, such that
  \begin{equation}
    \label{equ:lattice-connectedness}
    \begin{gathered}
      \lambda_0 = z_1\,, \quad \forall j=0,\dots, N-1,
      \quad\lambda_{j+1} = \lambda_j \plus{k_\h} \epsilon^{(j+1)}\,,
      \quad \text{ and } \quad \lambda_N = z_2\,.
    \end{gathered}
  \end{equation}
\end{lemm}
\begin{demo}
  Let $(G_\h,U)$ be an asymptotic chart for $\Lh$ associated with
  $k_\h$ (see Lemma~\ref{lemm:labelling}). Let $c\in \tilde B$. Let
  $\tilde V\subset U$ be an open cube centered at
  $\xi=G_0^{-1}(c)$. Let us first prove the `lattice-connectedness'
  property~\eqref{equ:lattice-connectedness} for
  $\tilde B=G_0(\tilde V)$.  Let $\hat V\subset U$ be another open
  cube containing $\tilde V$. It follows from the existence and the
  asymptotic expansion of $G_\h^{-1}$ (Item~\ref{item:6} of
  Lemma~\ref{lemm:chart-pot-pourri}) that
  $G_\h^{-1}(\tilde B)\subset \hat V$, if $\h$ is small enough.
  Therefore, if $V\subset U$ is another open cube containing $\hat V$,
  Lemma~\ref{lemm:labelling} implies that, for $\h$ small enough,
  $\h k_\h(\lambda)\subset V$ for any $\lambda\in\Lh\cap \tilde B$,
  and $k_\h^{-1}((\h^{-1}V) \cap \ZM^n )\subset \hat B$. In
  $V\cap \h\ZM^n$, any two points $\xi_1,\xi_2$ can obviously be
  joined by acting by the canonical basis of the lattice. Applying
  $k_\h^{-1}$, this shows the required
  property~\eqref{equ:lattice-connectedness}. The number $N$ of steps
  in~\eqref{equ:lattice-connectedness} can be as large as
  $\mathcal{O}(1/\h)$, but the set of involved $\varkappa$'s in
  $\{-1,0,1\}^n$ being finite, we do get an $\h_0>0$ for which the
  result is uniform for all $\h\in\cI\cap\interval[open left]0{\h_0}$.

  If $\tilde B\subset B$ is a general connected open set, with compact
  closure contained in $B$, we can cover it by a finite number of
  deformed cubes of the form $G_0(\tilde V)\Subset\hat B$, as
  above. If $\h$ is small enough, the intersection of two such
  deformed cubes is either empty or contains a point in
  $\Lh$. Therefore the connectedness property holds for the union of
  two deformed cubes with a common point, and by induction for the
  union of all deformed cubes, hence for $\tilde B$.
\end{demo}
\begin{demo}[of Proposition~\ref{prop:linear-labelling}]
  Let $\h_0=\min_{\varkappa\in \{-1,0,1\}^n}\h_\varkappa$. We may also
  assume that $\h_0$ is small enough so that
  Lemma~\ref{lemm:lattice-connected} holds. Hence,
  applying~\eqref{equ:ktilde} to~\eqref{equ:lattice-connectedness} we
  obtain
  $\tilde k_\h(\lambda_{j+1}) = \tilde k_\h(\lambda_j) +
  \epsilon^{(j+1)}$, and hence
  \[
    \tilde k_\h(z_2) = \tilde k_\h(z_1) + \sum_{j=1}^N
    \epsilon^{(j)}\,.
  \]
  On the other hand, if we consider the good labelling $k_\h$, the
  translation invariance property~\eqref{equ:translation} gives
  similarly
  \[
    k_\h(z_2) = k_\h(z_1) + \sum_{j=1}^N \epsilon^{(j)}\,.
  \]
  Therefore,
  $\tilde k_\h(z_2)-k_\h(z_2) = \tilde k_\h(z_1)-k_\h(z_1)$, for all
  $z_1,z_2\in \tilde B$ and $\h\leq \h_0$; so $\tilde k_\h - k_\h$ is
  equal to a constant $K(\h)\in\ZM^n$ on $\tilde B$.  Hence
  $\tilde k_\h = k_\h - K(\h)$ is a linear labelling for $\Lh$ on
  $\tilde B$.
\end{demo}

\bigskip

Theorem~\ref{theo:bs} says that the joint spectrum of a quantum
integrable system is an asymptotic lattice near any regular point of
the momentum map; in other words, the joint spectrum possesses ``good
quantum numbers'' $(k_1,\dots,k_n)$, which is formalized in
Definition~\ref{defi:labelling} as the existence of a good labelling.
This property was used in~\cite{san-mono} to define the notion of
\emph{quantum monodromy}.

In this paper, we shall need to go one step further, namely we will
ask the question: can one construct a good labelling from the joint
spectrum only? Indeed, a good labelling is not unique, according to
the following, straightforward, result:
\begin{lemm}\label{lem:nonunique}
  If $G_\h$ is an asymptotic chart for $\Lh$, defined on an open set
  $U$, then for any orientation preserving linear transformation with
  integer coefficients $A\in\textup{SL}(n,\ZM)$, the map
  $\tilde G_\h:= G_\h\circ A$ is another asymptotic chart for $\Lh$,
  defined on $A^{-1}U$, and corresponding to the new good labelling
  $\tilde k_\h=A^{-1}\circ k_\h$.
\end{lemm}

\begin{prop}
  \label{prop:changement-lineaire}
  If $\bar k_\h^{(1)}$ and $\bar k_\h^{(2)}$ are two linear labellings
  for an asymptotic lattice $(\Lh,\cI, B)$, then for any connected
  open set $\tilde B\Subset B$, there exists a unique matrix
  $A\in\textup{SL}(n,\ZM)$, independent of $\h$, a family
  $(\varkappa_\h)_{\h\in\cI}$ in $\ZM^n$, and $\h_0>0$ such that
  \[
    \forall \h\in\cI\cap\interval0{\h_0}, \qquad \bar k_\h^{(2)} = A
    \circ \bar k_\h^{(1)} + \varkappa_\h\, \qquad \text{ on } \quad
    \Lh\cap \tilde B.
  \]
\end{prop}
\begin{demo}
  \noindent\textbf{(a)}~
  Let $G^{(1)}_\h$ and $G^{(2)}_\h$ be asymptotic charts associated
  with the linear labellings $\bar k_\h^{(1)}$ and $\bar k_\h^{(2)}$ ,
  respectively, and let $k^{(1)}_\h$ and $k^{(2)}_\h$ be the
  corresponding good labellings.  Let $\lambda_0\in\Lh$ and choose it
  as an ``origin'' of the linear labellings: up to changing their
  constant term, one can assume that
  $\bar k^{(1)}_\h(\lambda_0) = \bar k^{(2)}_\h(\lambda_0) =
  0\in\ZM^n$.  Let $\mathcal{B}=(e_1,\dots,e_n)$ be the canonical
  basis of $\ZM^n$, and let $(\lambda^{(j)}_1,\dots,\lambda^{(j)}_n)$,
  $j=1,2$, be the corresponding points in $\Lh$ for the two
  labellings, namely:
  \begin{equation}
    \label{equ:linear-basis}
    \bar k^{(j)}_\h(\lambda^{(j)}_i) = e_i.
  \end{equation}
  To make sure these points do exist, one can impose $\lambda_0$ to
  stay in a $\O(\h)$ neighborhood of a fixed point $c\in \tilde B$;
  then by Item~\ref{item:6} of Lemma~\ref{lemm:chart-pot-pourri},
  $(G^{(j)}_\h)^{-1}(\lambda_0)$ will stay in an $\O(\h)$ neighborhood
  of $\xi^{(j)}:=(G_0^{(j)})^{-1}(c)$. Hence,
  by~\eqref{equ:chart-labelling}, the same holds for
  $\h k^{(j)}_\h(\lambda_0)$. Denote
  $k^{(j)}_0 := k^{(j)}_\h(\lambda_0)$. Let
  $\tilde U^{(j)}_0\Subset (G^{(j)})^{-1}(B)$ containing $\xi^{(j)}$,
  and let $\h$ be small enough so that $\tilde U^{(j)}_0$, for
  $j=1,2$, contains $\h k^{(j)}_0 + \h\mathcal{B}$. By
  Item~\ref{item:def-asym-latt2b} of Definition~\ref{defi:labelling},
  the points $\lambda^{(j)}_i$ are now well-defined
  by~\eqref{equ:linear-basis}, because
  $\bar k^{(j)}_\h(\lambda^{(j)}_i) = \bar k^{(j)}_\h(\lambda^{(j)}_i)
  - \bar k^{(j)}_\h(\lambda_0) = k^{(j)}_\h(\lambda^{(j)}_i) -
  k^{(j)}_\h(\lambda_0)$, and hence $\lambda^{(j)}_i$ is defined by
  \begin{equation}
    \label{equ:linear-basis2}
    \h k^{(j)}_\h(\lambda^{(j)}_i) = \h k^{(j)}_0 + \h e_i \in \tilde
    U_0.
  \end{equation}

  \noindent\textbf{(b)}~ 
  Now, let $v^{(j)}_i=(\lambda^{(j)}_i - \lambda_0)$, $i=1,\dots,n$,
  be the corresponding ``basis vectors'' of the asymptotic lattice. By
  the uniform Taylor formula~\eqref{equ:G-taylor2},
  and~\eqref{equ:G0-bounded},
  \begin{equation}
    \label{equ:vi}
    v^{(1)}_i = \h (G_\h^{(1)})'(\h k^{(1)}_0) \cdot e_i + \O(\h^2) =
    \h (G_0^{(1)})'(\xi^{(1)}) \cdot e_i + \O(\h^2)\,.
  \end{equation}
  The point $\lambda_i^{(1)}$ can also be labelled by $k^{(2)}_\h$;
  let
  \begin{equation}
    \label{equ:zi}
    z_i := k^{(2)}_\h(\lambda^{(1)}_i) - k^{(2)}_0 \in \ZM^n.
  \end{equation}
  Remember that $z_i$, contrary to $e_i$, depends on $\h$. Similarly,
  we have
  \[
    v^{(1)}_i = \h (G_0^{(2)})'(\xi^{(2)}) \cdot z_i + \O(\h^2)\,,
  \]
  and hence
  \begin{equation}
    \label{equ:change-linear-basis}
    e_i = [(G_0^{(1)})'(\xi^{(1)})]^{-1} (G_0^{(2)})'(\xi^{(2)}) \cdot z_i + \O(\h)\,.
  \end{equation}
  Let $d_2 = d_2(\h):= \det (z_1,\dots,z_n)\in\NM$. Let
  $\delta_{1,2}:=\det \left([(G_0^{(1)})'(\xi^{(1)})]^{-1}
    (G_0^{(2)})'(\xi^{(2)})\right)$. From~\eqref{equ:change-linear-basis},
  $1 = \delta_{1,2} d_2 + \O(\h)$. Now repeat the argument with
  exchanging the two labellings: we obtain a new integer
  $d_1 = d_1(\h)\in\NM$ such that $1 = \delta_{2,1} d_1 + \O(\h)$,
  with $\delta_{2,1} = \delta_{1,2}^{-1}$. Hence
  $d_1 d_2 = 1 + \O(\h)$, which implies $d_1 d_2 = 1$ if $\h$ is small
  enough, and hence $d_2 = d_1 = 1$. This shows that
  $(z_1,\dots, z_n)$ is an oriented $\ZM$-basis of $\ZM^n$; let
  $A\in\textup{SL}(n,\ZM)$ be the corresponding
  matrix. Equation~\eqref{equ:change-linear-basis} gives
  $[(G_0^{(2)})'(\xi^{(2)})]^{-1} (G_0^{(1)})'(\xi^{(1)}) = A +
  \O(\h)$, and since the left-hand side does not depend on $\h$, $A$
  must converge towards it as $\h\to 0$, but since $A$ has integer
  coefficients, $A$ must be constant for $\h$ small enough, and we
  have
  \begin{equation}
    \label{equ:changement-A}
    A =  [(G_0^{(2)})'(\xi^{(2)})]^{-1} (G_0^{(1)})'(\xi^{(1)})
    \in \textup{SL}(n,\ZM)\,.
  \end{equation}
  From~\eqref{equ:zi} and~\eqref{equ:linear-basis2}, we have
  \[
    k^{(2)}_\h(\lambda^{(1)}_i) = A\cdot e_i + k^{(2)}_0 = A\cdot
    k^{(1)}_\h(\lambda^{(1)}_i) - A \cdot k^{(1)}_0 + k^{(2)}_0\,,
  \]
  which says that
  \begin{equation}
    \label{equ:compare-linear-labellings}
    k_\h^{(2)}(\lambda) = A\circ k_\h^{(1)}(\lambda) + \varkappa_h
  \end{equation}
  with $\varkappa_\h:=k^{(2)}_0 - A\cdot k^{(1)}_0$, \emph{when $\lambda$
    is restricted to the set}
  $\{\lambda^{(1)}_1, \dots \lambda^{(1)}_n\}$. It remains to extend
  this equality to the whole set $\Lh$.
  
  \noindent\textbf{(c)}~
  We cannot directly adapt the argument, replacing $\lambda^{(1)}_i$
  by \emph{any} $\lambda\in\Lh$, because the corresponding integral
  vector $z_i$ from~\eqref{equ:zi} would be unbounded (typically, of
  size $1/\h$), which would make the Taylor formula unusable. Instead,
  we need a connectedness argument, as in
  Proposition~\ref{prop:linear-labelling}.  First, remark that we may
  now replace in~\eqref{equ:linear-basis2} the vector $e_i$ by a
  linear combination $e:=\sum n_i e_i$, where $n_i$ are bounded
  integers, defining a lattice point $\mu$ by
  \[
    \h k^{(1)}_\h(\mu) = \h k^{(1)}_0 + \h e \in \tilde U_0\,.
  \]
  Let $v:=\mu - \lambda_0$. From
  \[
    v = \h (G_0^{(1)})'(\xi^{(2)}) \cdot e + \O(\h^2)\,,
  \]
  and comparing with~\eqref{equ:vi}, we see that
  $v = \sum n_i v^{(1)}_i + \O(\h)$. Letting
  $z:= k^{(2)}_\h(\mu) - k^{(2)}_0$, we see from
  \[
    v = \h (G_0^{(2)})'(\xi^{(2)}) \cdot z + \O(\h^2)\,,
  \]
  that $z = \sum n_i z_i + \O(\h)$. Since the integers $n_i$ are
  bounded, we must have $z = \sum n_i z_i$ for $\h$ small enough. This
  gives
  \[
    k^{(2)}_\h(\mu) = z + k^{(2)}_0 = A\cdot e + k^{(2)}_0\,.
  \]
  Therefore, the equality~\eqref{equ:compare-linear-labellings} still
  holds for $\mu$, \emph{i.e.} for all points of the form
  $\lambda = \lambda_0 \plus{k^{(1)}_\h} \vec\nu$, where $\vec\nu$ is
  uniformly bounded in $\ZM^n$. Let us restrict such $\vec\nu$ to
  belong to $\{-1,0,1\}^n$. This means that the set of points
  satisfying~\eqref{equ:compare-linear-labellings} is invariant under
  the action of $\{-1,0,1\}^n$. From
  Lemma~\ref{lemm:lattice-connected}, this set must be
  $\Lh \cap \tilde B$.
\end{demo}

\subsubsection{$\h$-continuity and drift for good labellings}

If we want to fully understand the structure of asymptotic lattices,
one should take into account the ``constant term'' $\varkappa_\h$ that
is disgarded in linear labellings. In other words, do we have an
analogue of Proposition~\ref{prop:changement-lineaire} for good
labellings? It turns out that answering this question leads to the
introduction of two new interesting concepts attached to any
asymptotic lattice $(\Lh,\cI, B)$: the $\h$-continuity and the
drift. The former asks whether we have enough values of $\h\in\cI$ to
continuously follow each individual point in the asymptotic lattice,
as $\h\to 0$. When this is the case, we see that, generically, these
points ``pass through'' the observation window $B$ (like rain drops
observed though a window) as if they wanted to converge to some
unknown limit which may lie outside of $B$: this is what the drift
will measure.

However, while this question is very natural, it is not necessary for
the treatment of the rotation number, because, as we will see below
(Section~\ref{ss:quantrot}), the latter is defined in terms of
differences of eigenvalues, which cancels the constant
term. Therefore, the reader interested in our main results about
rotation numbers can directly jump to Section~\ref{ss:quantrot}).

\medskip



\begin{defi}\label{defi:continuous}
  An asymptotic lattice is called \emph{$\h$-continuous} if the set
  $\cI$, which the small parameter $\h$ belongs to, is such that the
  set $\{\frac{1}{\h_1}-\frac{1}{\h_2}; \quad \h_1,\h_2\in\cI\}$
  accumulates at zero, namely:

  For all $\epsilon>0$, there exists $\h_0>0$ such that for all
  $\h_1\in\overline{\cI}\setminus\{0\}$ with $\h_1<\h_0$, there exists
  $\h_2\in\cI$, $\h_2<\h_1$ such that
  \begin{equation}
    \label{equ:h-cont}  
    \frac{1}{\h_2} - \frac{1}{\h_1} < \epsilon.
  \end{equation}
\end{defi}
Of course, this property is satisfied if $\cI=(0,\h_0]$, because the
map $\h\to\frac{1}{\h}$ is continuous. It also holds if the closure
$\overline{\cI}$ contains $[0,\h_0]$ for some $\h_0>0$. However, it is
not satisfied if $\cI=\{1/k; \quad k\in \NM^*\}$, which is typical in
geometric quantization. For $\cI=\{1/k^\beta; \quad k\in\NM^*\}$, the
property is satisfied if and only if $\beta\in\interval[open]01$.

The following proposition can be seen as the quantum analogue of
Lemma~\ref{lemm:ItoI'}.
\begin{prop}\label{prop:A}
  Let $k_\h$ and $\tilde k_\h$ be two good labellings for the
  $\h$-continuous asymptotic lattice $\Lh$. Then there exists a unique
  orientation preserving transformation
  $\tau\in \textup{GA}^+_\ZM(n,\ZM)$, independent of $\h$, and
  $\h_0>0$, such that $\tilde k_\h = \tau\circ k_\h$ for all
  $\h\leq\h_0$, $\h\in\cI$.
\end{prop}
Here $\textup{GA}^+_\ZM(n,\ZM)$ is the orientation preserving integral
affine group
$\textup{GA}^+_\ZM(n,\ZM):=\textup{SL}(n,\ZM)\ltimes\ZM^n$.

\medskip

\begin{demo}
  Let $G_\h$ and $\tilde G_\h$ be the corresponding asymptotic charts
  (Equation~\eqref{equ:chart}).  Let $\tilde F_\h$ be a formal inverse
  of $\tilde G_\h$, as in item~\ref{item:5} of
  Lemma~\ref{lemm:chart-pot-pourri}, and define
  $K_\h:=\tilde F_\h \circ G_\h$. Let $c\in B$, and let
  $\xi=G_0^{-1}(c)$.  Let $k=k(\h)$ be a family in $\ZM^n$ such that
  $\h k(\h) = \xi + \O(\h)$ as $\h\to 0$.  From~\eqref{equ:chart},
  there exists another family $\tilde k = \tilde k (\h)\in\ZM^n$ such
  that
  \begin{equation}
    K_\h(\h k) = \h \tilde k + \O(\h^\infty).
    \label{equ:K}
  \end{equation}
  Since $K_\h$ admits a $\Cinf$ asymptotic expansion
  $K_\h=K_0 + \h K_1 + \cdots$, a uniform Taylor expansion gives, for
  any fixed integer vector $v\in\ZM^n$,
  \[
    K_\h(\h (k + v)) - K_\h(\h k) = \h K_h'(\xi)\cdot v + \O(\h^2) =
    \h K_0'(\xi)\cdot v + \O(\h^2) \;.
  \]
  Since $\h (k + v) = \xi+\O(\h)$ as well, we apply
  again~\eqref{equ:K} (with a possibly different integral vector
  $\tilde k$) and obtain, letting $\h\to 0$,
  \[
    K_0'(\xi)\cdot v \in \ZM\,.
  \]
  This shows that the matrix $K'_0(\xi)$ has integer entries. Swapping
  $G_\h$ and $\tilde G_\h$, and repeating the argument, we obtain that
  the inverse of $K'_0(\xi)$ has integer entries as well, meaning that
  $K'_0(\xi)\in\textup{GL}(n,\ZM)$. In particular,
  $\xi\mapsto K'_0(\xi)$ must be a constant matrix
  $A\in \textup{GL}(n,\ZM)$ in the connected open set $G_0^{-1}(B)$,
  which in turn implies that $K_0(\xi) = A\xi + \alpha$, for some
  $\alpha\in\RM^n$. Thus, the maps $G_0$ and $\tilde G_0$ must differ
  by an affine transformation in $\textup{GA}(n,\ZM)$:
  \begin{equation*}
    \tilde G_0^{-1}\circ G_0 (\xi) = A\xi +\alpha.
  \end{equation*}
  Let us assume, by contradiction, that $\alpha\neq 0$. We have
  \begin{align}
    \label{equ:k}
    \tilde k = &~ \frac{1}{\h}K_0(\h k) + K_1(\h k) + \O(\h)\nonumber \\
    = &~ A k + \frac{\alpha}{\h} +   K_1(\xi) + R(\h),
  \end{align}
  with $R(\h)=\O(\h)$ (because $\h k = \xi + \O(\h)$).  Let
  $\ell(\h):=\tilde k(\h) - A k(\h)$. Let $\h_0\in\cI$ be small enough
  so that $\sup_{\h\leq\h_0} \norm{R(\h)}<\frac{1}{4}$ and so that the
  $\h$-continuity property~\eqref{equ:h-cont} holds for
  $\epsilon=\frac{1}{4\abs{\alpha}}$.  Let
  \[
    \cI_0:=\{\h\in \cI \quad | \quad \h\leq \h_0, \quad \ell(\h) =
    \ell(\h_0)\},
  \]
  and let $\h_1:=\inf\cI_0$.  Since $\ell(\h)\sim\frac{\alpha}{\h}$,
  it must be unbounded as $\h\to 0$, which implies that $\h_1>0$.

  \paragraph{Case 1: $\h_1\in\cI_0$.} Then, by $\h$-continuity, one
  can find $\h_2\in\cI$, $\h_2<\h_1$, such that
  $\frac{1}{\h_2} - \frac{1}{\h_1}\leq \frac{1}{2\abs{\alpha}}$.
  Since $\h_2\not\in\cI_0$ and $\ell(\h)\in\ZM^n$, we have
  \[
    \abs{\ell(\h_1) - \ell(\h_2)} \geq 1.
  \]
  From~\eqref{equ:k} we get
  \[
    \abs{\ell(\h_1) - \ell(\h_2)} < \abs{\alpha}\left(\frac{1}{\h_2} -
      \frac{1}{\h_1}\right) + \frac{1}{2} < 1,
  \]
  a contradiction.

\paragraph{Case 2: $\h_1\not\in\cI_0$.}
Let $\eta:=\frac{\h_1^2}{4\abs{\alpha}}$. By definition of $\h_1$,
there exists $\h_2\in\cI_0$ such that $\h_1<\h_2\leq \h_1+\eta$. We
have
\[
  \h_2\leq \h_1+\eta = \h_1(1+\frac{\h_1}{4\abs{\alpha}}) \leq
  \frac{\h_1}{1-\frac{\h_1}{4\abs{\alpha}}}.
\]
Hence $\frac{1}{\h_1} - \frac{1}{\h_2}\leq \frac{1}{4\abs{\alpha}}$.
Since $\h_1\in\overline{\cI}$, by $\h$-continuity one can find
$\h_3\in\cI$, $\h_3<\h_1$, such that
$\frac{1}{\h_3} - \frac{1}{\h_1}\leq \frac{1}{4\abs{\alpha}}$. Hence
$\frac{1}{\h_1} - \frac{1}{\h_3}\leq \frac{1}{2\abs{\alpha}}$ and we
may conclude as in Case 1 that
\[
  \abs{\ell(\h_1) - \ell(\h_3)} < \abs{\alpha}\left(\frac{1}{\h_1} -
    \frac{1}{\h_3}\right) + \frac{1}{2} < 1,
\]
while $\abs{\ell(\h_1) - \ell(\h_3)} \geq 1$, a contradiction.

Consequently we must have $\alpha=0$.  Thus we can write
\[
  \ell(\h) = K_1(\xi) + \O(\h),
\]
which implies that $\ell(\h)$ converges to $K_1(\xi)$ as $\h\to 0$ and
hence that $K_1(\xi)\in\ZM^n$, forcing it to be a constant
$\ell\in\ZM^n$ on the connected open set $G_0^{-1}(B)$. We obtain
\[
  \tilde k = Ak + \ell + R(\h),
\]
thus $R(\h)\in\ZM^n$ . Since $R(\h)=\O(\h)$, we must have, for
$\h<\h_0$ small enough, $R(\h)= 0$, which finishes the proof of the
proposition.
\end{demo}
Proposition~\ref{prop:A} shows that, given a set $\Lh\subset B$, the
set of good labellings on subsets $\tilde B\subset B$, \emph{i.e.} for
which $(\Lh\cap\tilde B, \cI, \tilde B)$ is an asymptotic lattice,
form a flat sheaf over $B$. As a consequence, if $B$ is simply
connected and $B$ is covered by open subsets on which $\Lh$ admit a
good labelling, then $\Lh$ admits a good labelling on $B$ itself.

The $\h$-continuity property has a natural interpretation in terms of
the uniform continuity of individual points in $\Lh$, as $\h$
varies. Recall from the discussion before
Definition~\ref{defi:linear-labelling} that the choice of a good
labelling defines an $\h$-evolution of each individual point in
$\Lh$. If $\Lh$ is $\h$-continuous, Proposition~\ref{prop:A} implies
that this evolution is in fact intrinsic to $\Lh$ itself. However,
this is not the case in general. The heuristics are very simple: given
$\h\in\cI$, if the closest element to $\h$ in $\cI$ is of the form
$\h + \delta$, then the lattice point corresponding to $G_\h(k\h)$
gets displaced by a distance of order $\O(k\delta)$. If $k\delta$ is
of the same order as the mean spacing between points (\emph{i.e.}
$\O(\h)$), then there will be a confusion between the evolution of
this point with that of its closest neighbors on $\Lh$. To avoid this
confusion, we need $k\delta \ll \h$; since $k$ is in general of order
$1/\h$, this means $\delta \ll \h^2$, which precisely gives the
$\h$-continuity condition~\eqref{equ:h-cont}.

Let us now give a precise statement of this, which will be needed for
the inverse problem in the next section.
\begin{prop}\label{prop:continuity}
  Let $(\Lh, \cI, B)$ be an asymptotic lattice. There exists
  $\h_0>0, \delta>0$, and $\epsilon>0$ such that, if
  $\h_1,\h_2\in\cI\cap\interval[open left]{0}{\h_0}$ satisfy:
  \[
    \abs{\frac{1}{\h_1} - \frac{1}{\h_2}} < \epsilon
  \]
  then the following holds. Fix $\lambda_1\in\mathcal{L}_{\h_1}$ and
  let $\lambda_2\in\mathcal{L}_{\h_2}$ be defined by
  \[
    \lambda_2 = k_{\h_2}^{-1}(k_{\h_1}(\lambda_1)),
  \]
  for some good labelling $k_\h$. Then
  \begin{equation}
    \mathcal{L}_{\h_i} \cap B(\lambda_j,\delta\h_i) = \{\lambda_i\} \qquad
    \forall (i,j)\in\{1,2\}^2.
    \label{equ:unique-in-ball}  
  \end{equation}
  Here $B(\lambda_j,\delta\h_i)$ denotes the Euclidean ball centered
  at $\lambda_j$, of radius $\delta\h_i$.
\end{prop}
In other words, we fix $\lambda_1\in\mathcal{L}_{\h_1}$, and we
consider the evolution of $\Lh$ as $\h$ moves from $\h_1$ to $\h_2$;
then the closest element to $\lambda_1$ in $\mathcal{L}_{\h_2}$ is
unique and is precisely the natural evolution of $\lambda_1$ obtained
by fixing its integer label $k_{\h_1}(\lambda_1)$.

\begin{demo}
  First of all, by Remark~\ref{rema:N=1} we may choose $\delta_0>0$
  such that $\Lh\cap B(\lambda,\delta_0\h) = \{\lambda\}$ for all
  $\h<\h_0$ and all $\lambda\in\Lh$. In particular,
  ~\eqref{equ:unique-in-ball} holds when $i=j$ for any
  $\delta\leq \delta_0$.  Without loss of generality we may assume
  $\h_2< \h_1$. Let $(G_\h, U)$ be an asymptotic chart associated with
  the good labelling $k_\h$. Let $k_1 = k_{\h_1}(\lambda_1)$. We have
  \[
    \lambda_2 = G_{\h_2}(\h_2 k_1) + \O(\h_2^\infty) = G_{\h_2}(\h_2
    k_1) + \O(\h_1^\infty)\,,
  \]
  and hence
  $\lambda_2 - \lambda_1 = G_{\h_2}(\h_2 k_1) - G_{\h_1}(\h_1 k_1) +
  \O(\h_1^\infty)$.  Let $\tilde U\Subset U$ be such that
  $G_0^{-1}(\tilde U)$ contains $\overline{B}$. Taking $\h_0$ small
  enough, we may assume that $G_\h$ is invertible on $\tilde U$ for
  all $\h\in\cI\cap\interval[open left]{0}{\h_0}$, see
  Lemma~\ref{lemm:chart-pot-pourri}, Item~\ref{item:6}. Therefore
  $\h_1 k_1\in B$ and hence is bounded. Thus, there exists $M>0$ such
  that
  \begin{equation*}
    \norm{\h_2 k_1 - \h_1 k_1} \leq \epsilon \h_1 \h_2 k_1 \leq \epsilon M
    \h_2.
  \end{equation*}
  Hence we may apply Item~\ref{item:4} of
  Lemma~\ref{lemm:chart-pot-pourri} and get
  $\norm{G_{\h_2}(\h_2 k_1) - G_{\h_2}(\h_1 k_1)} \leq L_0\epsilon M
  \h_2$. Using now Item~\ref{item:3} of the same lemma, we obtain
  \begin{align}
    \norm{\lambda_2 - \lambda_1} 
    & \leq L_0\epsilon M \h_2 + C(\abs{\h_2 -
      \h_1} + \h_1^2 + \h_2^2) + \O(\h_1^\infty)\nonumber\\
    &  \leq L_0\epsilon M \h_2 +
      C(\epsilon+2)\h_1^2  + \O(\h_1^\infty).
      \label{equ:l1-l2}  
  \end{align}
  Since $\h_1\leq \h_2(1+\epsilon\h_1)$, we see that, if $\epsilon$
  and $\h_0$ are small enough, the right-hand side
  of~\eqref{equ:l1-l2} is less than $\delta_0\h_2/3$. Choosing finally
  $\delta = \delta_0/2$, we see that~\eqref{equ:unique-in-ball} must
  hold. Indeed, if $\mu_1\in\mathcal{L}_{\h_1}$ and
  $\mu_1\neq \lambda_1$, then
  $\norm{\mu_1 - \lambda_1}\geq \delta_0 \h_1$, hence
  $\norm{\mu_1 - \lambda_2} \geq \norm{\mu_1 - \lambda_1} -
  \norm{\lambda_2 - \lambda_1} \geq \delta_0\h_1 - \delta_0\h_2/3 >
  \delta \h_1$, and similarly if $\mu_2\in\mathcal{L}_{h_2}$ and
  $\mu_2\neq \lambda_2$, then
  $\norm{\mu_2 - \lambda_1} \geq \delta_0 \h_2 - \delta_0\h_2/3 >
  \delta \h_2$.
\end{demo}

We conclude this section by making precise the relationship between
the natural evolution of $\h$-continuous asymptotic lattices given by
Proposition~\ref{prop:continuity} above, the ``drift'' mentioned
before Definition~\ref{defi:linear-labelling}, and the asymptotic
expansion of the lattice points in $\h$.

\begin{prop}
  \label{prop:asymptotic-exp}
  Let $(\Lh, \cI, B)$ be an asymptotic lattice. Let $c\in B$, and let
  $\lambda_\h\in\Lh$ be such that $ \lambda_\h = c + \O(\h)$.  Let
  $k_\h$ be a good labelling for this lattice.  There there exists
  $\alpha_c\in\RM^n$, $A_c\in\textup{M}_n(\RM)$, and $\beta_c\in\RM^n$
  such that
  \begin{equation}
    \lambda_\h = \alpha_c + \h(A_c\cdot k_\h(\lambda_\h) + \beta_c) +
    \O(\h^2).
    \label{equ:asymptotic-exp}
  \end{equation}
  Moreover,
  \[
    \alpha_c = c - G_0'(\xi_0)\cdot \xi_0, \quad A_c = G_0'(\xi_0),
    \quad \beta_c = G_1(\xi_0),
  \]
  where $G_\h = G_0 + \h G_1 + \cdots$ is the asymptotic chart
  associated with $k_\h$, and $\xi_0 = G_0^{-1}(c)$.
\end{prop}
\begin{demo}
  Let $k = k_\h(\lambda_\h)$ for simplifying notation. Using a uniform
  Taylor formula~\eqref{equ:G-taylor2}, we have
  \[
    G_\h(\h k) = G_0(\xi_0) + G'_0(\xi_0)\cdot (\h k - \xi_0) + \h
    G_1(\xi_0) + \O(\norm{\h k - \xi_0}^2) + \h \O(\norm{\h k -
      \xi_0}) + \O(\h^2).
  \]
  Using the boundedness of $G_\h^{-1}$ (Item~\ref{item:5} or
  \ref{item:6} of Lemma~\ref{lemm:chart-pot-pourri}), we have
  \[
    \norm{\h k - \xi_0} = \O(\norm{G_\h^{-1}(\h k) - c}).
  \]
  Since $\lambda_\h = G_\h(\h k) + \O(\h^\infty)$, this gives
  \begin{equation}
    \lambda_\h = \alpha_c + \h(A_c\cdot k_\h(\lambda_\h) + \beta_c
    + \O(\norm{\lambda_\h - c})) + \O(\norm{\lambda_\h - c}^2) + 
    \O(\h^2),
    \label{equ:asymptotic-exp2}
  \end{equation}
  which establishes the result.
\end{demo}

Assume now that the asymptotic lattice $\Lh$ is $\h$-continuous; we
see from \eqref{equ:asymptotic-exp} that the natural evolution of a
point $\lambda\in\Lh$, as $\h$ varies slightly, which corresponds to
freezing the integer $k$, is to move ``as if it wanted to tend to
$\alpha_c$''; and, in general $\alpha_c \neq c$.
\begin{defi}
  \label{defi:drift}
  Let $(\Lh, \cI, B)$ be an $\h$-continuous asymptotic lattice. Let
  $c\in B$ and $\xi_0 := G_0^{-1}(c)$. We call the quantity
  \[
    \delta_c := G_0'(\xi_0)\cdot \xi_0 \in\RM^n
  \]
  the \emph{drift} of the asymptotic lattice at $c$. It does not
  depend on the choice of a good labelling.
\end{defi}
\begin{demo}
  Let show that the drift is indeed well defined. Let $\tilde G_\h$ be
  another asymptotic chart of $\Lh$. Since $\Lh$ is $\h$-continuous,
  we may apply Proposition~\ref{prop:A}, and obtain
  $\tau\in \textup{GA}^+_\ZM(n,\ZM)$ such that
  \[
    G_\h(\h k_\h) = \tilde G_\h(\h\tau(k_\h)) + \O(\h^\infty)
  \]
  for any indices $k_\h$ such that $\h k_\h$ belongs to $G_0^{-1}(B)$.
  Whence (with $\tau (k) = Ak + \ell$, \emph{i.e.} $A=\tau'$)
  \[
    \begin{aligned}
      G_0(\h k_\h) &= \tilde G_0(\h\tau(k_\h)) + \O(\h) \\
      & = \tilde G_0(\h(A k_\h +\ell)) + \O(\h) \\
      & = \tilde G_0(\h A k_\h ) + \O(\h) \\
      & = \tilde G_0\circ A (\h k_\h) + \O(\h)\,.
    \end{aligned}
  \]
  We deduce that $G_0(\xi)=\tilde G_0\circ A (\xi)$ for any
  $\xi\in G_0^{-1}(B)$, which means $G_0=\tilde G_0\circ A$. Let
  $\tilde \xi_0 = \tilde G_0^{-1}(c)$, we have $\tilde\xi_0=A(\xi_0)$
  and $G_0'(\xi_0)\cdot v = \tilde G_0'(A(\xi_0)) \cdot A (v)$, hence
  $G_0'(\xi_0)\cdot \xi_0 = \tilde G_0'(\tilde\xi_0)) \cdot
  \tilde\xi_0$.
\end{demo}

\begin{rema}\label{rema:euler}
  When the asymptotic lattice $(\Lh, \cI, B)$ is the intersection of
  $B$ with a joint spectrum of a quantum integrable system of
  commuting pseudodifferential operators, we may take
  $\cI=\interval[open left]{0}{\h_0}$, which is $\h$-continous, and
  the drift is well-defined. In this case, the drift has a natural
  geometric definition. Since $M=T^*X$ is a cotangent bundle, there is
  a globally defined primitive for the symplectic form, namely the
  Liouville 1-form $\alpha$. The formula
  $I_j(\Lambda) = \frac{1}{2\pi}\int_{\gamma_j}\alpha$, where
  $(\gamma_1,\dots,\gamma_n)$ is a basis of loops on the Liouville
  torus $\Lambda$, defines action variables that depend only on this
  choice of basis. Hence, the integral affine manifold $\textup{B}_r$
  is ``exact'' in the sense that its structure group can be reduced to
  the linear group $\textup{GL}(n,\ZM)$ (instead of the affine
  group). (This fact was used in~\cite{san-mono} to show that quantum
  monodromy is linear instead of affine.)  Therefore, the tangent
  bundle $T\textup{B}_r$ has a natural Euler vector field
  $E=\sum_j I_j\deriv{}{I_j}$. Now, consider the momentum map
  $F:M\to\RM^n$. By the action-angle theorem, we have a smooth local
  diffeomorphism $G_0$ such that $F=G_0(I)$; in other words,
  $G_0:\textup{B}_r\to\RM^n$ is the natural map induced by $F$ on the
  leaf space. The push-forward of $E$ by this map is
  $\dd{}_I G_0\cdot E = G_0'(I)\cdot I$, \emph{i.e.} the drift of
  Definition~\ref{defi:drift}.
\end{rema}
Since the drift is well-defined, we can expect to find a way to
recover it directly from the asymptotic lattice, without choosing a
particular good labelling. This is the content of the following lemma.
\begin{prop}
  \label{prop:drift}
  Let $(\Lh, \cI, B)$ be an $\h$-continuous asymptotic lattice.  Let
  $c\in B$, and let $\lambda_\h\in\Lh$ be such that
  $ \lambda_\h = c + \O(\h)$.  For $\h_1,\h_2\in\cI$, let
  $\lambda(\h_1,\h_2)$ be a closest element to $\lambda_{\h_1}$ in
  $\mathcal{L}_{\h_2}$.  We define the divided difference
  \begin{equation}
    \label{eq:defdelta}
    \delta(\h_1,\h_2) :=
    \h_1\,\frac{\lambda_{\h_1} - \lambda(\h_1,\h_2)}{\h_1-\h_2}\,.
  \end{equation}
  For any $N>2$, there exists positive numbers $\h_0$,
  $\varepsilon_0$, and $C>0$ such that for all
  $\varepsilon\le\varepsilon_0$, $\h_1\in\cI$, $\h_1\le\h_0$ and
  $\h_2\in\cI$ such that
  \begin{equation}
    \label{equ:good-h1-h2}
    \h_2<\h_1\,, \quad \h_1 - \h_2 \geq \h_1^N\,, \quad \text{ and }
    \quad \frac{1}{\h_2} - \frac{1}{\h_1} < \varepsilon\,,
  \end{equation}
  then $\lambda(\h_1,\h_2)$ is unique, and
  \[
    \abs{\delta(\h_1,\h_2) - \delta_c} \leq C \h_1\,,
  \]
  where $\delta_c$ is the drift of $\Lh$ at $c$.
\end{prop}
As we shall see below (Lemma~\ref{lemm:h-continuous-refined}), for an
$\h$-continuous asymptotic lattice, it is always possible to find
couples $(\h_1,\h_2)$ satisfying the
requirements~\eqref{equ:good-h1-h2} while being arbitrarily small.
Hence, the conclusion of this lemma is that we can recover the drift
$\delta_c$ from the limit of the divided
difference~\eqref{eq:defdelta}, as $\h_1\to 0$.

\begin{demo}[of Proposition~\ref{prop:drift}]
  Let $k_\h$ be a good labelling for $\Lh$, and $G_\h$ an associated
  asymptotic chart.  From Proposition~\ref{prop:continuity} there
  exist $\delta>0$, $\varepsilon_0>0$ and $\h_0>0$ such that for all
  $\h_1,\h_2\in\cI\cap\interval[open left]{0}{\h_0}$ with
  $\abs{h_1^{-1}-\h_2^{-1}}<\varepsilon_0$, the element
  $\lambda(\h_1,\h_2)$ is uniquely defined, and
  $\norm{\lambda(\h_1,\h_2) - \lambda_{\h_1}}<\delta\h_2$. Moreover,
  $k_{\h_1}(\lambda_{\h_1}) = k_{\h_2}(\lambda (\h_1,\h_2))$. For ease
  of notation, let us denote by $k$ this integer. We have
  \begin{equation}
    \lambda_{\h_1} = G_{\h_1}(\h_1 k) + \O(\h_1^\infty), \qquad \lambda
    (\h_1,\h_2) = G_{\h_2}(\h_2 k) + \O(\h_2^\infty).
    \label{equ:lambda-h1h2}
  \end{equation}
  If follows from Item~\ref{item:2} of
  Lemma~\ref{lemm:chart-pot-pourri} that $c = G_0(\h_1 k) + \O(\h_1)$,
  and hence, because $\h_2\in(\frac{\h_1}{1+\h_1\epsilon}, \h_1)$,
  there exists a constant $C>0$ such that
  \begin{equation*}
    \norm{\h_1 k - \xi_0} \leq C \h_1, \quad
    \text{ and } \norm{\h_2 k - \xi_0} \leq C \h_1,
  \end{equation*}
  where $\xi_0 = G_0^{-1}(c)$.  On the other hand, from
  Item~\ref{item:3} of Lemma~\ref{lemm:chart-pot-pourri}, there exists
  $C_N>0$ such that
  \begin{equation*}
    \|G_{\h_1}(\h_1 k) - G_{\h_2}(\h_1 k)\| \leq
    C_N (|\h_1-\h_2| + \h_1^N + \h_2^N). 
  \end{equation*}
  (Note that we could have got rid of the terms $\h_1^N + \h_2^N$ by
  choosing an asymptotic chart that is smooth in $\h$, see
  Remark~\ref{rema:h-smooth}).  Hence, since $\h_1-\h_2 \geq \h_1^N$,
  \begin{equation}
    \label{equ:G12}
    \h_1 \frac{\|G_{\h_1}(\h_1 k) - G_{\h_2}(\h_1 k)\|}{\h_1 - \h_2}
    \leq 3C_N \h_1
  \end{equation}
  From Item~\ref{item:4} of Lemma~\ref{lemm:chart-pot-pourri}, we can
  write
  \begin{equation*}
    G_{\h_2}(\h_1 k) - G_{\h_2}(\h_2 k) =
    (\h_1 - \h_2)G'_{\h_2}(\h_2 k)\cdot(k)
    + \O(\abs{\h_1- \h_2}^2 \norm{k}^2)\,,
  \end{equation*}
  which, in view of~\eqref{equ:G12}, and using that
  $\h_j k = \xi_0 + \O(\h_1)$ for $j=1,2$, we obtain
  \begin{align*}
    \h_1 \frac{G_{\h_1}(\h_1 k) - G_{\h_2}(\h_2 k)}{\h_1 - \h_2}
    & = G'_{\h_2}(\h_2 k)\cdot(\h_1 k) + \O(\h_1 \abs{\h_1- \h_2} \norm{k}^2)\\
    ~& = G'_{\h_2}(\xi_0)\cdot(\xi_0) + \O(\h_1)
       + \O(\h_1 \epsilon\h_1\h_2 \norm{k}^2)\\
    ~& = \delta_c + \O(\h_2) + \O(\h_1) + \O(\epsilon\h_2\norm{\xi_0}^2)\\
    ~& = \delta_c + \O(\h_1)\,.
  \end{align*}
  Finally, using~\eqref{equ:lambda-h1h2} and the fact that
  $\h_1-\h_2 \geq \h_1^N$, this yields the desired estimate.
\end{demo}

The following lemma shows that the requirement~\eqref{equ:good-h1-h2}
can always be met.
\begin{lemm}
  \label{lemm:h-continuous-refined}
  Assume $\cI$ is $\h$-continuous. For any $\epsilon>0$, for any
  $N>2$, there exists $\h_0>0$ such that for all
  $\h_1\in\cI\cap\interval[open left]{0}{\h_0}$, there exists
  $\h_2\in\cI$ such that
  \[
    \h_1 < \h_2, \quad \frac{1}{\h_2} - \frac{1}{\h_1} < \epsilon
    \quad \text{ and } \h_1 - \h_2 \geq \h^N\,.
  \]
\end{lemm}
\begin{demo}
  By contradiction, assume that the statement of the lemma does not
  hold.  Let
  $\h_3 := \inf \{\h\in\cI; \h<\h_1, \frac{1}{\h} - \frac{1}{\h_1} <
  \epsilon\}$. By Definition~\ref{defi:continuous} the set in question
  is not empty, hence $\h_3\ge 0$, and
  $\h_3\in\interval{\frac{\h_1}{1+\h_1\epsilon}}{\h_1}$. By the
  negation of the lemma, one must have $\h_1 - \h_3 < \h_1^N$, and
  hence $\frac{1}{\h_3}-\frac{1}{\h_1} < \h_1^{N-2}(1+\h_0\epsilon)$.
  Applying Definition~\ref{defi:continuous} (with $\epsilon$ replaced
  by $\frac{\epsilon}{2}$) to this
  $\h_3\in\overline{\cI}\setminus\{0\}$, we obtain $\h_4\in\cI$ with
  $\h_4<\h_3$ and
  $\frac{1}{\h_4} - \frac{1}{\h_3}<\frac{\epsilon}{2}$. Hence
  $\frac{1}{\h_4}- \frac{1}{\h_1} < \frac{\epsilon}{2} +
  \h_0^{N-2}(1+\h_0\epsilon)$. If $\h_0$ was small enough, the
  right-hand side is less than $\epsilon$, thus $\h_4$ contradicts the
  definition of $\h_3$.
\end{demo}

\subsection{Quantum rotation number}
\label{ss:quantrot}
In this section, we define a spectral quantity from the joint spectrum
of a quantum integrable system, which will be the natural analogue of
the rotation number $w_I(\Lambda)$. As we shall see, making this
quantity a purely spectral invariant is not obvious, because, in the
same way as the classical rotation number $w_I(\Lambda)$ depended on
the choice of the action variables $I$
(Definition~\ref{defi:rotation}), the quantum rotation number will
depend on the choice of a good labelling.

Let us consider a quantum Hamiltonian $\hat H$, and assume that it is
completely integrable, in the sense that there exists a quantum
integrable system $(P_1, \dots, P_n)$ (see Section~\ref{sec:bs}) with
$P_n=\hat H$, with proper classical momentum map
$F=(f_1,\dots,f_{n-1}, f_n=H):M \to \RM^n$ (here $f_j$ is the
principal symbol of $P_j$).  Let $c\in F(M)$ be a regular value of
$F$, and assume that $F^{-1}(c)$ is connected.  From
Theorem~\ref{theo:bs-bis}, we know that the joint spectrum $\Sigma_\h$
is an asymptotic lattice near $c$.
\begin{defi}
  \label{defi:quant-rotat-numb}
  Let $\lambda\mapsto k_\h(\lambda)$ be a good labelling for the joint
  spectrum $\Sigma_\h$ in a ball $B\subset \RM^n$. Let
  $\tilde B\Subset B$ be an open set.  Let
  $\lambda\subset \Sigma\cap \tilde B$, and let
  $k=k_\h(\lambda)\in\ZM^n$ be the corresponding label.  Denote
  \[
    \lambda =: (\lambda^{(1)}_k(\h),\dots, \lambda^{(n)}_k(\h))\,,
  \]
  and let $E_k(\h):=\lambda^{(n)}_k(\h)$ be the corresponding
  eigenvalues of $\hat H$.  We define the \emph{quantum rotation
    number}, for the good labelling $k_\h$, to be
  \begin{equation}
    \label{equ:def-q-rot}
    [ \hat w_\h ](k) := [ E_{k+e_1}(\h) - E_{k}(\h) : \cdots :
    E_{k+e_n}(\h) - E_k(\h) ] \in \RP^{n-1}\,,
  \end{equation}
  where $(e_1,\dots,e_n)$ is the canonical basis of $\ZM^n$.
\end{defi}
The restriction to the smaller open subset $\tilde B$ ensures that the
labels $(k+e_j)$ do correspond to joint eigenvalues in $B$, when
$\h < \h_0$ and $\h_0$ is small enough.  When $n=2$, as in the classical case
(Equation~\eqref{equ:ratio}), it is often more convenient to think of
the rotation number as an element of the 1-point compactification
$\overline\RM$,
\begin{equation}
  \hat w_\h(k) := \frac{E_{k+e_1}(\h) - E_{k}(\h)}{E_{k+e_2}(\h) -
    E_{k}(\h)} \in\overline{\RM}.
  \label{equ:quantum-rotation-number-real}
\end{equation}
The following result shows that, once a good labelling is known, the
classical rotation number can be recovered from the quantum rotation
number, in the semiclassical limit.
\begin{theo}
  \label{theo:q-rotation-number}
  Let $c\in \RM^n$ be a regular value of $F$, and assume that the
  Liouville torus $\Lambda:=F^{-1}(c)$ is connected.  Let $\lambda_\h$
  be a joint eigenvalue in $\Sigma_\h$ such that
  \[
    \lambda_\h = c + \O(\h).
  \]
  Let $G_\h$ be an asymptotic chart for $\Sigma_\h$ in a neighborhood
  of $c$, let $k=k(\h)\in\ZM^n$ be the corresponding label for
  $\lambda_\h$ and let $I:=G_0^{-1}\circ F$ be the associated action
  variables (see Equation~\eqref{equ:G0}). Then
  \[
    [\hat w_\h](k) = [w_I](\Lambda) + \O(\h).
  \]
\end{theo}
The term $\O(\h)$ is relative to the topology of $\RP^{n-1}$, and is
uniform if $c$ varies in a compact subset. The existence of
$\lambda_\h$ in the statement of the theorem is guaranteed by
Lemma~\ref{lemm:lattice-dense}.

\medskip

\begin{demo}
  Let $g_\h$ be the last component of $G_\h$, so that
  $E_{k}(\h) = g_\h(\h k) + \O(\h^\infty)$. We have the
  $\Cinf$ asymptotic expansion
  \[
    g_h(\xi) = g_0(\xi) + \h g_1(\xi) + \cdots
  \]
  where $g_0$ is the last component of $G_0$; thus, $H=g_0(I)$. By
  definition
  \[
    [w_I](\Lambda) = [\partial_1 g_0(I(\Lambda)) : \cdots :
    \partial_n g_0(I(\Lambda))]\,.
  \]
  A Taylor formula (Item~\ref{item:4} of
  Lemma~\ref{lemm:chart-pot-pourri}) gives
  \begin{align*}
    E_{k+e_j}(\h) - E_{k}(\h)
    & = g_\h(\h(k+e_j)) - g_\h(\h k) + \O(\h^\infty)\\
    & = \h\partial_j g_\h(\h k)  + \O(\h^2)\\
    & = \h\partial_j g_0(\h k)  + \O(\h^2),
  \end{align*}
  because $\h k$ is bounded. In fact, 
  $\h k =  G_0^{-1}(c) +\O(\h)= I(\Lambda) + \O(\h)$; hence
  $E_{k+e_j}(\h) - E_{k}(\h) = \h\partial_j g_0(I(\Lambda)) +
  \O(\h^2)$, which yields
  \[
    [\hat w_\h](k) = [\partial_1 g_0(I(\Lambda)) + \O(\h) : \cdots  :
    \partial_n
    g_0(I(\Lambda)) + \O(\h)],
  \]
  which gives the result, because
  $(\partial_1 g_0(I(\Lambda)), \dots, \partial_n g_0(I(\Lambda)))
  \neq 0$.
\end{demo}

\subsection{Quantum rotation number for semitoric systems}
\label{sec:q-semitoric}

If $F=(J,H)$ is a classical semitoric system, we have shown in
Section~\ref{sec:semitoric-case} that the rotation number is well
defined as an angle in $\RM/\ZM$. We will show the quantum analogue
here, namely that for a quantum semitoric system, the rotation number
can be defined with no ambiguity on the choice of a good labelling.

\begin{defi}[\cite{san-alvaro-survey}]
  A quantum integrable system $(\hat J, \hat H)$ is called semitoric
  if the corresponding classical system $(J,H)$ given by the principal
  symbols is semitoric.
\end{defi}
Examples include Laplace-Beltrami operators on surfaces of revolution,
Schrödinger operators with axi-symmetric potential (like the quantum
spherical pendulum), or toric
Laplacians~\cite{dryden-guillemin-senadias12}. On the Toeplitz side
(which we do not consider explicitly here, see
Remark~\ref{rema:toeplitz}), many spin-coupling systems (or
Jaynes-Cummings) that are used in quantum optic or chemistry are
semitoric systems, see for
instance~\cite{sadovski-zhilinski,cushman-prl,lefloch-pelayo19} and
references therein.
\begin{rema}
  It would be very interesting to have a purely spectral
  characterization of a quantum semitoric system. The quantum analogue
  of the Hamiltonian $\SM^1$ action should be reflected in the fact
  that the spectrum of $\hat J$ coming from a bounded region of the
  joint spectrum is close to an arithmetic sequence of the type
  $\alpha+\h(j+\mu)$, $j\in\ZM$ (see Proposition~\ref{prop:J-spectrum}
  below).

  However the fact that a point $c\in\RM^2$ is a regular value of $F$
  seems more delicate to obtain in a purely spectral way. Moreover,
  the semitoric hypothesis also impacts the singularity types at the
  boundary of the image of $F=(J,H)$ (see
  Proposition~\ref{prop:elliptic-labelling}). We don't address these
  issues here. Instead, we assume the semitoric nature of the system,
  and from this we try and recover the rotation number.
\end{rema}

\begin{figure}
  \centering \includegraphics[width=0.6\linewidth]{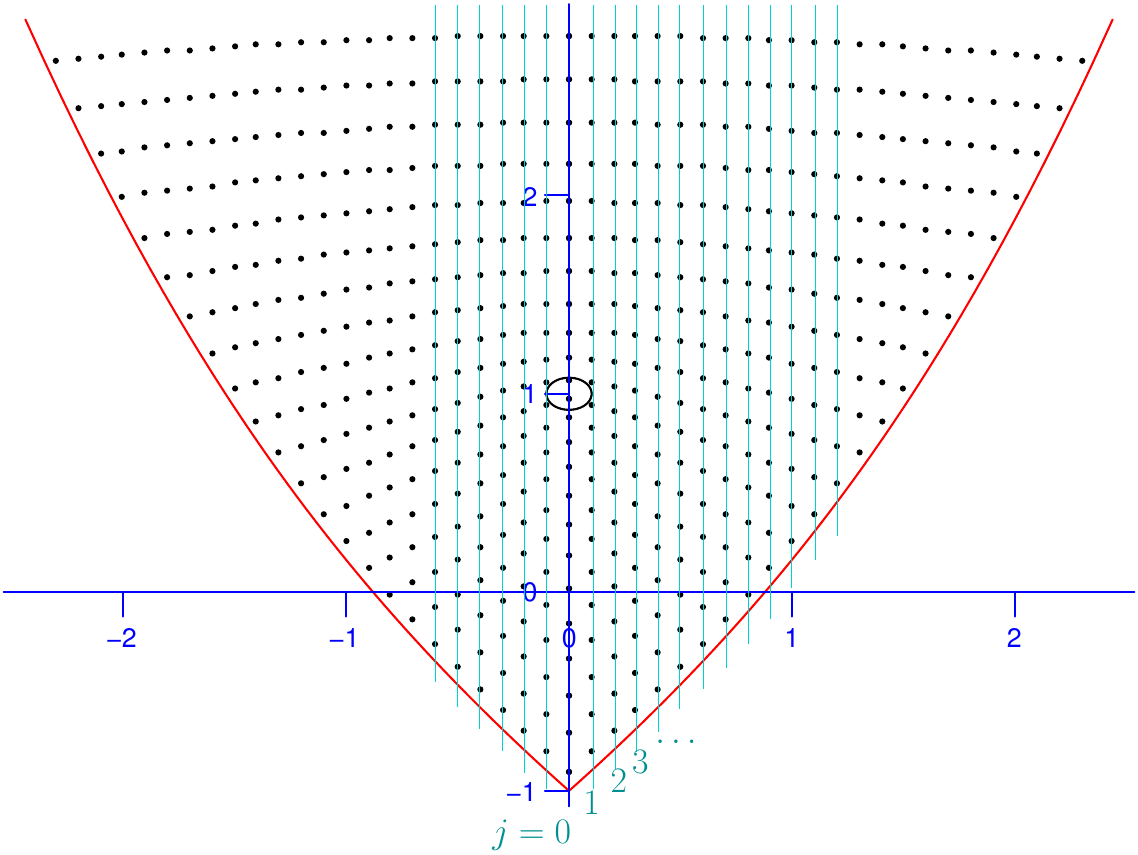}
  \caption{Joint spectrum of the Spherical Pendulum. Joint eigenvalues
    are organized along vertical lines indexed by $j$.}
  \label{fig:pendulum}
\end{figure}

In the semitoric case, Proposition~\ref{prop:asymptotic-exp} can be
improved.
\begin{prop}
  \label{prop:J-spectrum}
  Let $(\hat J, \hat H)$ be a semitoric quantum integrable system
  defined for $\h\in\interval[open left]{0}{\h_0}$.  Let $c\in\RM^2$
  be a regular value of $F$ and assume that $F^{-1}(c)$ is connected,
  where $F=(J,H)$ is the classical momentum map. Then there exist a
  ball $B$ around $c$, $\alpha\in\RM$, $\mu\in \RM$, and a good
  labelling $\lambda\mapsto (j,k)$ of the joint spectrum $\Sigma_\h$
  in $B$ such that, uniformly for
  $\lambda= (J_{j,k}(\h), E_{j,k}(\h))\in \Sigma_\h\cap B$,
  \begin{equation}
    J_{j,k}(\h) = \alpha+ \h(j+\mu + \O(\lambda-c)) + \O(\h^2) .
    \label{equ:J-spectrum}
  \end{equation}
\end{prop}
Loosely speaking, this proposition says that, in a small ball around
$c$, the joint spectrum of a semitoric system is organized along
regularly spaced vertical lines $J=\alpha+\h(j+\mu)$, and the quantum
number $j\in\ZM$ labels these lines. See Figure~\ref{fig:pendulum}.
If $J$ is proper (and $\hat J$ has no subprincipal symbol), the result
follows from the cluster structure of the spectrum of
pseudodifferential operators with periodic characteristics, see for
instance~\cite{DG,dozias,colin-bica,hall-hitrik-sjostrand}. In our
case, we do not impose the properness of $J$ and there is no
restriction on subprincipal symbols; the result is still valid because
we restrict the joint spectrum to a small $B$ (it would not hold for
the usual spectrum of $\hat J$ alone).

\begin{lemm}\label{lemm:semitoric-chart}
  Under the assumptions of Proposition~\ref{prop:J-spectrum}, let $B$
  be a bounded, simply connected open subset of regular values around
  $c$, such that $(\Sigma_\h,\interval[open left]0{\h_0}, B)$ is an
  asymptotic lattice. Then this asymptotic lattice admits a
  \emph{semitoric asymptotic chart}, \emph{i.e.} an asymptotic chart
  $G_\h\sim G_0 + \h G_1 + \cdots$ such that the first component
  $G_0^{(1)}$ of
  $G_0 = (G_0^{(1)}(\xi_1,\xi_2), G_0^{(1)}(\xi_1,\xi_2)) : U \to
  \RM^2$ satisfies, for all $(\xi_1,\xi_2)\in U$,
  \begin{equation}
    \dd G_0^{(1)}(\xi_1,\xi_2) = \dd \xi_1.
    \label{equ:dG01}
  \end{equation}
\end{lemm}
\begin{demo}
  Let $G_\h$ be an asymptotic chart for $\Sigma_\h$ near $c$,
  see~\eqref{equ:bs} and Definition~\ref{defi:lattice}, and let
  $I=(I_1,I_2)$ be the corresponding action coordinates:
  $F=G_0\circ I$.  Since $(J,H)$ is semitoric, there exist oriented
  action coordinates near $\Lambda:=F^{-1}(c)$ of the form
  $(J,J_2)$. Hence by Lemma \ref{lemm:ItoI'} there is an affine
  transformation $\tau\in \textup{GA}(n,\ZM)$ such that
  $(I_1,I_2) = \tau(J,J_2) = A(J,J_2)+z$, where
  $A\in\textup{SL}(2,\ZM)$ and $z\in\RM^2$. Then the map
  $\hat G_\h(\xi):= G_\h(A\xi)$ is an asymptotic chart by Lemma
  \ref{lem:nonunique} and we have:
  \[
    \dd {\hat G_0} = \dd G_0 \circ A = \dd G_0 \circ \dd \tau = \dd
    (G_0 \circ \tau)\,.
  \]
  Since $G_0(\tau(J,J_2)) = (J,H)$, this implies
  $\dd(G_0\circ \tau)^{(1)} = \dd\xi_1$, and gives the result.
\end{demo}
Equation~\eqref{equ:dG01} means that, up to a constant, the chart
$G_\h$ is associated with semitoric action variables
(Lemma~\ref{lemm:semitoric-action}).

\medskip

\begin{demo}[of Proposition~\ref{prop:J-spectrum}]
  Let us consider the proof of
  Proposition~\ref{prop:asymptotic-exp}. In view of
  Lemma~\ref{lemm:semitoric-chart}, if we
  project~\eqref{equ:asymptotic-exp2} on the first component, the term
  $\O(\norm{\lambda_\h - c}^2)$ disappears, because the first
  component of $G_0$ is an affine map. This gives the required
  estimate.
\end{demo}

\begin{figure}[ht]
  \centering
  \includegraphics[width=0.5\linewidth]{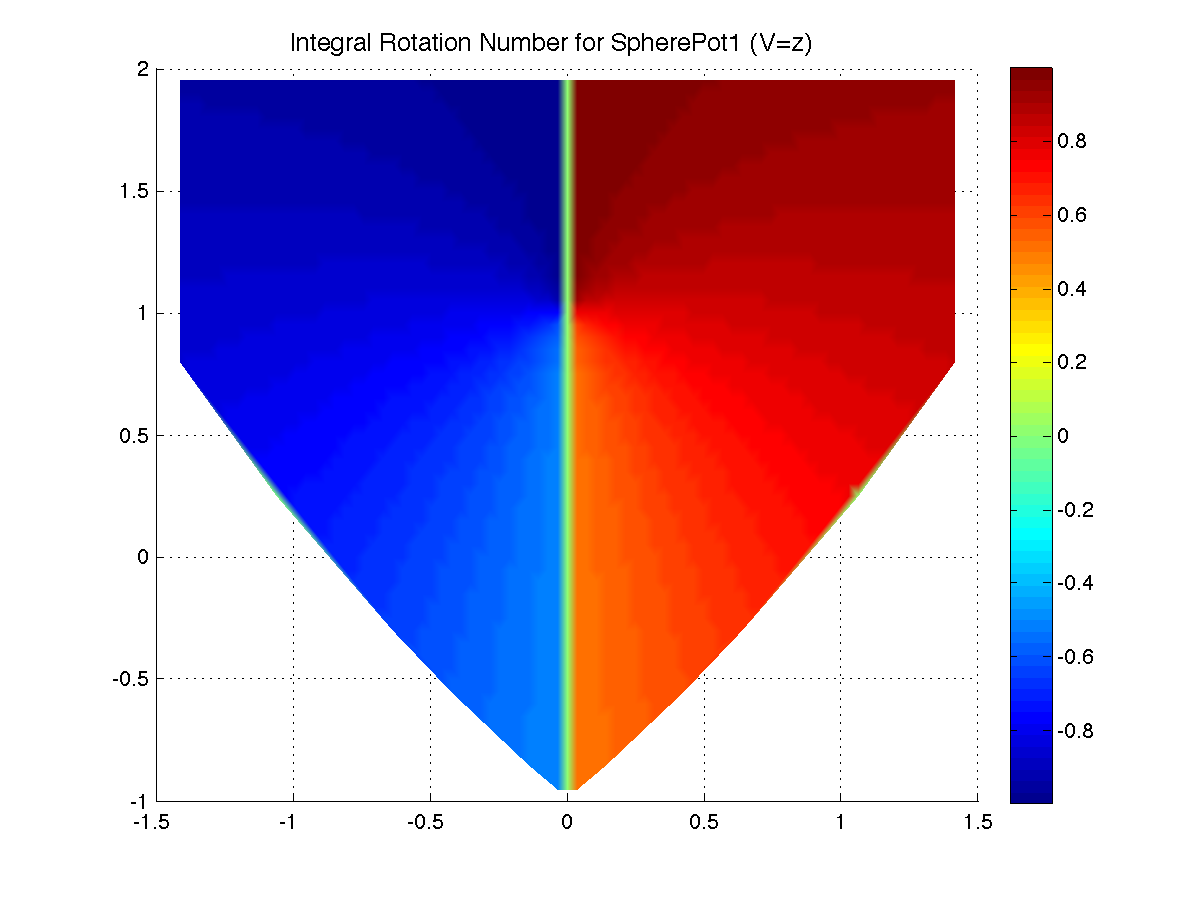}%
  \includegraphics[width=0.5\linewidth]{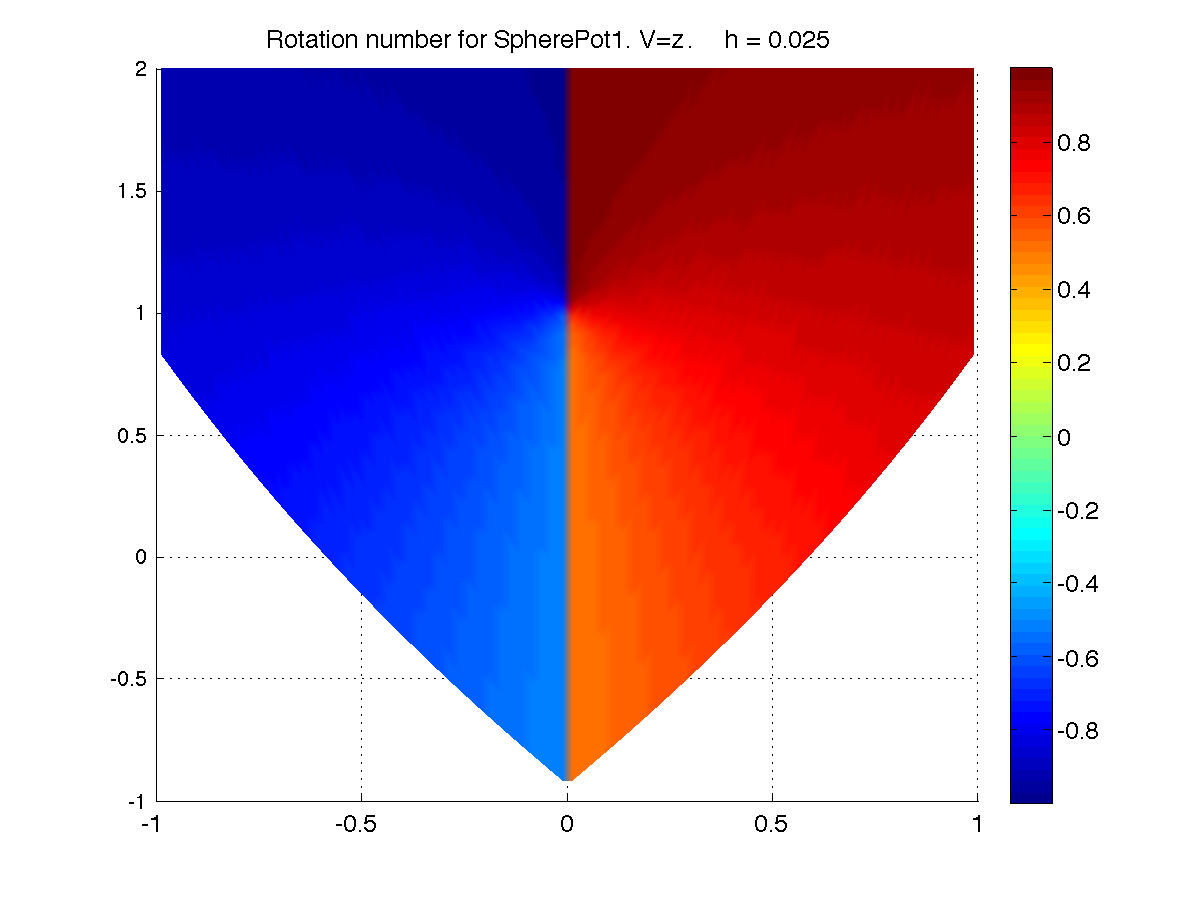}
  \caption{%
    Classical (left) and quantum (right) rotation numbers for an
    axisymmetric Schrödinger operator on the sphere $\mathbb{S}^2$
    with potential $V=z$ \emph{(quantum spherical pendulum)}.}
  \label{fig:numerics}
\end{figure}

\begin{rema}
  Comparing with Proposition~\ref{prop:asymptotic-exp}, we see that
  the number $\alpha$ in~\eqref{eq:J1} is equal to
  $c^{(1)} - \xi_0^{(1)}$, \emph{i.e.} the first component of
  $c-\xi_0$. Thus $\xi_0^{(1)}$ is the first component of the
  \emph{drift} of the joint spectrum
  (Definition~\ref{defi:drift}). Recall that $\xi_0^{(1)}$ is also the
  value of the action integral along the $\mathbb{S}^1$-cycle on the
  torus $\Lambda_c$.  This is consistent with Remark~\ref{rema:euler},
  since $J$ is an action variable and can be completed into a set of
  local affine coordinates $(J=I_1,I_2)$ of $\textup{B}_r$.
\end{rema}
\begin{theo}
  \label{theo:q-semitoric}
  Let $(\hat J, \hat H)$ be a semitoric quantum integrable system,
  with momentum map $F=(J,H)$.  Let $c\in \RM^2$ be a regular value of
  $F$, and assume that the Liouville torus $\Lambda:=F^{-1}(c)$ is
  connected.  Let $G_\h=G_0+\h G_1+\cdots $ be an asymptotic chart for
  $\Sigma_\h$ near $c$, such that the associated good labelling
  $\lambda\mapsto (j,k)$ satisfies Equation~\eqref{equ:J-spectrum}.
  Then~\eqref{equ:dG01} holds, \emph{i.e.} the first component of
  $G_0$ satisfies, for all $(\xi_1,\xi_2)$ near $G_0^{-1}(c)$:
  \begin{equation*}
    \dd G_0^{(1)}(\xi_1,\xi_2) = \dd \xi_1.
  \end{equation*}
  Moreover, let $\lambda_\h$ be a joint eigenvalue in $\Sigma_\h$ such
  that
  \[
    \lambda_\h = c + \O(\h)\,.
  \]
  Then the corresponding quantum rotation number
  (Equation~\eqref{equ:quantum-rotation-number-real}) satisfies:
  \[
    \hat w_\h(j,k) = w(\Lambda) + \O(\h) \mod \ZM\,,
  \]
  where $w(\Lambda)$ is the semitoric rotation number in the sense of
  Definition~\ref{defi:semitoric-rotation}.
\end{theo}
\begin{demo}
  By assumption, the first component of $\lambda_\h$ has the
  asymptotic expansion:
  \begin{equation}
    \label{eq:J1}
    J_{j,k}(\h) =  \alpha+ \h(j+\mu + \O(\lambda-c)) + \O(\h^2) \,,
  \end{equation}
  where $\alpha,\mu\in\RM$ do not depend on $\h$.  From
  Theorem~\ref{theo:q-rotation-number}, we have
  \begin{equation}
    \label{equ:semitoric-w-proj}
    [\hat w_\h](j,k) = [w_I](\Lambda) + \O(\h),
  \end{equation}
  where $I=G^{-1}_0\circ F$. From Lemma~\ref{lemm:semitoric-chart}, we
  introduce a semitoric asymptotic chart $\hat G_\h:=G_\h\circ A$,
  where $A \in \textup{SL}(2,\ZM)$ is such that $\hat G_\h$
  satisfies \begin{equation}
    \label{equ:dG01-hat}
    \dd{\hat G_0}^{(1)}(\xi_1,\xi_2) = \dd \xi_1.
  \end{equation}
  From the proof of Proposition~\ref{prop:J-spectrum}, $\hat G_\h$
  satisfies Equation~\eqref{equ:J-spectrum}: there exists
  $\alpha',\mu'$ in $\RM$ such that, $\forall \h<\h_0$,
  $J_{j,k}(\h) =\alpha' + \h(j'+\mu' + \O(\lambda-c)) + \O(\h^2)$;
  here $j'=aj + bk$, where $a$, $b$ are the integers forming the fist
  line of $A^{-1}$.  Comparing with Equation~\eqref{eq:J1}, we obtain,
  since $\lambda_\h - c = \O(\h)$,
  \begin{equation*}
    \alpha' +  \h(aj+bk+\mu') = \alpha + \h(j+\mu ) +
    \O(\h^2)\,,
  \end{equation*}
  hence there exists a constant $C>0$ such that, for all $\h\in\cI$,
  \begin{equation}
    \abs{\alpha'-\alpha + \h\big((a-1)j + b k + \mu' - \mu\big)}\leq C \h^2.
    \label{equ:two-semitoric-labellings}
  \end{equation}
  This holds for all $\h$-dependent couples $(j,k) = (j(\h), k(\h))$
  that label a joint eigenvalue $\lambda_\h\in\Sigma_\h$ near $c$. If
  $u,v$ are given integers, independent of $\h$, then the joint
  eigenvalue labelled by $(j(\h)+u,k(\h)+v)$ differs from $\lambda_\h$
  by $\O(\h)$. Hence~\eqref{equ:two-semitoric-labellings} must hold
  for this new label as well, and writing the triangle inequality we
  obtain, for all fixed $(u,v)\in\ZM^2$,
  \[
    \h\abs{(a-1)u + b v} \leq 2 C_{u,v} \h^2.
  \]
  Choosing $(u,v)=(0,1)$ or $(1,0)$, as soon as $2C_{u,v}\h<1$, this
  implies $a=1$ and $b=0$ (and hence $\alpha = \alpha'$ and
  $\mu' = \mu$ from~\eqref{equ:two-semitoric-labellings}). This
  entails that the matrix $A$ takes the form $A=
  \begin{pmatrix}
    1 & 0\\r & 1
  \end{pmatrix}$.  Hence Equation~\eqref{equ:dG01} follows
  from~\eqref{equ:dG01-hat}.  This proves that, up to a constant, the
  action variables $I$ are in fact semitoric. Thus, from
  Proposition~\ref{prop:semitoric} we obtain
  \[
    w(\Lambda) = w_I(\Lambda) \mod \ZM,
  \]
  which, together with~\eqref{equ:semitoric-w-proj} finishes the proof
  of the Theorem (recall that in the semitoric case, the direction
  $[w_I(\Lambda)]$ can never be vertical).
\end{demo}
In order to illustrate this result, we have produced a numerical
comparison between classical and quantum rotation numbers, in the case
of the spherical pendulum (an axisymmetric Schrödinger operator on the
sphere $\mathbb{S}^2$), see Figure~\ref{fig:numerics}.

Theorem~\ref{theo:q-semitoric} will be used in
Section~\ref{sec:semitoric-algorithm} to obtain, algorithmically, a
good labelling for the joint spectrum of a quantum semitoric system,
based on the fact that the ``global label $j$'' corresponding to the
integer in~\eqref{equ:J-spectrum} is easy to detect. However, even for
a semitoric system, assigning another quantum number $k$ to each joint
eigenvalue can be more delicate. While it is in principle always
possible near any regular value $c$ of the moment map $(J,H)$, thanks
to Theorem~\ref{theo:bs-bis}, there is no global recipe: due to a
possibly non-trivial monodromy, a good labelling for the joint
spectrum may simply not globally exist~\cite{san-mono}. There is,
however, a common situation where a natural second quantum number
shows up: consider, for a fixed $j$, the corresponding vertical
spectral band $V_j$. If the spectrum inside $V_j$ is either bounded
from above or from below, and if this bound corresponds to an elliptic
singularity --- which is generically the case ---, then labelling $k$
in increasing order from the bounded side will provide a good
labelling. To understand this, let $c=(x,y)$ be a rank-one elliptic
singularity of $F$, where $x$ is a regular value of $J$. This means
that the restriction of $H$ to the submanifold $J^{-1}(x)$ admits a
Morse-Bott non-degenerate critical point $m$ with $H(m)=y$. Let us
assume that $F^{-1}(c)$ is connected; by the theory of non-degenerate
singularities of integrable systems~\cite{zung-I}, this fiber must
then be a circle of critical points (the corresponding point in the
reduced space $J^{-1}(x)/\SM^1$ is a standard Morse non-degenerate
critical point). Letting $x$ vary in a small neighborhood, we thus
obtain a smooth cylinder of critical points in $M$, whose critical
values form a smooth curve through $c$ in $\RM^2$, which belongs to
the boundary of $F(M)$. This situation is called a \emph{simple
  $J$-transversal elliptic singularity}. From the viewpoint of the
energy $H$, there are two situations, where $c$ is either a local
maximum or minimum of $H$ restricted to $J^{-1}(x)$. For simplicity we
shall only deal with the minimum case (which we call `positive' in the
statement below). This is the case of the Spherical Pendulum,
Figure~\ref{fig:pendulum}; of course the maximum case is completely
similar.
\begin{prop}
  \label{prop:elliptic-labelling}
  Let $(\hat J, \hat H)$ be a semitoric quantum integrable system with
  principal symbols $(J,H)$. Let $c\in\RM^2$ be a simple
  \emph{positive} $J$-transversal elliptic critical value of
  $F$. Then, for $\h$ small enough, the joint eigenvalues
  $\lambda\in\Sigma_\h$ in a neighborhood of $c$ belong to the union
  of disjoint vertical bands $V_j$ given by the equation
  \begin{equation}
    \label{equ:elliptic-labelling}
    x = \alpha+ \h(j+\mu + \O(\lambda-c)) + \O(\h^2), \qquad j\in\ZM
  \end{equation}
  where $\alpha,\mu\in\RM$ are fixed.  In each vertical band $V_j$,
  the $y$ coordinates of the joint eigenvalues are distinct and
  bounded from below. We label them in increasing $y$-order by a
  non-negative integer $k\in\NM$. Then, for any regular value $c'$ of
  $F$, close to $c$, the labels $(j,k)$ form a good labelling of
  $\Sigma_\h$ near $c'$.
\end{prop}
Contrary to the rest of the article, for this result the standard
action-angle theorem (and its semiclassical version) is not enough. We
need to resort to the microlocal study of the spectrum near a simple
transversally elliptic singularity, which was done in~\cite[Theorem
5.2.4]{san-panoramas}.
\begin{theo}[\cite{san-panoramas}]
  \label{theo:bs-elliptic}
  Let $(\hat J, \hat H)$ be a quantum integrable system, with momentum
  map $F=(J,H)$, and let $c$ be a simple transversally elliptic
  critical value of $F$. Then the joint spectrum $\Sigma_\h$ in a
  neighborhood of $c$ (independent of $\h$) can be described as
  follows:
  \begin{enumerate}
  \item joint eigenvalues have multiplicity one, in the sense of item
    1. of Theorem~\ref{theo:bs};
  \item $\Sigma_\h$ is an ``asymptotic half lattice'' in the sense
    that item 2. of Theorem~\ref{theo:bs} holds:
    \begin{equation}
      \lambda = G_\h(\h k_1, \h k_2) + \O(\h^\infty),
      \label{equ:bs-elliptic}      
    \end{equation}
    when replacing $(k_1,k_2)\in\ZM^2$ by $(k_1,k_2)\in\ZM\times\NM$,
    and replacing Equation~\eqref{equ:G0} by
    \begin{equation}
      \label{equ:G0-elliptic}
      F = G_0(\xi_1,q_2)
    \end{equation}
    where $q_2(x,\xi) = (x_2^2+\xi_2^2)/2$. Here the local coordinates
    $(x_1,\xi_1,x_2,\xi_2)\in T^*\SM^1\times \RM^2$ near
    $\SM^1\times\{\beta\}\times (0,0)$ for some $\beta\in\RM$,
    describing a neighborhood of the circle $F^{-1}(c)$, are
    symplectic, and $G_0$ is a local diffeomorphism from
    $(\RM^2,(\beta,0))$ to a neighborhood of $\overline B$.
  \end{enumerate}
\end{theo}
\begin{demo}[of Proposition~\ref{prop:elliptic-labelling}]
  It follows from~\eqref{equ:G0-elliptic} that $J=g(\xi_1,q_2)$ for a
  smooth $g$. Since $J$, $\xi_1$ and $q_2$ all have $2\pi$-periodic
  flows, there must exist integers $a,b$ such that
  $dJ=a d\xi_1 + b dq_2$. Since the $J$-action is effective, $a$ and
  $b$ must be co-prime. The hypothesis of $J$-transversality implies
  that $a\neq 0$. Thus $\xi_1=a^{-1}J - ba^{-1} q_2 + \textup{const}$,
  and the same argument implies that $a^{-1}\in\ZM$ and hence
  $a=\pm 1$. Up to composing by the symplectomorphism
  $(x_1,\xi_1)\mapsto (-\xi_1,x_1)$ (and replacing $k_1$ by $-k_1$
  in~\eqref{equ:bs-elliptic}), one may assume that $a=1$. Thus
  $(\xi_1,q_2) = \tau(J,q_2)$, where $\tau\in \textup{GA}(n,\ZM)$ with
  linear part $A=
  \begin{pmatrix}
    1 & -b\\0 & 1
  \end{pmatrix}$.  Arguing as in the proof of
  Proposition~\ref{prop:J-spectrum}, we let
  $\hat G_\h(x,y):= G_\h(A(x,y))$, and the joint spectrum is described
  as
  \[
    \lambda = G_\h(\h k_1, \h k_2) + \O(\h^\infty) = \hat
    G_\h(\h(j,k_2)) + \O(\h^\infty),
  \]
  with $k_1=j-b k_2$. Then
  $\hat G_\h = \hat G_0 + \h \hat G_1 + \cdots$, with
  \[
    \hat G_0 (x,y)= G_0\circ A (x,y)= G_0 (x - b y + \alpha, y),
  \]
  for some constant $\alpha\in\RM$. But since the first component of
  $G_0$ is $g$ and $J=g(\xi_1,q_2) = g(J-b q_2,q_2)$, we deduce that
  for all $(x,y)$ near $A^{-1}(\beta,0) = (\beta,0)$,
  \[
    \hat G_0(x,y) = (x+\alpha,f(x,y))
  \]
  for some smooth function $f$ with $\partial_y f\neq 0$.  We conclude
  as in the proof of Proposition~\ref{prop:J-spectrum} that the first
  component of $\lambda$ takes the form
  \[
    J_{j,k_2}(\h) = \alpha + \h j + \h \mu + \h \O(\lambda-c) +
    \O(\h^2).
  \]
  Hence we obtain the description of the vertical bands $B_{j}$
  in~\eqref{equ:elliptic-labelling}, which are disjoint if
  $\abs{\lambda-c}$ is small enough (independently of $\h$) and $\h$
  itself is small enough. Finally, if $j$ is fixed, the second
  component of the joint eigenvalues is given by
  \begin{equation*}
    E_{j,k_2}(\h) = f(\h j, \h k_2) + \h \O(\lambda-c) + \O(\h^2).
  \end{equation*}
  Thus,
  $ E_{j,k_2+1} - E_{j,k_2} = \h \partial_y f(c) + \h \O(\lambda-c) +
  \O(\h^2)$.  For a fixed value of $J$, $H$ is assumed to be minimal
  at the elliptic critical point. Since $q_2\geq 0$, we must have
  $\partial_{q_2}H \geq 0$ when $q_2=0$, which implies
  $\partial_y f >0$. Hence $ E_{j,k_2+1} - E_{j,k_2} \geq \h/C$
  for some constant $C>0$.

  It remains to prove that $(j,k_2)$ is a good labelling away from the
  critical value. Since $(k_1,k_2)=A(j,k_2)$, this is equivalent to
  showing that $(k_1,k_2)$ is a good labelling. Since item \emph{1.}
  in Definition~\ref{defi:lattice} is already given by
  Theorem~\ref{theo:bs-elliptic}, we must only prove item \emph{2.}
  Let $c'$ be a regular value of $F$. Let $B$ be a ball around $c$ in
  which the description of the joint spectrum by
  Theorem~\ref{theo:bs-elliptic} is valid:
  \[
    \lambda = G_\h(\h k_1, \h k_2) + \O(\h^\infty),
  \]
  where $(k_1,k_2)\in\ZM\times\NM$, and $G_\h$ is defined in a
  neighborhood of $(\beta,0)$. Let $B'\subset B$ be a ball of regular
  values around $c'$. We claim that ${G_\h}{\restr_{G_0^{-1} (B')}}$
  is an asymptotic chart for $\Sigma_\h$ in $B'$. Indeed, let
  $(\beta',\gamma'):= G_0^{(-1)}(c')$; necessarily $\gamma'>0$ because
  $c'$ is regular. Now assume that $(k_1,k_2)\in\ZM^2$ is such that
  $\h (k_1,k_2)\in G_0^{-1}(B')$; then $\h k_2$ is close to $\gamma'$,
  and in particular, $k_2>0$. Hence, by~\eqref{equ:bs-elliptic}, there
  exists a joint eigenvalue in $\lambda\in B'$ such that
  $ \lambda = G_\h(\h k_1, \h k_2) + \O(\h^\infty)$. Therefore,
  according to Definition~\ref{defi:lattice} (item \emph{2.}), $G_\h$
  is an asymptotic chart for $\Sigma_\h$ in $B'$.
\end{demo}

\subsection{Labelling algorithms}
\label{sec:labelling-algorithm}

In the inverse spectral problem, we are given a (portion of a) joint
spectrum $\Sigma_\h$; we know that $\Sigma_\h$ is an asymptotic
lattice, but we don't have access, \emph{a priori}, to any asymptotic
chart. In Theorem~\ref{theo:q-rotation-number}, we showed that the
classical rotation number can be obtained from the joint spectrum,
\emph{provided} eigenvalues are properly labelled. Hence, in order to
recover the rotation number from $\Sigma_\h$ itself, as a bare
discrete set of points, we have to find a way to construct an
appropriate labelling. In view of formula~\eqref{equ:def-q-rot}, a
\emph{linear labelling} (Definition~\ref{defi:linear-labelling}) is
enough for this purpose.

In this section, we show how to obtain, in an algorithmic and robust
way, a linear labelling from a bare 2-dimensional asymptotic lattice  $(\Lh, \cI, B)$.

We proceed in two steps: first, we construct an algorithm that, for
any fixed value of $\h$, assigns integer labels to the points of
$\Lh$. However, these labels may not have the required continuity with
respect to $\h$ that would entitle them to be linear labellings.  To
remedy this, in a second step, we exhibit another algorithm that,
given any sequence $(\h_i)$ in $\cI$ tending to $0$, correct the
previously obtained labels to make them a true linear labelling.

In the semitoric case, the simple structure of the joint spectrum,
organized in vertical bands (Theorem~\ref{theo:q-semitoric}), makes it
easier. In particular, the second step is not necessary, see
Section~\ref{sec:semitoric-algorithm}.

We wrote our algorithms in the two-dimensional case, which corresponds
to our initial motivations. We believe that the general case could be
dealt with in a similar way for most parts, although some steps like
Theorem~\ref{theo:basis} would require a different approach.

\subsubsection{An algorithm for fixed $\h$}
\label{sec:an-algorithm-fixed}

The core part of our labelling algorithm can be executed for any fixed
value of $\h\in\cI$.  Actually, it can be applied to any finite set
$\Lh\subset B$. If $c\in B$, we will use the expression ``choose a
closest point to $c$'' to indicate that we have to choose a point in
$\Lh$ that minimizes the distance to $c$. After such a point has been
picked up and labelled, we remove it from $\Lh$; thus, subsequent
calls to ``choose a closest point to...'' implicitly indicate that
this point should be chosen within $\Lh$ minus all the already
labelled points. This is of course always possible as long as this new
set is not empty.

\begin{figure}
  \centering \includegraphics[width=0.3\linewidth]{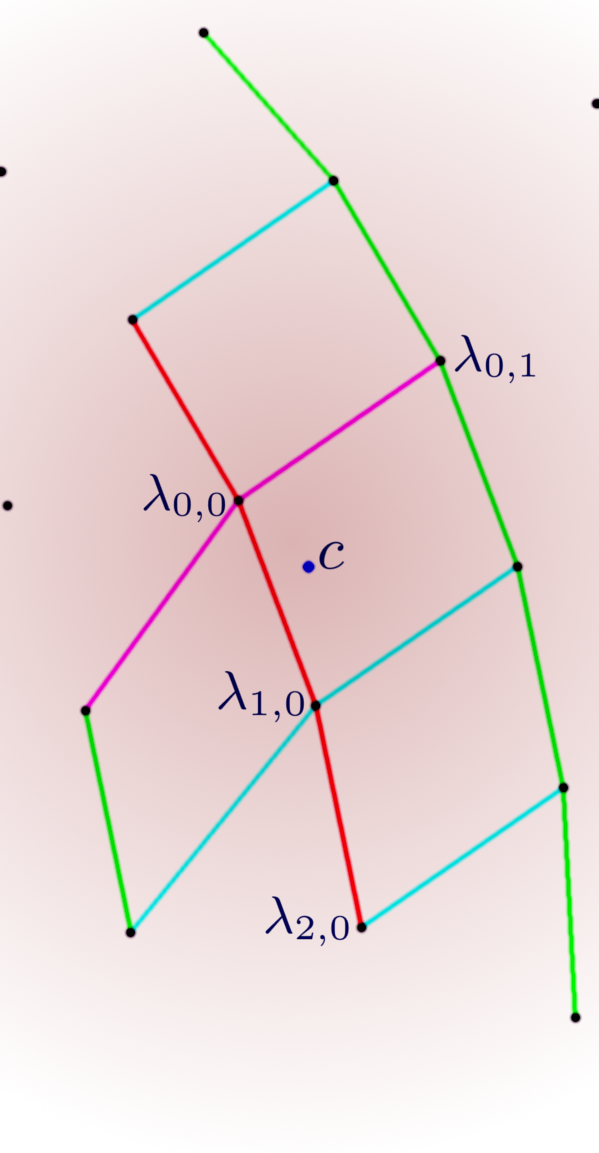}
  \caption{The labelling algorithm}
  \label{fig:algo}
\end{figure}

If $(n,m)$ is a label for a point $\lambda\in\Lh$, we shall denote
this point $\lambda=\lambda_{n,m}$. See Figure~\ref{fig:algo}. The
complete algorithm consists in the twelve following steps.
\begin{enumerate}
\item \label{item:algo-c} Choose an open subset $B_0\Subset B$, and
  fix $c\in B_0$.
\item
  Choose a closest point to $c$. Label it as $(0,0)$.

\item \label{item:algo-1-0} Choose a closest point to $\lambda_{0,0}$
  (in the set $\Lh\setminus \{\lambda_{0,0}\}$).  Label it as $(1,0)$.

\item \label{item:algo-n-0}

  Choose a closest point to
  $2\lambda_{1,0} - \lambda_{0,0} =\lambda_{1,0} +
  (\lambda_{1,0}-\lambda_{0,0})$ and label it as $(2,0)$.

  Continuing in this fashion (if $\lambda_{n-1,0}$ is chosen, take
  $\lambda_{n,0}$ to be a closest point to
  $\lambda_{n-1,0} +(\lambda_{n-1,0} - \lambda_{n-2,0})$), label
  points $\lambda_{n,0}$, $n>0$, until the next point lies outside of
  $B_0$.

  Label $\lambda_{n,0}$ for negative $n$ in the same way, starting
  from the closest point to $2\lambda_{0,0} - \lambda_{1,0}$.

\item \label{item:algo-0-1} Choose a closest point to $\lambda_{0,0}$
  not already labeled and label it as $(0,1)$.

\item \label{item:algo-1-1} Label a closest point to
  $\lambda_{0,1} + (\lambda_{1,0} - \lambda_{0,0})$ as $(1,1)$.

\item Use the points $\lambda_{0,1}$, $\lambda_{1,1}$ to repeat the
  process in step 4, labelling as many points $\lambda_{n,1}$ as
  possible (if $\lambda_{n-1,1}$ is chosen, take $\lambda_{n,1}$ to be
  a closest point to
  $\lambda_{n-1,1} +(\lambda_{n-1,1} - \lambda_{n-2,1})$).

\item \label{item:algo-0-2} Label a closest point to
  $2\lambda_{0,1} - \lambda_{0,0}$ as $\lambda_{0,2}$. Repeat steps
  6-7 to label all points $\lambda_{n,2}$.

\item Continuing as above, label all points $\lambda_{n,m}$, $m>0$
  which lie in the given neighborhood.

\item Label a closest point to $2\lambda_{0,0} - \lambda_{0,1}$ as
  $(0,-1)$.
\item Repeat steps 6,7,8,9 with negative $m$ indices.

\item\label{item:algo-determinant} Finally, if the determinant of the
  vectors
  $(\lambda_{1,0} - \lambda_{0,0},\lambda_{0,1} - \lambda_{0,0})$ is
  negative, switch the labelling $\lambda_{n,m} \mapsto \lambda_{n,-m}$
  (in order to make it oriented).
\end{enumerate}

Let us prove now some properties of this algorithm when
$(\Lh, \cI, B)$ is a given asymptotic lattice.  We use the notation of
Definition~\ref{defi:lattice}; in particular, $G_\h: U \to \RM^n$ is
an asymptotic chart, and $k_h$ is the corresponding good labelling. Of
course, both of them are \emph{a priori} unknown. Let $B_0\Subset B$
be an open subset containing $c$, and let $\tilde B_0$ be an open
subset such that $B_0\Subset \tilde B_0 \Subset B$. Let $\tilde\h_0$
be given by Item~\ref{item:def-asym-latt2b} of
Definition~\ref{defi:lattice} with
$\tilde U_0:= G_0^{-1}(\tilde B_0)$. Because the results of this
section are made to be directly implementable on a computer, we shall
try to write all estimates as explicitly as possible. Let
$\tilde U\Subset U$ and $\h_0\in\interval[open left]0{\tilde\h_0}$ be
such that $(G_\h)\restr_{\tilde U}$ is invertible onto a neighborhood
of $\overline{B}$ for all $\h\in\cI\cap\interval[open left]{0}{\h_0}$,
see Item~\ref{item:6} of Lemma~\ref{lemm:chart-pot-pourri}. From now
on in this section, every $\h$ is tacitly assumed to belong to
$\cI\cap\interval[open left]{0}{\h_0}$.

Let $U_0=G_0^{-1}(B_0)$ and $U_\h=G_\h^{-1}(B_0)$. Then
$U_0\Subset G_0^{-1}(B)$ and for all $\h\leq\h_0$,
$U_\h\Subset \tilde U$. It follows from the asymptotic expansion of
$G_\h^{-1}$ that there exists $C>0$ such that
\[
  \forall \h\leq \h_0, \quad U_\h \subset U_0 + B(0,C \h)\,,
\]
and hence for any $R>0$, there exists
$\h_0^{[0]}\in\interval[open left]{0}{\h_0}$ such that
\hequa
\begin{equation}
  \label{equ:domains}
  \forall \h \leq \h_0^{[0]},
  \quad U_\h + B(0,R \h) \subset \tilde U_0 \quad
  \text{ and } \quad
  U_\h + B(0,R \h) \subset \tilde U.
\end{equation}
\nequa
As in Lemma~\ref{lemm:labelling}, let $L_F$ be an upper bound on the
Lipschitz constant of $G_\h^{-1}$ on a neighborhood of
$\overline{B}$. Choose $\epsilon\in(0,\frac{1}{L_F})$, and let
$\h_0^{[1]}\in\interval[open left]{0}{\h_0}$ be small enough to
verify~\eqref{equ:h0-lambda-unique}, so that~\eqref{equ:unique-ball}
holds.

Note that, with the exception of the three points $\lambda_{0,0}$,
$\lambda_{1,0}$, and $\lambda_{0,1}$, all points are constructed by
the following process:
\begin{enumerate}[label=(\roman*)]
\item \label{item:unique1} Choose a point $\mu_1\in B_0\cap\Lh$.
\item Choose a ``vector'' $\vec v=\mu_0-\mu_{-1}$ (the difference
  between two previously constructed points $\mu_0$ and
  $\mu_{-1}\in B_0\cap\Lh$).
\item \label{item:unique3} Identify a closest point $\mu_2\in\Lh$ to
  $\mu_1+\vec v$.
\end{enumerate}
The following lemma shows that this process is uniformly well-defined
if these points lie in ball of size $\O(\h)$ and $\h_0$ is small
enough.  Then, Step~\ref{item:unique3} amounts to picking up the
natural parallel transport defined in~\eqref{equ:translation}:
\[
  \mu_2 = \mu_1 \plus{k_\h} (k_\h(\mu_0) - k_\h(\mu_{-1})).
\]
\begin{lemm}\label{lem:unique}
  Given $\epsilon\in(0,\frac{1}{L_F})$, there exists $L>0$ such that
  the following holds.  Choose $\h\leq\h^{[1]}_0$ and four points
  $\mu_i\in\Lh$ satisfying \ref{item:unique1}-\ref{item:unique3}
  above, and let $k_i := k_\h(\mu_i)$ be the corresponding
  multi-integers in $\ZM^2$, for $i\in\{-1,0,1,2\}$. Let $\rho>0$ be
  such that $\mu_1$ and $\mu_{-1}$ belong to the ball
  $B(\mu_0,\rho\h)$. If $\h\leq \h_0^{[0]}$ defined
  by~\eqref{equ:domains} with some $R > 2L_F\rho$, and if for some
  $N> 2$, $\h$ satisfies the inequalities
  \hequa
  \begin{equation}
    \label{eq:hbound}
    3C_N\h^{N-1} < \frac{R}{L_F} - 2\rho \quad \text{ and }
    \quad   \h \big( L(\rho + C_N \h^{N-1})^2 + 4C_N
    \h^{N-2}\big) \leq \epsilon\,,
  \end{equation}
  \nequa
  where $C_N$ is defined in~\eqref{equ:chart-labelling}, then $\mu_2$
  is unique and $k_2 = k_1 + k_0 -k_{-1}$.
\end{lemm}
\begin{demo} We wish to consider the point $\lambda_2\in\Lh$ whose
  label is $k_1 + k_0-k_{-1}$. First of all we show that
  $\h(k_1 + k_0-k_{-1})\in\tilde U_0$.  For all $\lambda\in\Lh$, we
  get from~\eqref{equ:chart-labelling}
  \[
    \norm{\h k-G_\h^{-1}(\lambda)} \leq L_F C_N\h^N\,.
  \]
  Therefore, $\h k_j$, for $j\in\{-1,0,1\}$, belongs to
  $U_\h + B(0,L_F C_N\h^N) \subset\tilde U$ provided
  $L_F C_N\h^{N-1} < R$, which is a consequence of the first
  inequality in~\eqref{eq:hbound}.  Thus, we may apply $G_\h^{-1}$ and
  use~\eqref{equ:chart-labelling} to obtain, for $j\in\{-1,1\}$,
  \begin{equation}\label{equ:k0k-1}
    \norm{\h k_0 - \h k_j} \leq L_F\norm{\mu_0 - \mu_j} + 2 L_F
    C_N\h^N < 2L_F(\h\rho + C_N\h^N)\,.
  \end{equation}
  Using the inequality for $j=-1$ we get
  $\h(k_1 + k_0-k_{-1})\in U_\h + B(0,L_F(3C_N\h^N + 2\h\rho))$.
  From~\eqref{equ:domains} and the first inequality
  in~\eqref{eq:hbound}, we get $\h(k_1 + k_0-k_{-1})\in\tilde U_0$.
  Since $U_0\Subset G_0^{-1}(B)$, by Item~\ref{item:def-asym-latt2b}
  of Definition~\ref{defi:lattice} we may define
  $\lambda_2:= k_\h^{(-1)}(k_1 + k_0-k_{-1}) \in \Lh$, \emph{i.e.},
  for all $N\geq 0$,
  \[
    \norm{\lambda_2 - G_\h(\h(k_1 + k_0-k_{-1}))} \leq C_N \h^N\,.
  \]
  We have, using a Taylor expansion at $\h k_1$,
  see~\eqref{equ:G-taylor2},
  \begin{align*}
    \norm{\lambda_2 - (\mu_1 + \vec v)}
    & \leq \norm{G_\h(\h(k_1 + k_0-k_{-1})) -  (\mu_1 + \vec v)} +  C_N \h^N\\
    & \leq \norm{G_\h(\h k_1) - \mu_1 + G'_\h(\h k_1)\cdot(\h k_0-\h k_{-1}) - \vec v}
      + L_1 \norm{\h k_0-\h k_{-1}}^2 + C_N \h^N \\
    & \leq \norm{G'_\h(\h k_1)\cdot(\h k_0-\h k_{-1}) - \vec v}
      + L_1 \norm{\h k_0-\h k_{-1}}^2 + 2C_N \h^N \,.
  \end{align*}
  We can also write
  $\norm{\vec v - G_\h(\h k_0) + G_\h(\h k_{-1})} \leq 2 C_N \h^N$.
  By Taylor expanding at $\h k_0$,
  \begin{equation*}
    \norm{ G_\h(\h k_0) - G_\h(\h k_{-1}) - G'_\h(\h k_0)\cdot(\h
      k_0-\h k_{-1})} \leq L_1\|\h k_0-\h k_{-1}\|^2\,.
  \end{equation*}
  Hence
  \begin{align*}
    \norm{\lambda_2 - (\mu_1 + \vec v)}
    & \leq  \norm{(G'_\h(\h k_0) - G'_\h(\h k_1))\cdot
      (\h k_0-\h k_{-1})} + 2L_1 \norm{\h k_0-\h k_{-1}}^2 + 4C_N \h^N\,.
  \end{align*}
  If follows from Item~\ref{item:1} of
  Lemma~\ref{lemm:chart-pot-pourri} that there exists a constant
  $L_2>0$ such that
  \[
    \forall [\xi_1,\xi_2]\subset \tilde U, \forall v\in\RM^n, \quad
    \norm{G_\h'(\xi_2)\cdot v - G_\h'(\xi_1) \cdot v} \leq L_2
    \norm{\xi_2 - \xi_1} \norm{v}.
  \]
  This gives
  \[
    \norm{(G'_\h(\h k_0) - G'_\h(\h k_1))\cdot (\h k_0-\h k_{-1})}
    \leq L_2 \norm{\h k_0-\h k_1}\norm{\h k_0-\h k_{-1}}\,,
  \]
  and hence
  \[
    \norm{\lambda_2 - (\mu_1 + \vec v)} \leq \norm{\h k_0-\h
      k_{-1}}\left(2L_1\norm{\h k_0-\h k_{-1}} + L_2\norm{\h k_0-\h
        k_{1}}\right) + 4C_N \h^N\,.
  \]
  Using~\eqref{equ:k0k-1} again,
  \[
    \norm{\lambda_2 - (\mu_1 + \vec v)} < \h^2\big(4(2L_1+L_2)
    L_F^2(\rho + C_N \h^{N-1})^2 + 4C_N \h^{N-2}\big)\,.
  \]
  Hence, for $\h$ small enough: precisely, as soon as
  \begin{equation*}
    \h\big(4(2L_1+L_2) L_F^2(\rho + C_N \h^{N-1})^2 + 4C_N
    \h^{N-2}\big) \leq \epsilon \,,
  \end{equation*}
  we may apply~\eqref{equ:unique-ball} to see that the closest point
  to $(\mu_1 + \vec v)$ in $\Lh$ must be $\lambda_2$, which proves the
  lemma with $L = 4(2L_1+L_2) L_F^2$.
\end{demo}

Of course, as in the end of Remark~\ref{rema:N=1}, one can simplify
estimates by choosing specific values, for instance
$\epsilon=\frac{1}{3L_F}$, $R=3L_F\rho$, and then we may
replace~(\ref{eq:hbound}) by the stronger assumption that there exist
$N>1$ such that
\[
  C_N \h^{N-1} \leq \min (\frac{\rho}{3}\;;\; \frac{\epsilon}{8})
  \quad \text{ and } \h \leq \frac{\epsilon}{4 L \rho^2}\,.
\]

\begin{coro}\label{cor:unique}
  With $\lambda_{m,n}$ the collection of points constructed in the
  algorithm without Step~\ref{item:algo-determinant}, let
  $k_{n,m} = k_\h(\lambda_{m,n})$. Set
  \begin{equation}
    \label{eq:z1z2}
    \vec z_1 = k_{1,0}-k_{0,0} \quad\text{ and }
    \quad \vec z_2 = k_{0,1}-k_{0,0}\,.
  \end{equation}
  There exists $\h_0^{[2]}>0$ (see Remark~\ref{rema:h2} below) such
  that for all $\h\leq\h_0^{[2]}$, we have
  \begin{equation}
    \label{eq:knm}
    k_{n,m}-k_{0,0} = n\,\vec z_1 + m\,\vec z_2\,.
  \end{equation}
  Moreover, for all $N\geq 1$,
  \begin{equation}
    \label{equ:lambda01-lambda00}
    \norm{\lambda_{1,0}-\lambda_{0,0}} \leq \h\rho_0\,,\quad
    \norm{\lambda_{0,1}-\lambda_{0,0}} \leq \h\rho_0\,,\quad
    \text{ with }  \quad  \rho_0 := L_0 + 2C_N\h^{N-1}\,,
  \end{equation}
  where $L_0$ is an upper bound on the Lipschitz constant of $G_\h$ on
  $\tilde U$, see~(\ref{equ:G-taylor1}), and
  \begin{equation}
    \label{equ:zj-bound}
    \norm{\vec z_1} \leq \tilde\rho\,, \quad \norm{\vec z_2} \leq
    \tilde\rho\,, \quad \text{ with } \quad \tilde\rho :=
    L_F(L_0+ 4C_N\h^{N-1})\,.
  \end{equation}
  After Step \ref{item:algo-determinant}, $\vec z_1$ and $\vec z_2$
  are possibly swapped.
\end{coro}
\begin{demo}
  Let $\lambda_1 = k_\h^{(-1)}(k_{0,0}+\vec e_1)$, where
  $\vec e_1=(1,0)$. Necessarily, $\lambda_1\neq \lambda_{0,0}$; hence,
  by Step \ref{item:algo-1-0},
  $\|\lambda_1-\lambda_{0,0}\|\geq\|\lambda_{1,0}-\lambda_{0,0}\|$.
  Since (using~\eqref{equ:chart}),
  \begin{equation}
    \label{equ:lambda1-lambda00}
    \norm{\lambda_1-\lambda_{0,0}} \leq \norm{G_\h(\h(k_{0,0}+\vec
      e_1))-G_\h(\h k_{0,0})} + 2C_N\h^N \leq L_0 \h + 2C_N\h^N\,;
  \end{equation}
  this shows the first inequality in~\eqref{equ:lambda01-lambda00}.
  We now proceed with Step~\ref{item:algo-n-0} of the algorithm by
  applying Lemma~\ref{lem:unique} to the triple
  $(\lambda_{0,0}, \lambda_{1,0}, \lambda_{1,0})$, with
  $\rho = \rho_0$, and thus obtain a unique
  $\lambda_{2,0} = \lambda_{0,0}\plus{k_\h} \vec z_1$.
  From~\eqref{equ:lambda1-lambda00} we immediately obtain
  \[
    \begin{aligned}
      \h\|\vec z_1\| & \leq L_F \|G_\h(\h k_{1,0})-G_\h(\h k_{0,0})\|
      \leq L_F \|\lambda_{1,0}-\lambda_{0,0}\| + 2L_FC_N\h^N \\
      &\leq L_FL_0 \h + 4L_FC_N\h^N = \h \tilde\rho\,.
    \end{aligned}
  \]
  This proves the first inequality of~\eqref{equ:zj-bound}. In order
  to repeat the application of Lemma~\ref{lem:unique} to the triple
  $(\lambda_{1,0}, \lambda_{2,0}, \lambda_{2,0})$ we need to estimate
  $\norm{\lambda_{2,0} - \lambda_{1,0}}$. Since
  $k_{2,0} = k_{1,0} + \vec z_1$, we get
  \begin{equation}
    \label{equ:dist-lambda02}
    \norm{\lambda_{2,0} - \lambda_{1,0}} \leq
    L_0 \h\norm{\vec z_1} + 2C_N\h^N \leq \h \rho_1
  \end{equation}
  with $\rho_1 := L_0\tilde\rho + 2C_N\h^{N-1}$. Thus we may apply
  Lemma~\ref{lem:unique} with $\rho = \rho_1$, and obtain that
  $\lambda_{3,0}$ is labelled by $k_{0,3} = k_{0,0} + 2 \vec z_1$.
  Therefore, we can estimate
  $\norm{\lambda_{3,0} - \lambda_{2,0}}\leq \h \rho_1$ exactly as
  in~\eqref{equ:dist-lambda02}. Repeating this process with the same
  $\rho=\rho_1$, we complete Step~\ref{item:algo-n-0} to obtain all
  $\lambda_{n,0}$ as long as they belong to $B_0$, and obtain, for all
  $n$, $\norm{\lambda_{n,0} - \lambda_{n-1,0}}\leq \h \rho_1$ and
  \begin{equation}
    \label{eq:kn0}
    k_{n,0}-k_{0,0} = n\,\vec z_1, \quad \text{ \emph{i.e.} }
    \quad \lambda_{n,0} = \lambda_{0,0} \plus{k_\h} n\vec z_1\,.
  \end{equation}
 
  Next, we consider $\lambda_{(0,1)}$ from Step~\ref{item:algo-0-1}.
  If $\lambda_1$ was not labelled in Step~\ref{item:algo-n-0}, we have
  $\|\lambda_1-\lambda_{0,0}\|\geq\|\lambda_{0,1}-\lambda_{0,0}\|$,
  which leads to the same estimates of as above,
  namely~\eqref{equ:lambda01-lambda00} and~\eqref{equ:zj-bound}
  hold. If, on the contrary, $\lambda_1$ was labelled in Step
  \ref{item:algo-n-0}, then by \eqref{eq:kn0} there exists a non-zero
  integer $n$ such that $\vec e_1 = n\vec z_1$.  Thus the point
  $\lambda_2 = k_\h(k_{0,0}+\vec e_2)$ was not labelled in Step
  \ref{item:algo-n-0}. Therefore
  $\|\lambda_2-\lambda_{0,0}\|\geq\|\lambda_{0,1}-\lambda_{0,0}\|$ and
  we the above estimates still hold, with the same $\rho$, which
  proves~\eqref{equ:lambda01-lambda00} and~\eqref{equ:zj-bound}.
  Therefore, we may continue to follow the algorithm and use Lemma
  \ref{lem:unique} at each step, proving the corollary.
\end{demo}
\begin{rema}
  \label{rema:h2} Let us investigate how small $\h_0^{[2]}$ should be
  for Corollary~\ref{cor:unique} to hold. We do not try to have an
  optimal bound, but rather to check what the geometric constraints
  are. First, we have the upper bound $\h_0^{[1]}$ defined
  by~\eqref{equ:h0-lambda-unique}, which ensures that $k_\h$ is well
  defined and one-to-one. This one is quite weak in principle, because
  we are free to increase the exponent of $\h$ to make it smaller.
  Then we apply Lemma~\ref{lem:unique} with $\rho = \rho_0$, which
  gives another upper bound given by~\eqref{eq:hbound}, where we are
  free to choose another $N>2$. Essentially, this means
  $\h\lesssim \frac{\epsilon}{L \rho_0} \sim
  \frac{1}{L_0(L_1+L_2)L_F^3}$. Another application of
  Lemma~\ref{lem:unique} with $\rho = \rho_1\sim L_F L_0^2$ gives a
  new bound~\eqref{eq:hbound}, which, roughly speaking, imposes
  $\h\lesssim \frac{1}{L_0^2(L_1+L_2)L_F^4}$, which is \emph{a priori}
  stronger than the previous one, at least if $L_F,L_0\geq 1$ and we
  neglect the term $C_N\h^N$. However, at each application of the
  lemma we are free to optimize by choosing a different $N$. In both
  cases one needs to select $R> 2L_F\rho$, which gives yet another
  bound by adjusting $\h_0^{[0]}\sim o(1/R)$,
  see~\eqref{equ:domains}. Given $R$, this last bound only depends on
  the size of the domain $B_0$ within $B$; we can improve it if
  necessary by choosing a smaller $B_0$.
\end{rema}

We now consider in more details the construction of the first three
points $(\lambda_{0,0}, \lambda_{1,0}, \lambda_{0,1})$ near
$c$. Naturally, we assume that $\h$ is small enough so that $\Lh$
contains at least three points.  Our aim is to prove that $\vec z_1$
and $\vec z_2$ form a $\ZM$-basis of $\ZM^2$.  After $\lambda_{0,0}$
is chosen, $\vec v_1= \lambda_{1,0}-\lambda_{0,0}$ and
$\vec v_2 = \lambda_{0,1}-\lambda_{0,0}$ are chosen by a minimization
process. We have to be careful with the fact that the vectors
$\vec z_i$ do not satisfy the same minimization properties because
$G_\h'(\h k_{0,0})$ is not orthogonal, in general.
  
\begin{lemm}\label{lemm:coprime}
  There exists $\h_0^{[3]}>0$ (given by~\eqref{equ:h0-r0}
  and~\eqref{equ:hbar3} below) such that if $\h\leq \h_0^{[3]}$, then
  the entries of $\vec z_1$ are co-prime integers.
\end{lemm}
\begin{demo}
  Assume that there exists $\vec z_0\in\ZM^2\setminus\{0\}$ and
  $n\in\NM^*$ such that $\vec z_1 = n \vec z_0$, where $n$ may depend
  on $\h$. Let us first show that
  $\h(k_{0,0} + \vec z_0)\in \tilde U_0$ (notice that $\tilde U_0$ is
  not assumed to be convex).  From Lemma~\ref{lemm:lattice-dense}, for
  any $N$ we have
  \[
    \norm{c - \lambda_{0,0}} \leq (M + L_0)\h + C_N\h^N.
  \]
  Likewise, if we let $k_{0,0}:= k_\h(\lambda_{0,0})$ and
  $\xi := G_0^{-1}(c)$, then $\xi\in\tilde B_0$ and we have
  \begin{align*}
    \norm{\xi - \h k_{0,0}}
    & \leq L_F \left(\norm{G_\h(\xi) - \lambda_{0,0}} +
      \norm{\lambda_{0,0} - G_\h( \h k_{0,0})}\right) \nonumber \\
    & \leq L_F\left(M \h + \norm{c - \lambda_{0,0}} + C_N\h^N\right)\nonumber \\
    & \leq L_F(2M + L_0 + 2C_N\h^{N-1})\h.
  \end{align*}
  Therefore, using~\eqref{equ:zj-bound},
  \[
    \norm{\h(k_{0,0} + \vec z_0) - \xi} \leq \norm{\h k_{0,0}- \xi} +
    \h\norm{\vec z_0} \leq L_F(2M + L_0 + 2C_N\h^{N-1})\h +
    \h\tilde\rho/{n}
  \]
  Since $\tilde B_0$ is open and independent of $\h$, there is $r_0>0$
  such that $B(\xi, r_0)\subset \tilde U_0$.  Assume
  $\h (L_F(2M + L_0 + 2C_N\h^{N-1}) + {\tilde \rho}) < r_0$,
  \emph{i.e.}
  \hequa
  \begin{equation}
    \label{equ:h0-r0}
    \h L_F(2M + 2L_0 + 6C_N\h^{N-1}) < r_0\,.
  \end{equation}
  \nequa
  Then $\h (k_{0,0} + \vec z_0)\in B(\xi, r_0)\subset \tilde U_0$. We
  may now let $\mu:=k_\h^{-1}(\h(k_{0,0} + \vec z_0))$ be the
  corresponding element in $\Lh$.  By~\eqref{equ:G-taylor2}, for any
  $\vec z\in\RM^2$ such that $\norm{\vec z}\leq \norm{\vec z_1}$,
  \begin{equation}
    \label{equ:norme-reste}  
    \norm{G_\h(\h k_{0,0} + \h \vec z) - G_\h(\h k_{0,0}) - G_\h'(\h
      k_{0,0})\cdot(\h \vec z)} \leq L_1\norm{\h \vec z}^2\,
  \end{equation}
  and hence
  \begin{align}
    \norm{\lambda_{1,0} - \lambda_{0,0}}
    & = \norm{G_\h'(\h k_{0,0})\cdot(\h \vec z_1)} + \sigma_1 \nonumber \\
    \norm{\mu - \lambda_{0,0}}
    & = \norm{G_\h'(\h k_{0,0})\cdot(\h \vec z_0)} + \sigma_0 \,,
  \end{align}
  with $\abs{\sigma_j}\leq L_1\norm{\h \vec z_1}^2 + 2C_N \h^N$, for
  $j=0,1$.  From the algorithm we know that
  $\norm{\lambda_{1,0} - \lambda_{0,0}} \leq \norm{\mu -
    \lambda_{0,0}}$ and therefore
  \[
    (n-1) \norm{G_\h'(\h k_{0,0})\cdot(\h \vec z_0)} \leq {\sigma_0 -
      \sigma_1}.
  \]
  From Lemma~\ref{lemm:chart-pot-pourri}, there exists $\Gamma>0$
  independent of $\h$ such that $\norm{G_\h'^{-1}}\geq \Gamma$ on
  $\tilde U$, and hence
  \[
    (n-1) \leq \frac{\sigma_0 - \sigma_1}{\Gamma\norm{\h \vec z_0}}
    \leq \frac{2L_1\h^2\tilde\rho^2+ 4C_N\h^N}{\Gamma\h}\,.
  \]
  Thus, if
  \hequa
  \begin{equation}
    \label{equ:hbar3}
    \h(2L_1\tilde\rho^2+ 4C_N\h^{N-2})<\Gamma,
  \end{equation}
  \nequa
  then we must have $n=1$.
\end{demo}

\begin{lemm}
  For all $\h \leq \h^{[3]}_0$, the vectors $\vec z_1$ and
  $\vec z_2$~\eqref{eq:z1z2} defined in Corollary~\ref{cor:unique} are
  linearly independent.
\end{lemm}
\begin{demo}
  If $\vec z_1$ and $\vec z_2$ are colinear, there exists
  $\sigma\in\QM$ such that $\vec z_2=\sigma\vec z_1$. Since the
  coefficients of $\vec z_1$ are co-prime, we must have
  $\sigma=n\in\ZM$.  Writing,
  $k_{0,1} = k_{0,0} + \vec z_2 = k_{0,0} + n\vec z_1$, we obtain from
  Corollary~\ref{cor:unique} that $\lambda_{0,1} = \lambda_{n,0}$,
  which contradicts Step~\ref{item:algo-0-1} of the algorithm.
\end{demo}

We arrive at the main result of this section.
\begin{theo}\label{theo:basis}
  There exists $\h_0^{[4]}>0$ (defined in~\eqref{equ:h04}) such that
  $\vec z_1$ and $\vec z_2$ defined in Corollary~\ref{cor:unique} form
  a $\ZM$-basis of $\ZM^2$ for any $\h\leq \h^{[4]}_0$.
\end{theo}
\begin{demo}
  Let $D=D(\h):=\det(\vec z_1,\vec z_2)$. We know that $D\neq0$ by the
  previous lemma.  We want now to prove that $D=\pm1$. By Pick's
  formula, the number of integral points in the closed convex hull of
  the three points $k_{0,0}$, $k_{1,0}$ and $k_{0,1}$ (other than
  these vertices), where points in the boundary count half, is
  $\frac{|D|-1}{2}$. Hence if $|D|\geq 2$, there is at least one such
  integral point, call it $\ell\in\ZM^2$, set
  $\vec z_3 = \ell - k_{0,0}$, and let $\mu= k_\h^{-1}(\h\ell)$ be the
  corresponding element of $\Lh$.
  \\[1ex]
  (i) Since the components of $\vec z_1$ are co-prime by
  Lemma~\ref{lemm:coprime}, $\vec z_3$ cannot be colinear to
  $\vec z_1$, because it then would be equal to it, and $\ell$ would
  be a vertex of the triangle. Therefore, by
  Corollary~\ref{cor:unique}, $\mu$ cannot be one of the
  $\lambda_{n,0}$. This implies
  \begin{equation}
    \label{equ:lambda-mu}
    \norm{\lambda_{0,1}-\lambda_{0,0}} \leq \norm{\mu -\lambda_{0,0}}.
  \end{equation}
  \\
  (ii) Set $\vec u_i = G_h'(\h k_{0,0})\cdot\vec z_i$, $i=1,2,3$.  As
  $\vec z_3$ is a convex linear combination of $\vec z_1$ and
  $\vec z_2$, the same holds between $\vec u_3$, and $\vec u_1$,
  $\vec u_2$. Using~\eqref{equ:norme-reste}, we see that there exists
  a constant $C>0$ that can be made explicit (namely
  $C\geq 2L_1\tilde \rho + 2C_N\h^{N-2}$) such that the ball
  $B(\mu,\h^2 C)$ intersects the interior of the triangle
  $(\lambda_{0,0}, \lambda_{1,0}, \lambda_{0,1})$. Let $\tilde\mu$ be
  a point of this intersection. Since
  $\norm{\vec z_j - \vec z_3}\geq 1$ for $j=1,2$, because these are
  non zero integer vectors, the point $\tilde\mu$ must stay away from
  a ball of size $\delta \h$ of the vertices, for some
  $\delta\in\interval[open]0{1/L_F}$ independent of $\h$ and which can
  also be made explicit.
  \\[1ex]
  (iii) We now need some elementary triangle estimates in order to
  bound from below the distance
  $\norm{\lambda_{1,0}-\lambda_{0,0}} - \norm{\tilde\mu -
    \lambda_{0,0}}$.  By construction of the algorithm,
  $\|\lambda_{1,0}-\lambda_{0,0}\| \leq
  \|\lambda_{0,1}-\lambda_{0,0}\|$. Hence the angle at the vertex
  $\lambda_{0,1}$ is strictly less than $\frac{\pi}{2}$, which implies
  that the orthogonal projection $H$ of $\lambda_{0,0}$ onto the line
  $(\lambda_{1,0}\lambda_{0,1})$ is located on the strict half line
  starting at $\lambda_{0,1}$ and containing $\lambda_{1,0}$. Thus
  there exists $r>0$ be such that the ball $B(\lambda_{0,1},r)$ does
  not contain $H$ (see Figure~\ref{fig:triangle}).
  \begin{figure}[ht]\label{fig:triangle}
    \centering \includegraphics[width=0.2\linewidth]{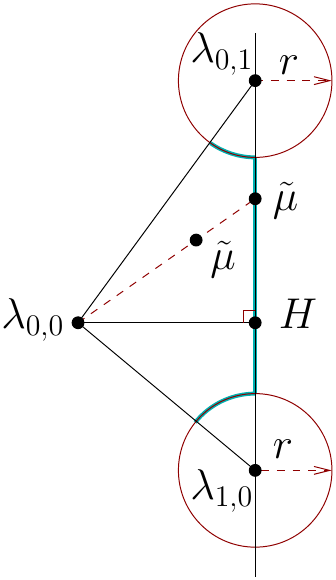}
    \caption{The point $\tilde\mu$ is pushed onto the thickened path.}
  \end{figure}
  In fact any $r\leq \norm{\lambda_{1,0} - \lambda_{0,1}}/2$ fulfills
  this requirement. We choose
  $r=\min(\norm{\lambda_{1,0} - \lambda_{0,1}}/2, \delta\h)$ so that,
  in addition,
  $\tilde \mu\not\in B(\lambda_{0,1},r)\cup B(\lambda_{1,0},r)$. Let
  us push $\tilde\mu$ along the ray $(\lambda_{0,0}\mu)$ until it
  reaches either the boundary of one of the balls $B(\lambda_{0,1},r)$
  or $B(\lambda_{1,0},r)$, or the opposite edge of the triangle,
  $E:=(\lambda_{1,0}\lambda_{0,1})$. We may call $\tilde\mu$ again
  this new point. Doing this, the distance
  $d:=\norm{\tilde\mu - \lambda_{0,0}}$ can not decrease. When
  $\tilde{\mu}$ is on the boundary of a ball, it is clear that $d$
  increases as we move $\tilde\mu$ on the circle towards the edge
  $E$. So, for any configuration of $\tilde\mu$, there is a position
  of $\tilde\mu$ on $E$ that produces a larger or equal distance
  $d$. Now it is easy to study the variation of $d = d(x)$ as a
  function of the abscissa $x$ of $\tilde{\mu}$ along $E$. Taking
  $x=0$ for the point $H$ and $x=:b>0$ at $\lambda_{0,1}$, and using
  $d(x) = \sqrt{x^2 + h^2}$, with $h = \norm{H-\lambda_{0,0}}$, we
  obtain for $x\in(0,b-r)$ that $d$ is increasing and
  \[
    d(b) - d(x) \geq \frac{r b }{2d(b)} \geq \frac{r^2}{2d(b)}.
  \]
  If the point $\lambda_{1,0}$ has a negative abscissa $c<0$, the
  distance $d$ will have a local maximum also at $c$, and we can
  repeat the argument above to get, when $\abs{c}>r$,
  \[
    \forall x\in (c+r,0), \qquad d(c) - d(x) \geq \frac{r \abs{c}
    }{2d(c)} \geq\frac{r^2}{2d(c)}.
  \]
  Recall that
  $d(c)= \|\lambda_{1,0}-\lambda_{0,0}\| \leq
  \|\lambda_{0,1}-\lambda_{0,0}\|= d(b)$.  Thus,
  \[
    \forall x\in (c+r,b-r), \qquad d(b) - d(x) \geq \frac{r^2
    }{2d(b)}\,.
  \]
  This finally gives
  \begin{equation}\label{equ:triangle}
    \norm{\lambda_{0,1}-\lambda_{0,0}} - \norm{\tilde\mu -
      \lambda_{0,0}} \geq \frac{r^2}{2\norm{\lambda_{0,1} - \lambda_{0,0}}}\,.
  \end{equation}
  \\[1ex]
  (iv) From~\eqref{equ:lambda-mu} we can write
  \begin{align*}
    \norm{\lambda_{0,1} -\lambda_{0,0}}
    & \leq \norm{\mu -\lambda_{0,0}}
      \leq \norm{\tilde\mu -\lambda_{0,0}} + C\h^2 \\
    & \leq \norm{\lambda_{0,1}-\lambda_{0,0}} -
      \frac{r^2}%
      {2\norm{\lambda_{0,1} - \lambda_{0,0}}} + C\h^2
      \quad \text{ by~\eqref{equ:triangle}}.
  \end{align*}
  Therefore
  \[
    \frac{r^2}{2\norm{\lambda_{0,1} - \lambda_{0,0}}} \leq C\h^2.
  \]
  Since $\norm{k_{0,1} - k_{1,0}}\geq 1$, we have
  $\norm{\lambda_{0,1} - \lambda_{1,0}} \geq \h /L_F - 2 C_N\h^N$,
  while~\eqref{equ:lambda01-lambda00} says
  $\norm{\lambda_{0,1} - \lambda_{0,0}} \leq L_0\h + 2C_N\h^N$. We
  obtain
  \[
    \frac{\min({(\frac{1}{L_F} - 2C_N\h^{N-1})}^2, \delta^2)}{2(L_0 +
      2C_N\h^{N-1})} \leq C\h\,
  \]
  which is of course impossible if $\h$ is small enough, namely if one
  takes $\h\leq \h_0^{[4]}$ with
  \hequa
  \begin{gather}
    \label{equ:h04}
    2C_N(\h_0^{[4]})^{N-1} < \frac{1}{L_F} \quad \text{ and }
    \quad \delta < \frac{1}{L_F} - 2C_N(\h_0^{[4]})^{N-1}\,\nonumber\\
    \text{ and } \quad \h_0^{[4]} < \frac{\delta^2}{2C(L_0 +
      2C_N(\h_0^{[4]})^{N-1})}\,.  
  \end{gather}
  \nequa

  For all $\h\leq\h_0^{[4]}$, we conclude that $|D|=1$. This concludes
  the proof of the theorem.
\end{demo}
\begin{rema}
  In order to have a rough idea of the size of the various bounds on
  $\h$, following up Remark~\ref{rema:h2}, we may neglect the terms
  $C_N\h^N$ and get from~\eqref{equ:h0-r0} and~\eqref{equ:hbar3} and
  the approximation $\Gamma\sim 1/L_0$ that
  $\h_0^{[3]}\lesssim \min\left(\h_0^{[2]}\,;\,
    \frac{r_0}{2L_F(M+L_0)}\,;\, \frac{1}{2L_0^3L_1L_F^2}\right)$.
  Then, from~\eqref{equ:h04} and $C\sim 2L_0^2L_1L_F^2$, we get
  $ \h_0^{[4]}\lesssim \min\left(\h_0^{[3]}\,;\,
    \frac{1}{4L_0^3L_1L_F^4} \right)$.
\end{rema}

\begin{defi}\label{defi:affine-basis}
  When the vectors $(\vec z_1,\vec z_2)$ form a $\ZM$-basis of
  $\ZM^2$, the triple $(\lambda_{0,0}, \lambda_{1,0}, \lambda_{0,1})$
  will be called an \emph{affine basis} of $\Lh$ at $\lambda_{0,0}$.
\end{defi}

The algorithm does not necessarily label \emph{all} points of
$\Lh\cap B_0$. Indeed, by construction, the set of produced labels
$(n,m)$ is of the form
\[ {\cal E} := \{(n,m)\in\ZM^2; \quad m_\textup{min} \leq m \leq
  m_\textup{max}; \quad n_\textup{min}(m) \leq n \leq
  n_\textup{max}(m)\}\,,
\]
where $m_\textup{min}$, $m_\textup{max}$, and $n_\textup{min}(m)$ and
$n_\textup{max}(m)$, for
$m\in\{m_\textup{min}, \dots, m_\textup{max}\}$, may depend on
$\h$. For a given $m$, $n_\textup{max}(m)$ is the smallest positive
integer produced by Step~\ref{item:algo-n-0} such that
$\lambda_{n,m}\in B_0$ and $\lambda_{n+1,m}\not\in B_0$. The integer
$n_\textup{min}(m)$ is defined in a similar way, $m_\textup{max}$ is
the smallest positive integer produced by Step~\ref{item:algo-0-2}
such that $\lambda_{0,m}\in B_0$ but $\lambda_{0,m+1}\not \in B_0$,
and $m_\textup{min}$ is defined in a similar way. Thus, there is no
reason why a set of the form ${\cal E}$ would fill up $\Lh\cap B_0$
entirely. However, since ${\cal E}$ will fill up all integral points
of any convex set $V\subset\tilde U_0$, the algorithm is guaranteed to
label all points of $\Lh$ in \emph{some} $\h$-independent open subset
containing $c$. It would be interesting to improve the algorithm in
order to make sure that it explores the whole connected component of
$B_0$.

Note that, after the last step (orientation test) of the algorithm,
the basis $\vec v_1$, $\vec v_2$ is made direct, so is the case for
the basis $\vec z_1$, $\vec z_2$ because by convention
$\det G_0'(\xi)>0$.  It follows from Theorem~\ref{theo:basis} and
Corollary~\ref{cor:unique} that for each $\h$ small enough, we can
define a matrix $Z_\h\in \textup{SL}(2,\ZM)$ such that
$Z_\h(\vec e_1,\vec e_2) = (\vec z_1,\vec z_2)$, and the ``labelling''
$ \lambda_{n,m}\mapsto (n,m)$ of the algorithm is such that
\[
  \lambda_{n,m} = G_\h(\h(k_{0,0} + n\vec z_1 + m\vec z_2)) +
  \O(\h^\infty) = G_\h\circ Z_\h (\h(\varkappa+n\vec e_1 + m\vec e_2)) +
  \O(\h^\infty)
\]
where $\varkappa=Z_\h^{-1}k_{0,0}\in\ZM^2$. However, this does not
produce a linear labelling for $\Lh$ because in general the map
$\h \to Z_\h$ will not be constant (and hence, not continuous), even
for arbitrary small values of $\h$. In order to produce a linear
labelling, we should find a way to ``detect'' the matrix $Z_\h$, which
will allow to correct the initial algorithm and make it smooth in
$\h$. This is the aim of the following sections.

\subsubsection{A semitoric algorithm}
\label{sec:semitoric-algorithm}

If a quantum system is known to be semitoric, the joint spectrum is an
asymptotic lattice with a special property given by
Proposition~\ref{prop:J-spectrum}.  In this case, we expect a well
designed algorithm not to provide \emph{any} labelling, but rather a
\emph{semitoric labelling}, \emph{i.e.} associated with a semitoric
chart, see Lemma~\ref{lemm:semitoric-chart}.

Let $(j,k)$ be the semitoric labelling of
Theorem~\ref{theo:q-semitoric}.  From~\eqref{equ:J-spectrum}, we know
that for any $\epsilon\in(0,\tfrac{1}{2})$, there exists a ball $B'$
around $c$ such that, if $\h_0$ is small enough, $\Sigma_\h\cap B'$ is
contained in a union of disjoint vertical strips of width
$2\epsilon\h$:
\begin{equation}\label{equ:union-strips}
  \Sigma_\h \cap B' \subset \bigcup_{j\in\ZM} V_j(\h),
\end{equation}
where
\begin{equation}\label{equ:strip-j}
  V_j(\h) = [\alpha+\h(j+\mu)-\h\epsilon, \alpha+\h(j+\mu)+\h\epsilon]
  \times \RM \subset \RM^2.
\end{equation}
The precise size of $B'$ depends on the variations of the subprincipal
term $G_1$ of a semitoric asymptotic chart
$G_\h\sim G_0 + \h G_1 + \cdots$.

Now, from the data of $\Lh$, perform the generic algorithm of the
previous section, in order to obtain
an oriented affine basis
$(\lambda_{0,0}, \lambda_{1,0}, \lambda_{0,1})$ of $\Lh$
(Theorem~\ref{theo:basis}, Definition~\ref{defi:affine-basis}, and
Step~\ref{item:algo-determinant}). Let $V$ be the vertical strip of
width $\h^{3/2}$, vertically centered at $\lambda_{0,0}$, and define
$\mu\in\Lh$ as the nearest point to $\lambda_{0,0}$ located in $V$ and
\emph{above} $\lambda_{0,0}$. The existence of a semitoric chart
(Lemma~\ref{lemm:semitoric-chart}) and the
decomposition~\eqref{equ:union-strips}-\eqref{equ:strip-j} imply that
$\mu$ exists in an $\O(\h)$ neighborhood of $\lambda_{0,0}$, belongs
to the strip $V_{j_0}$ that contains $\lambda_{0,0}$, and is unique if
$\h$ is small enough. Therefore, we know from the analysis of the
general algorithm (Corollary~\ref{cor:unique}) that there exists
bounded co-prime integers $(n,m)\in\ZM^2$ such that
\[
k_\h(\mu) - k_{0,0} =  n\,\vec z_1 + m\,\vec z_2\,.
\]
In practice, the integers $(n,m)$ can be found by expressing
$(\mu - \lambda_{0,0})/\h$ on the basis $(v_1,v_2)$ and then rounding
its coefficients to their nearest integers. Finally, we choose
$(n',m')\in\ZM^2$ such that $n'm - m'n = 1$, and define $\mu'\in\Lh$
by
\[
k_\h(\mu') =  n'\,\vec z_1 + m'\,\vec z_2\,.
\]
We obtain in this way a new oriented affine basis of $\Lh$ given by
\begin{equation}\label{equ:semitoric-affine-basis}
(\mu' - \lambda_{0,0}, \mu - \lambda_{0,0})\,.
\end{equation}
\begin{prop}\label{prop:semitoric-algorithm}
  The affine basis~\eqref{equ:semitoric-affine-basis} is associated
  with a semitoric asymptotic chart.
\end{prop}
\begin{demo}
  The uniqueness of $\mu$ ensures that its label, in a semitoric chart
  $\hat G_\h$ associated with a labelling $\hat k_\h$, is equal to
  $\hat k_\h(\lambda_{0,0}) + (0,1)$. Hence the label of $\mu'$ must
  be of the form $\hat k_\h(\lambda_{0,0}) + (1,\ell)$, for some
  $\ell\in\ZM$. Hence the affine
  basis~\eqref{equ:semitoric-affine-basis} is obtained from the one
  associated with the chart $\hat G_\h$ by composition with the matrix
  $A =
  \begin{pmatrix}
    1 & 0 \\ \ell & 1
  \end{pmatrix}\,.  $ This matrix $A$ preserves the semitoric property
  of the chart.
\end{demo}

\subsubsection{An algorithm for a sequence of values of $\h$}
In our way to reconstructing a linear labelling from the data of the
sets $\Lh$ and of a window $B$ with a distinguished point $c\in B$, we
introduce a ``algorithm with uniform labelling'' on a decreasing
sequence $\h_j\in\cI$, $j\geq1$, tending to $0$. The algorithm is
inductive: once the labelling
$k_{\h_j}: \Lambda_{\h_j} \ni \lambda^{(j)} \mapsto (n^{(j)}, m^{(j)})
\in\ZM^2$ is known, together with a matrix $S_j\in\textup{M}_2(\ZM)$,
it will produce the labelling
$k_{\h_{j+1}} : \Lambda_{\h_{j+1}}\to \ZM^2$ (and the matrix
$S_{j+1}$).  Thus, it theoretically defines $k_{\h_j}$ for all $j$. Of
course, in practice, if one wants to obtain the labelling $k_{\h_j}$
for a specific $j$, it is enough to stop at the step $j$.

The algorithm works by running the previous algorithm of
Section~\ref{sec:an-algorithm-fixed} with all values $\h_j$, for all
$j=1,\dots$, and self-adjusting the resulting labelling for each
$j$. In order to have a more efficient implementation, if we know in
advance that we want to stop at a specific step $j=j_\textup{stop}$,
it is in fact not necessary to compute the full labellings for the
values of $j$ less than $j_\textup{stop}$; for these values, it is
enough to find the correct ``affine basis'', which corresponds to
Steps~\ref{item:algo-c} to~\ref{item:algo-0-1} of the algorithm of
Section~\ref{sec:an-algorithm-fixed}. Thus, the new algorithm, with
exit test at $j=j_\textup{stop}$, works as follows.
\begin{enumerate}[label=\emph{\alph*})]
\item Choose an open subset $B_0\Subset B$, and fix $c\in B_0$. Let
  $S_0 = \textup{Id}\in \textup{M}_2(\RM)$. Let $j=1$.
\item\label{item:algo2-base} Apply steps~\ref{item:algo-c}
  to~\ref{item:algo-0-1} of the algorithm of
  Section~\ref{sec:an-algorithm-fixed} with $\h = \h_j$. This defines
  points $\lambda^{(j)}_{n,m}$. In particular, we have an origin
  $\lambda^{(j)}_{0,0}$, and the first generating vectors
  $v^{(j)}_1 = \lambda^{(j)}_{1,0}-\lambda^{(j)}_{0,0}$ and
  $v^{(j)}_2 = \lambda^{(j)}_{0,1}-\lambda^{(j)}_{0,0}$.
\item Define $T_j\in \textup{M}_2(\RM)$ to be the matrix formed by the
  column vectors:
  \[
    T_j := (\h_j^{-1} \,v^{(j)}_1 , \; \h_j^{-1}\,v^{(j)}_2 ).
  \]
  If $T_j$ is not invertible or $j=1$, increase $j$ by one, let
  $S_j := S_{j-1}$, and go back to Step~\ref{item:algo2-base}. If $T_j$
  is invertible, make it oriented as in
  Step~\ref{item:algo-determinant} of the previous algorithm.
\item Define $A_j:=T_{j}^{-1} T_{j-1}$, and let
  $A^\sharp_j\in \textup{M}_2(\ZM)$ be the matrix obtained by rounding
  the entries of $A_j$ to their ``nearest integer'' (in the usual,
  unique way). If $\det A^\sharp_j \neq 1$, define $S_j :=
  S_{j-1}$. If $\det A^\sharp_j = 1$, define
  \[
    S_j := S_{j-1}\,(A^\sharp_j)^{-1}\,.
  \]
\item In case of a concrete implementation, if $j < j_\textup{stop}$,
  increase $j$ by one and go back to
  Step~\ref{item:algo2-base}. Otherwise, finish the previous
  algorithm, \emph{i.e.} perform Steps~\ref{item:algo-1-1}
  to~\ref{item:algo-determinant}. Let $\lambda\mapsto (n,m)$ be the
  resulting labelling.
\item The new labelling of $\mathcal{L}_{\h_j}$ is the map
  $\lambda^{(j)}\mapsto (\tilde n,\tilde m)$ given by a linear
  transformation acting on the labelling $\lambda$ of the previous
  step according to
  \[
    \begin{pmatrix}
      \tilde n \\ \tilde m
   \end{pmatrix} =  S_j \begin{pmatrix} n \\ m
    \end{pmatrix}\,.
  \]
  In other words,
  $\lambda_{n,m}=\lambda^{(j)}_{S_j(m,n)}=\lambda^{(j)}_{\tilde
    m,\tilde n}$.
\end{enumerate}

\begin{theo}
  \label{theo:algo} Let $(\Lh, \cI, B)$ be an asymptotic lattice,
  where $B\subset\RM^2$. Let $\h_j\in\cI$, $j\geq1$, be a decreasing
  sequence tending to $0$.  Then the previously described algorithm
  produces a linear labelling of the asymptotic lattice
  $(\Lh, \cI', B)$, where $\cI'=\{\h_j, j\in\NM^*\}$.
\end{theo}
\begin{demo}
  We interpret the newly introduced objects with respect to an
  asymptotic chart $G_\h$ (which is known to exist, but is
  unknown). Let us denote by $Z_j$ the matrix formed by the column
  vectors $z^{(j)}_1=\vec z_1$ and $z^{(j)}_2=\vec z_2$ defined in
  Corollary \ref{cor:unique}\eqref{eq:z1z2} for $\h=\h_j$ (initial
  algorithm):
  \[
    Z_j := (z^{(j)}_1 , \; z^{(j)}_2 ).
  \]
  Formula~(\ref{equ:G-taylor2}), in view of~\eqref{equ:zj-bound},
  gives
  \[
    T_j = G'_0(\xi)\,Z_j + \O(\h_j)\,.
  \]
  By Theorem~\ref{theo:basis}, $Z_j$ is unimodular and hence $T_j$ is
  invertible (with positive determinant) if $\h_j$ is small enough,
  which happens for all $j\geq j_0$, for some $j_0$. Thus,
  \[
    A_j := T_{j}^{-1} T_{j-1} = Z_{j}^{-1} Z_{j-1} + \O(\h_{j-1})\,.
  \]
  Therefore, if $j_0$ is large enough, we obtain for all $j\geq j_0$,
  \[
    A^\sharp_j = Z_{j}^{-1} Z_{j-1} \quad \text{ and } \quad \det
    A^\sharp_j = 1\,.
  \]
  Set $\tilde Z_j = Z_j\, S_j^{-1}$. The matrix $\tilde Z_j$ has
  integer coefficients and satisfies by definition of the new
  labelling for $\h=\h_j$
  \[
    \tilde Z_j \begin{pmatrix} \tilde n \\ \tilde m
    \end{pmatrix} = Z_j \begin{pmatrix} n \\ m
    \end{pmatrix}.
  \]
  We check that the sequence $\tilde Z_j$ is stationary as
  $j\to\infty$:
  \[
    \tilde Z_j = Z_j\, S_j^{-1} = Z_j\,A^\sharp_j\,S^{-1}_{j-1}
    = Z_j\,Z_{j}^{-1} Z_{j-1}\,S^{-1}_{j-1} = Z_{j-1}\,S^{-1}_{j-1} =
    \cdots = \tilde Z_{j_0}\,.
  \]
  By \eqref{eq:knm}, in terms of the good labelling $k_\h$ associated
  with $G_\h$
  \[
    k_\h(\tilde\lambda^{(j)}_{\tilde n,\tilde m}) = \tilde n\,\vec z_1
    + \tilde m\,\vec z_2 + k^{(j)}_{0,0},
  \]
  \emph{i.e.}
  \[
    \tilde Z_{j_0}^{-1} \, k_\h(\tilde\lambda^{(j)}_{n,m}) = n\,\vec
    e_1 + m\,\vec e_2 + \tilde Z_{j_0}^{-1}\,k^{(j)}_{0,0}.
  \]
  Thus we have constructed a linear labelling $\lambda\to(n,m)$ for
  $(\Lh, \cI', B)$, where $(n,m)$ is such that
  $\lambda=\tilde\lambda^{(j)}_{n,m}$, associated with the asymptotic
  chart $\tilde G_\h := G_\h \circ \tilde Z_{j_0}$.
\end{demo}

\subsection{The inverse problem for the rotation number}
\label{sec:inverse}

We can finally apply our algorithms to the initial question, because
having a linear labelling is actually enough for recovering the
rotation number. In Theorem~\ref{theo:q-rotation-number} we have
obtained the classical rotation number as an $\O(h)$ limit of the
quantum rotation number; this suggests that one can actually expect
the recovery process to be robust with respect to smaller (namely,
$\O(\h^2)$) perturbations of the asymptotic lattice.  Previous results
for special cases of semiclassical operators show that this
expectation is very natural~\cite{san-inverse,san-lefloch-pelayo:jc}.

\begin{theo}
  \label{theo:A}
  Let $(\hat J, \hat H)$ be a quantum integrable system with principal
  symbols $F:=(J,H)$. Let $\Sigma_\h:=\Sigma_\h(\hat J, \hat H)$ be
  the joint spectrum of $(\hat J, \hat H)$. Let $c\in\RM^2$ be a
  regular value of $F$, such that the fiber $F^{-1}(c)$ is compact and
  connected. Then there exists a small ball $B$ around $c$ for which
  the following statements hold.

  \begin{itemize}
  \item From the knowledge of the family of discrete sets
    $B\cap \Sigma_{\h}\subset \RM^2$ where $\h$ varies in a sequence
    $\mathcal{I}' = (\h_j)_{j\in\NM^*}$ accumulating at zero, one can
    recover the classical rotation number $[w]$ in $B$ modulo the
    natural action of Möbius transformations (Lemma~\ref{lemm:ItoI'}).
  \item This recovery is stable under $\O(\h^2)$ perturbations of the
    joint spectrum, and is algorithmic; precisely, there exists a
    constructive algorithm with the following property: given any such
    data $(\Sigma_\h, \cI', B)$, there exists a choice of action
    coordinates $I$ in $B$ such that for any $b\in B$, the algorithm
    produces a sequence $[w_{\h_i}](b)\in\RP^{1}$, where
    $\h_i\in\cI'$, such that
    \begin{equation}
      [w_{\h_i}](b) = [w_I](\Lambda) + \O(\h_i)
      \label{equ:recover-w}
    \end{equation}
    where $\Lambda=F^{-1}(b)$.
  \end{itemize}
\end{theo}


\begin{demo}
  We first assume that
  $(\Sigma_\h := \Sigma_\h(\hat J, \hat H))_{\h\in\cI'}$ is given
  exactly, without any $\O(\h^2)$ error term.  By
  Theorem~\ref{theo:bs-bis}, there exists a ball $B$ around $c$,
  contained in the set of regular values of $F$, and $\h_0>0$, such
  that $(\Sigma_\h\cap B, \interval[open]0{\h_0}, B)$ is an asymptotic
  lattice.

  Let $\bar k_\h$ be a linear labelling of a neighborhood
  $\tilde B\Subset B$ of $c$ constructed by the algorithm of
  Theorem~\ref{theo:algo}. Let $b\in \tilde B$ and let $\lambda\in\Lh$
  be such that $\lambda = b + \O(\h)$ (for instance, choose a closest
  point to $b$). Let $(\bar n, \bar m) = \bar k_\h(\lambda)$ be the
  corresponding label obtained from the algorithm, \emph{i.e.}
  $\lambda = \lambda_{\bar n, \bar m}$.
  
  Let $k_\h$ be the good labelling associated with $\bar k_\h$ (it is
  not known from the algorithm).  There exist $\h$-dependent integers
  $n_0,m_0$ such that for any $\lambda\in\Sigma_\h$, if we let
  $(n,m) = k_\h(\lambda)$ and $(\bar n, \bar m) = \bar k_\h(\lambda)$,
  then
  \begin{equation}
    \label{equ:label-shift}
    (n,m) = (\bar n, \bar m) + (n_0, m_0).
  \end{equation}
  Recall from~\eqref{equ:quantum-rotation-number-real} that the
  quantum rotation number $ \hat w_\h(n,m)$ is by definition
  \[
    \hat w_\h(n,m) = \frac{E_{n+1,m}(\h) - E_{n,m}(\h)}{E_{n,m+1}(\h)
      - E_{n,m}(\h)}\,,
  \]
  where for all $n,m$, we denote by $E_{n,m}(\h)$ the second component
  of $\lambda = k_\h^{-1}(n,m)\in\Sigma_\h$.  Hence
  \[
    \hat w_\h(n,m) = \frac{\pi_2\lambda_{\bar n + 1,\bar m} -
      \pi_2\lambda_{\bar n,\bar m}}%
    {\pi_2\lambda_{\bar n,\bar m + 1} - \pi_2\lambda_{\bar n,\bar
        m}}\,,
  \]
  where $\pi_2:\RM^2 \to \RM$ is the projection onto the second
  factor; note that this quantity can be computed directly from the
  algorithm. We denote it by $[w_{\h_i}](b)$. It just remains to apply
  Theorem~\ref{theo:q-rotation-number} using an asymptotic chart $G_\h$ 
  associated with the good labelling $k_\h$: for the action variables
  $I=G_0^{-1}(F)$, defined in a fixed neighborhood of the torus
  $\Lambda_c = F^{-1}(c)$, containing $\Lambda_b = F^{-1}(b)$, the
  classical rotation number is
  \[ [w_I(\Lambda_b)] = [\hat w_\h(n,m)] + \O(\h)\,
  \]
  which gives~\eqref{equ:recover-w}, and hence can be recovered in the
  limit as $\h\to 0$, $\h\in\cI'$.

  It remains to investigate the effect of the error term $\O(\h^2)$ in
  the knowledge of the joint spectrum. By~\eqref{equ:unique-ball},
  this will not affect, if $\h$ is small enough, the choice of all
  points of the form $\mu_2$ in Lemma~\ref{lem:unique}, that is, once
  the affine basis $(\lambda_{0,0}, \lambda_{1,0}, \lambda_{0,1})$ is
  chosen (Definition~\ref{defi:affine-basis}). In contrast, the choice
  of this affine basis \emph{can} depend on the error term, because we
  are minimizing distances of order $\h$. Fortunately, this error
  won't affect estimates like~\eqref{equ:norme-reste}, provided we
  accept to make the constant $L_1$ larger, which is harmless. Thus,
  the perturbed triple
  $(\lambda'_{0,0}, \lambda'_{1,0}, \lambda'_{0,1})$ is still an
  affine basis for $\h$ small enough, and the rest of the algorithm
  (Lemma~\ref{lem:unique}) goes through, leading to a labelling
  $\lambda'_{\bar n, \bar m}$ of the perturbed joint spectrum near
  $c$. Let $\lambda_{\bar n,\bar m}$ be the unique point in
  $\Sigma_\h$ that is $\O(\h^2)$-close to $\lambda'_{\bar n,\bar m}$,
  using~\eqref{equ:unique-ball} again. We have
  $\lambda_{\bar n,\bar m} = b + \O(\h)$, and because of that
  uniqueness, $\lambda_{\bar n,\bar m}\mapsto (\bar n,\bar m)$ is a
  linear labelling of $\Sigma_\h$. As a above, we can now introduce
  the good labelling for which~\eqref{equ:label-shift} holds. The
  resulting ``perturbed quantum rotation number'' will be computed as
  \[
    \hat w'_\h(n,m) = \frac{\pi_2\lambda'_{\bar n + 1, \bar m} -
      \pi_2\lambda'_{\bar n, \bar m}}%
    {\pi_2\lambda'_{\bar n, \bar m + 1} - \pi_2\lambda'_{\bar n, \bar
        m}} = \frac{\pi_2\lambda_{\bar n + 1, \bar m}
      -\pi_2\lambda_{\bar n, \bar m}}%
    {\pi_2\lambda_{\bar n, \bar m + 1} - \pi_2\lambda_{\bar n, \bar
        m}} + \O(\h)\,.
  \]
  Hence, as above, we introduce the corresponding asymptotic chart
  $G_\h$ and conclude that
  \[
    [w_{I}(\Lambda_b)] = [\hat w'_\h(n,m)] + \O(\h)\,,
  \]
  with $I = G_0^{-1}\circ F$.
\end{demo}

We may now turn to the global problem which, for simplicity, we state
without the $\O(\h^2)$ perturbation. The notation follows
Proposition~\ref{prop:monodromy}.
\begin{theo} 
  Let $(\hat J, \hat H)$ be a quantum integrable system with principal
  symbols $F:=(J,H)$. Let $\textup{B}_r$ be the set of regular
  Liouville tori of $F$, and let $B_r = F(\textup{B}_r)$; we assume
  that all fibers $F^{-1}(c)$ for $c\in B_r$ are compact and
  connected. Then, from the joint spectrum $\Sigma_\h$ of
  $(\hat J, \hat H)$, one can construct a map
  $\hat \omega: \tilde{\textup{B}}_r\to\RP^1$ such that, if
  $[\tilde w]$ is the globalized rotation number defined in
  Proposition~\ref{prop:monodromy}, then
  \[
    \hat \omega = A\circ [\tilde w] + \O(\h)\,,
  \]
  for some fixed $A\in\textup{SL}(2,\ZM)$, and the remainder is
  locally uniform.  In particular, the globalized rotation number
  $[\tilde w]$ can be recovered from the joint spectrum, modulo a
  global Möbius transformation.
\end{theo}
\begin{demo}
  Because of the connectedness of fibers of $F$, we may identify
  $\textup{B}_r$ with its image $B_r$, and endow the latter with the
  integral affine structure given by action coordinates.  We construct
  the map $\hat\omega$ by a \v Cech cohomology argument. Fix
  $c_0\in {B}_r$, and let $\gamma:[0,1]\to {B}_r$ be a
  path starting at $\gamma(0)=c_0$. Let $c=\gamma(1)$. Let
  $\lambda\in\Sigma_\h$ be a nearest point to $c$, so that, by
  Lemma~\ref{lemm:lattice-dense},
  \[
    \lambda = c + \O(\h)\,.
  \]
  Applying Theorem~\ref{theo:algo} we get a neighborhood $B_0$ of
  $c_0$ equipped with a linear labelling $\bar k_{0,\h}$ of
  $(\Sigma_\h\cap B_0,\cI', B_0)$.  Cover the image $\gamma([0,1])$ by
  a finite union of small balls $B_1,\dots,B_{N}$, such that
  $c\in B_{N}$, $B_i \cap B_{i+1}\neq \varnothing$, for all
  $i=0,\dots, N$, and such that on each $B_i$, the algorithm of
  Theorem~\ref{theo:algo} produces a linear labelling $\bar k_{i,\h}$
  of $(\Sigma_\h\cap B_i, \cI', B_i)$. From
  Proposition~\ref{prop:changement-lineaire} applied to the
  restrictions of $\bar k_{i,\h}$ and $\bar k_{i+1, \h}$ on the
  asymptotic lattice $\Sigma_\h\cap B_i\cap{B_{i+1}}$, there exists a
  unique matrix $A_i\in\textup{SL}(2,\ZM)$ and a family
  $(\varkappa_{i,\h})_{\h\in\cI'}$ in $\ZM^2$ such that, for $\h$
  small enough,
  \begin{equation}
    \label{equ:Ai}
    \bar k_{i+1,\h} = A_i \circ \bar k_{i,\h} + \varkappa_{i,\h} \quad
    \text{ on } \quad \Sigma_\h\cap B_i\cap{B_{i+1}}\,.
  \end{equation}
  We define
  \begin{equation}
    \label{equ:global-quantum}
    \hat\omega(\gamma) = [\hat w_{N}](A_0^{-1} \circ A_1^{-1} \circ
    \cdots A_{N-1}^{-1}\circ \bar k_{N,\h}(\lambda))\,,
  \end{equation}
  where $[\hat w_{N}]$ is the quantum rotation number associated with
  $\bar k_{N,\h}$. It is constructed from the joint spectrum as in
  Theorem~\ref{theo:A}.  Note that the transition matrices $A_i$ can
  be detected from the algorithm by comparing the affine basis
  described in Theorem~\ref{theo:basis}.

  By Lemma~\ref{prop:changement-lineaire}, the sheaf that assigns to a
  point $c\in {B}_r$ a linear labelling (given by the
  algorithm) on a small neighborhood of $c$ has constant transition
  functions, when $\h$ is small enough, modulo the addition of
  $\h$-families in $\ZM^2$. This ensures that the cocyle condition for
  the linear part of the transition functions is satisfied, and hence
  that the definition in~\eqref{equ:global-quantum} is invariant by
  homotopy transformation of the path $\gamma$ with fixed endpoints,
  provided the initial labelling $\bar k_{0,\h}$ is fixed. Thus, it
  defines a map
  \[
    \hat\omega :\tilde{\textup{B}}_r \to \RP^1\,.
  \]
  This map will be modified by a global $\textup{SL}(2,\ZM)$
  transformation if one changes $\bar k_{0,\h}$. \emph{A priori}, the
  smallness of $\h$ for which~\eqref{equ:global-quantum} is defined
  depends on $\gamma$; however, if $\gamma$ stays in a compact region,
  one can use a fixed, finite covering by small balls on which the
  algorithm applies, and hence obtain a uniform $\h_0>0$ for
  which~\eqref{equ:global-quantum} holds for all $\h\leq \h_0$.

  From~\eqref{equ:Ai} we get, using~\eqref{equ:changement-A},
  \begin{equation}
    \label{equ:Ai2}
    A_i = [(G_{i+1,0}'(\xi_{i+1})]^{-1} (G_{i,0})'(\xi_i)\,.
  \end{equation}
  with $\xi_i := G_{i,0}^{-1}(c_i)$.

  Let $\tilde I: \tilde{\textup{B}}_r\to\RM^2$ be the global action
  variable used in Proposition~\ref{prop:monodromy} to define
  $[\tilde w]$. Along the path $\gamma$, the balls $B_i$ can be lifted
  to open sets $\tilde B_i\subset \tilde{\textup{B}}_r$ on which
  $\pi_i$, the restriction of
  $\pi:\tilde{\textup{B}}_r\to \textup{B}_r$ is a diffeomorphism. Since
  both $\tilde I \circ \pi_i^{-1}$ and $I_i:= G_{i,0}^{-1}\circ F$ are
  action variables above $B_i$, we must have
  \[
    \tilde I \circ \pi_i^{-1} = Z_i \circ I_i\,,
  \]
  for some $Z_i\in\textup{SL}(2,\ZM)$, which is of course independent
  of $\h$. Using~\eqref{equ:Ai2}, we obtain
  $A_i = Z_{i+1}^{-1}\circ Z_i$, for all $i=0,\dots, N-1$. Thus,
  $A_{N-1}\cdots A_1 A_0 = Z_N^{-1} Z_0$. For $\h$ small enough,
  $Z_0^{-1} Z_N \bar k_{N,\h}(\lambda)$ stays in the image of the
  linear labelling $\bar k_{N,\h}$; hence we can write
  \begin{align*}
    [\tilde w] (\gamma)
    & =  [\hat w_{N}](Z_0^{-1} Z_N \bar k_{N,\h}(\lambda))\\
    & = [w_{Z_0^{-1} Z_N I_N}(\Lambda)]  + \O(\h)
      \quad \text{ by Theorem~\ref{theo:q-rotation-number}}\\
    & = \trsp Z_0 \circ [w_{Z_N I_N}(\Lambda)]  + \O(\h)
      \quad \text{ by Lemma~\ref{lemm:ItoI'}}\\
    & = \trsp Z_0 \circ [\tilde w](\gamma)  + \O(\h)
  \end{align*}
  because by definition,
  $[\tilde w]\restr_{\tilde B_N} = [w_{\tilde I \circ \pi_N^{-1}}] =
  [w_{Z_N I_N}]$. 
\end{demo}

Finally, let us consider the semitoric case. Combining
Theorem~\ref{theo:q-semitoric} with the semitoric algorithm of
Section~\ref{sec:semitoric-algorithm}
(Proposition~\ref{prop:semitoric-algorithm}), we obtain:
\begin{coro}
  \label{coro:semitoric}
  Let $(\hat J, \hat H)$ be a semitoric quantum integrable system with
  principal symbols $F=(J,H)$. Let $c\in\RM^2$ be regular value of
  $F$. Then from the joint spectrum $(\Sigma_\h)_{\h\in\cI}$ in a
  neighborhood $B$ of $c$, where $\cI$ is a set accumulating at zero,
  one can algorithmically compute the quantum rotation
  number~\eqref{equ:quantum-rotation-number-real}, and this quantum
  rotation number is an $\O(\h)$ deformation of the classical
  semitoric rotation number $w(\Lambda)$, uniformly for all Liouville
  tori $\Lambda\subset F^{-1}(B)$.
\end{coro}

In the presence of an elliptic singularity we may apply
Proposition~\ref{prop:elliptic-labelling}, instead of the general
algorithm, to obtain the corresponding statement.
\begin{coro}
  \label{coro:semitoric-elliptic}
  Let $(\hat J, \hat H)$ be a semitoric quantum integrable system with
  principal symbols $F=(J,H)$. Let $c\in\RM^2$ be a simple
  $J$-transversal elliptic critical value of $F$. Then from the joint
  spectrum $(\Sigma_\h)_{\h\in\cI}$ in a neighborhood $B$ of $c$,
  where $\cI$ is a set accumulating at zero, using the labelling of
  Proposition~\ref{prop:elliptic-labelling}, one can compute the
  quantum rotation number~\eqref{equ:quantum-rotation-number-real},
  and this quantum rotation number is an $\O(\h)$ deformation of the
  classical semitoric rotation number $w(\Lambda)$, uniformly for all
  Liouville tori $\Lambda\subset F^{-1}(B)$.
\end{coro}

In both cases, the detection of the rotation number is robust with
respect to an $\O(\h^2)$ perturbation of the joint spectrum.

\medskip

\begin{rema}
  We would like to stress the difference between recovering the
  rotation number directly from the spectrum, in a constructive way
  (Theorem~\ref{theo:A}), and simply proving the abstract injectivity
  result that the joint spectrum determines the rotation number (which
  means that if two systems have the same joint spectrum, then they
  must share the same rotation number, modulo Möbius
  transformations). The latter, which is of course a consequence of
  the former, is actually quite easier (it doesn't require any
  algorithm); in view of Proposition~\ref{prop:monodromy}, it is a
  direct consequence of Proposition~\ref{prop:A} which implies that
  the integral affine structure of $\textup{B}_r$ is determined by the
  joint spectrum. This fact was already exploited in~\cite{san-mono}
  to prove that classical monodromy determines quantum monodromy.
\end{rema}
\begin{rema}
  It would be quite interesting (and challenging) to extend our
  results to the case where the fibers of $F$ could have several
  connected components. The microlocal derivation of Bohr-Sommerfeld
  rules should continue to hold, giving rise, on the spectral side, to
  a superposition of several asymptotic lattices. In the
  one-dimensional case, it was possible to separate each contribution
  by a Fourier analysis~\cite{san-inverse}; however, in that case,
  because of dimension 1, the labelling problem was obvious; the
  multidimensional case is certainly much more involved and, to our
  knowledge, completely open.
\end{rema}
\begin{rema}
  Quantum integrable systems, as defined in Section~\ref{sec:quantum},
  come from pseudo\--differen\-tial operators, and hence are defined for
  an interval of values of $\h\in\interval[open left]0{\h_0}$. For
  this reason, they are obviously $\h$-continuous, and one can apply
  Proposition~\ref{prop:drift} to detect their drift. However, the
  detection of the rotation number in Theorem~\ref{theo:A} does not
  use the $\h$-continuity. As a consequence, it can in principle be
  applied to quantum integrable systems defined by Berezin-Toeplitz
  operators on compact, prequantizable symplectic manifolds, using the
  Bohr-Sommerfeld theory developed in~\cite{charles-bs}.
\end{rema}

We hope that the formalism of asymptotic lattices, that we have tried
to develop here in a precise way, independently of any particular
quantization scheme, should help attacking other inverse problems for
quantum integrable systems. For instance, our analysis is already used
in~\cite{san-yohann21} to solve the inverse spectral conjecture for
semitoric systems; it may also prove useful for the geometric
quantization of Lagrangian fibrations, extending the first order
Bohr-Sommerfeld quantization that is often used instead. However, in
order to get a more complete picture of asymptotic lattices \emph{vs.}
joint spectra, one should include singularities of the moment map into
the picture.

\bibliographystyle{abbrv}%

\end{document}